\documentclass[12pt]{article}
\usepackage{amsthm,amsmath,amssymb,amscd,amsfonts,mathrsfs}\begin{document}

\bigbreak\centerline{\bf Section Extension from Hyperbolic Geometry of Punctured Disk}
\centerline{\bf and Holomorphic Family of Flat Bundles}

\bigbreak\centerline{\it Dedicated to Fabrizio Cantanese}

\bigbreak\centerline{Yum-Tong Siu\ %
\footnote{Partially supported by Grant 1001416 of the National Science
Foundation.  Written for the Festschrift of the sixtieth birthday of Fabrizio Cantanese} }

\bigbreak\noindent{\bf Abstract.}  The construction of sections of bundles with prescribed jet values plays a fundamental role in problems of algebraic and complex geometry.  When the jet values are prescribed on a positive dimensional subvariety, it is handled by theorems of Ohsawa-Takegoshi type which give extension of line bundle valued square-integrable top-degree holomorphic forms from the fiber at the origin of a family of complex manifolds over the open unit $1$-disk when the curvature of the metric of line bundle is semipositive.  We prove here an extension result when the curvature of the line bundle is only semipositive on each fiber with negativity on the total space assumed bounded from below and the connection of the metric locally bounded, if a square-integrable extension is known to be possible over a double point at the origin.  It is a Hensel-lemma-type result analogous to Artin's application of the generalized implicit function theorem to the theory of obstruction in deformation theory.  The motivation is the need in the abundance conjecture to construct pluricanonical sections from flatly twisted pluricanonical sections.  We also give here a new approach to the original theorem of Ohsawa-Takegoshi by using the hyperbolic geometry of the punctured open unit $1$-disk to reduce the original theorem of Ohsawa-Takegoshi to a simple application of the standard method of constructing holomorphic functions by solving the $\bar\partial$ equation with cut-off functions and additional blowup weight functions.

\bigbreak

\bigbreak\centerline{\bf Table of Contents}
{\footnotesize\begin{itemize}
\item[]Introduction
\item[Part I.] Extension Problem
\begin{itemize}
\item[\S1] Extension from hyperbolic geometry of punctured disk
\item[\S2] A priori estimates for extension and shifting of blowup weight function
\end{itemize}
\item[Part II.] Holomorphic Family of Flat Bundles
\begin{itemize}
\item[\S1] Natural metrics of family of flat line bundles
\item[\S2] Extension with controllably negative curvature
\item[\S3] Local linearity of subvariety defined by minimum number of independent flatly twisted pluricanonical sections
\item[\S4] Two descriptions of moduli, one by projective embedding and one by transition functions
\item[\S5] Algebraicity of subvariety in moduli described by transition functions
\item[\S6] Technique of Gelfond-Schneider, Lang, Bombieri, Brieskorn, and Simpson
\end{itemize}
\item[]References
\end{itemize}}

\bigbreak

\bigbreak

\bigbreak\noindent{\bf Introduction.}  The construction of sections of bundles with presribed jet values plays a fundamental role in problems of algebraic and complex geometry.  When the jet values are prescribed on a positive dimensional subvariety, it is handled by theorems of Ohsawa-Takegoshi type [Ohsawa-Takegoshi1987] which give extension of line bundle valued integrable top-degree holomorphic forms from the fiber at the origin of a family of complex manifolds over the open unit $1$-disk when the curvature of the metric of line bundle is semipositive.

\medbreak We prove here an extension result when the curvature of the line bundle is only semipositive on each fiber with negativity on the total space assumed bounded from below and the connection of the metric locally bounded, if a square-integrable extension is known to be possible over a double point at the origin (Theorem (II.2.10)).  It is a Hensel-lemma-type result [Hensel1897] analogous to Artin's application of the generalized implicit function theorem to the theory of obstruction in deformation theory [Artin1968, Artin1969, Wavrik1975].

\medbreak The motivation stems from one step in the details of the implementation of the analytic methods introduced and sketched in [Siu2010] for the proof of the abundance conjecture.  For the abundance conjecture one needs to construct holomorphic sections of holomorphic line bundles with only semi-positivity of the curvature of the line bundle along each fiber of some holomorphic family of complex manifolds instead of the usual semi-positivity of the curvature current on the total space of the family.  The reason is that in the family of flat line bundles on a compact complex algebraic manifold, though each flat line bundle has zero curvature on the fiber of the family, the line bundle over the total space of the family does not have semipositive curvature (see (II.2.1.1) and (II.2.5)).   This kind of Hensel-lemma-type extension is used to construct pluricanonical sections from flatly twisted pluricanonical sections.  For example, with the use of Gelfond-Schneider type arithmetic arguments (\S6 of Part II) and the two different descriptions of flat line bundles (\S4 of Part II), it yields the following structure result for the subvariety in the moduli space of flat bundles defined by the dimension of the space of flatly twisted pluricanonical sections (Theorem (II.1)).  For a compact complex algebraic manifold $X$ and positive integers $m$ and $q$, the subvariety $Z_{m,q}$ of flat line bundles $F$ such that $\dim_{\mathbb C}H^0\left(X, mK_X+F\right)\geq q$ is regular and is a finite union of translates of abelian subvarieties by torsion elements in the abelian variety of all flat line bundles.

\medbreak This result (Theorem (II.1)) is already in the literature ({\it e.g.,} [Campana-Peternell-Toma2007] and [Budur2009]) when the numerical Kodaira dimension of $X$ is zero, because it can be reduced to the case $m=1$ by taking an $m$-sheeted branched cover of $X$.  However, when the numerical Kodaira dimension of $X$ is greater than $0$, the $m$-canonical section on $X$ obtained from a holomorphic canonical section on the $m$-sheeted branched of $X$ may have poles and may not be holomorphic on $X$.  The pluricanoincal situation in the case of nonzero numerical Kodaira dimension is much more delicate.

\medbreak We would like to remark that in [Campana-Peternell-Toma2007, Remark 3.6] it is mentioned as plausible the expectation that the generalized Green-Lazarsfeld set consisting of all flat line bundles $F$ of $X$ with $\dim_{\mathbb C}H^p\left(X,mK_X+L\right)\geq q$ is a finite union of translates of subtori by torsion elements.  Theorem (II.1) here is the confirmation, for the case of $p=0$, of the expected plausible statement given in [Campana-Peternell-Toma2007, Remark 3.6].

\medbreak Here instead of proving first the more abstract and more general Theorem (II.2.10) and then deriving Theorem (II.1) with the use of Gelfond-Schneider type arithmetic arguments and the two descriptions of flat line bundles,  we first prove a special case of it (Proposition (II.2.7) and Proposition (II.2.8)) which is applicable for derivation of Theorem (II.1) and then present Theorem (II.2.10) as the more general formulation to which the arguments in the proof of Proposition (II.2.7) and Proposition (II.2.8) also apply.  The reason is that Theorem (II.1) is the primary motivation for Theorem (II.2.10) and, what is more important, is that the explicit expression of the curvature in the proof of Proposition (II.2.7) and the arguments necessitated and motivated by such an explicit expression of the curvature make the proof more transparent and more concrete.

\medbreak The core tool in the construction of holomorphic sections of line bundles with prescribed jets is extension results of Ohsawa-Takegoshi type when the given jets are prescribed on positive dimensional subvarieties.  This technique originated with the paper of Ohsawa and Takegoshi in 1987 [Ohsawa-Takegoshi1987].   In [Ohsawa-Takegoshi1987] Ohsawa and Takegoshi, in the context of a bounded Stein domains $\Omega$ in ${\mathbb C}^n\times{\mathbb C}$ (with coordinates $\left(z_1,\cdots,z_n\right)$ of ${\mathbb C}^n$ and $w$ of ${\mathbb C}$ respectively) and a plurisubharmonic function $\varphi$, succeeded in extending any holomorphic function $f$ on $\Omega\cap\left\{w=0\right\}$ with finite $L^2$ norm $\left\|f\right\|_\varphi$ with weight function $e^{-\varphi}$ to a holomorphic function $F$ on $\Omega$ whose $L^2$ norm $\left\|F\right\|_\varphi$ with weight function $e^{-\varphi}$ is bounded by $\left\|f\right\|_\varphi$ times a constant depending on the supremum $\sup_\Omega|w|$ of $|w|$ on $\Omega$.  The significance of the bound depending only on $\sup_\Omega|w|$ and not on the diameter of $\Omega$ is that only the positivity of curvature along the direction of $w$ is required by their method.  For the direction along ${\mathbb C}^n$, where the given section is already holomorphic in, only semi-positivity of the curvature is needed.  For their proofs Ohsawa and Takegoshi used results and strategies from the work of Andreotti-Vesentini [Andreotti-Vesentini1965], Donnelly-Fefferman [Donnelly-Fefferman1983], and Donnelly-Xavier [Donnelly-Xavier1985].

\medbreak Over the years a theory around the theorem of Ohsawa-Takegoshi [Ohsawa-Takegoshi1987] has been developed with various reformulations and generalizations, resulting in a rich collection of extension results (see {\it e.g.,} [Berndtsson1996, Chen2011, Kim2010, Manivel1993, Paun2007, Popovici2005, Siu1996, Straube-Zampieri2011, Takayama2006,  Varolin2008]).  In order to distinguish them from the original theorem of Ohsawa-Takegoshi, we call them extension results of Ohsawa-Takegoshi type.

\medbreak Extension results of Ohsawa-Takegoshi type have been successfully applied to the solution of some longstanding problems in algebraic geometry.  The first application is the original theorem of Ohsawa-Takegoshi itself given in [Angehrn-Siu1995] to prove the semicontinuity of multiplier ideal sheaves to provide for the first time a confirmation of the freeness part of Fujita's conjecture for general dimension but with a weaker bound.   Later in [Siu1998, Siu2002] the conjecture on the deformational invariance of plurigenera was solved by using as one ingredient the following extension result of Ohsawa-Takegoshi type.  For a holomorphic family $\pi:X\to\Delta$ of compact complex algebraic manifolds over the open unit $1$-disk and a holomorphic line bundle $L$ on $X$ with a metric $e^{-\varphi}$ of nonnegative curvature current, any element $f$ of $\Gamma\left(X_0, K_{X_0}+L\right)$ (where $X_0=\pi^{-1}(0)$ with $\left\|f\right\|=\int_{X_0}\left|f\right|^2 e^{-\varphi}$ finite can be extended to an element $F$ of $\Gamma\left(X, K_X+L\right)$ with $\int_X\left|F\right|^2 e^{-\varphi}$ bounded by $\left\|f\right\|$ times a constant which is universal.

\medbreak Because of its usefulness in problems of algebraic and complex geometry and because of the {\it ad hoc} nature of the techniques used in its proof, there have been attempts to understand the original theorem of Ohsawa-Takegoshi [Ohsawa-Takegoshi1987] by making it a special case of the general method of solving $\bar\partial$-equation with some appropriate cut-off functions and additional blowup weight functions.  For example, in [Siu1996] an alternative presentation of the argument for the original theorem of Ohsawa-Takegoshi was given to make it closer to the usual method of using $L^2$ estimates of $\bar\partial$ with some appropriate cut-off functions and additional blowup weight functions.  The alternative proof given in [Siu1996] uses the two-weight-function approach, but it is still quite different from usual method of using $L^2$ estimates of $\bar\partial$ with some appropriate cut-off functions and additional blowup weight functions.

\medbreak In this note we use the hyperbolic geometry of the punctured open unit $1$-disk to construct cut-off functions and additional blowup weight functions to make the original theorem of Ohsawa-Takegoshi a straightforward application of the usual method of using $L^2$ estimates of $\bar\partial$ with some appropriate cut-off functions and additional blowup weight functions (see \S1 of Part I).

\medbreak With hindsight it turns out that the two-weight-function proof given in [Siu1996] can also be naturally interpreted in terms of the hyperbolic geometry of $\Delta-\left\{0\right\}$.  In \S2 of Part I we present this natural interpretation by formulating the {\it a priori} estimate for extension in terms of the current defined by the initial fiber and by introducing the technique of shifting a blowup weight function from the right-hand side of the $\bar\partial$ equation to its solution.

\medbreak The two different descriptions of the moduli of flat line bundles and the Gelfond-Schneider type arithmetic argument used to derive Theorem (II.1) from Proposition (II.2.7) are given in \S4, \S5, and \S6 of Part II. These arguments are needed to handle the conclusion in Theorem (II.1) about torsion elements of the moduli of flat bundles.

\medbreak There are two ways of describing the moduli of flat line bundles, one by modification of a connection by a holomorphic $1$-form and the other by transition functions.  The second description is obtained by integrating (the holomorphic $1$-form in) the first description.  The existence of torsion flat line bundles in the translates of abelian subvarieties in Theorem (II.1) results from the usual arithmetic argument which places a great constraint on $Z_{m,q}$ from its definition over the algebraic closure of the field of rational numbers in both descriptions of the moduli (after using the specialization technique of replacing algebraically independent numbers in the coefficients of the defining polynomials of $X$ by appropriate algebraic numbers to replace $X$ by another $X$ so that the conclusion holds for $X$ if it holds for the new $X$).  The intermediate result of $Z_{m,q}$ being a finite union of translates of abelian varieties is used so that by replacing the full abelian variety of the moduli of flat line bundles by an appropriate abelian subvariety whose intersection with $Z_{m,q}$ is zero-dimensional, we can reduce the general case to the special case where $Z_{m,q}$ consists of a finite number of algebraic points and we can then apply the arithmetic argument the finite number of algebraic points.  This arithmetic argument, of constraints in two different kinds of algebraic moduli related by integration, was first introduced by Gelfond and Schneider [Gelfond1934, Schneider1934] in their independent solution of the seventh problem of Hilbert [Hilbert1900] and was later generalized by Lang [Lang1962, Lang1965, Lang1966] and Bombieri [Bombieri1970, Bombieri-Lang1970] and used for monodromy and torsion results by Brieskorn [Brieskorn1970] and Simpson [Simpson1993].  The arithmetic argument consists of passing to limit of the application of the First Main Theorem of Nevanlinna to a sequence of appropriate functions which are polynomials of functions used in the description which is the integral of the other description.

\medbreak We carry out our proof of Theorem (II.1) in a direct, explicit, elementary setting in analysis to make the presentation clearer and more concrete and to minimize the obscuration from complicated notations and terminology.

\medbreak We would like to remark that Theorem (II.2.10), which concerns extension even without semi-positivity of the curvature current if extension is known over a double point, is analogous to bounding the dimension of the first cohomology as the obstruction to solving the $\bar\partial$ equation in terms of the extent of failure of positivity of the curvature of the line bundle involved, for example, as quantitatively formulated in the form of holomorphic Morse inequality in the work of Demailly [Demailly1985] (see (II.2.9) below).  Here we take advantage of the special form of our metric with non semipositive curvature to limit the obstruction of extension to the extendibility over the double point.  More useful and significant would be the development of a general theory linking quantitatively the extent of failure of semipositivity of curvature to the order of the multiple over which extension needs to be assumed to guarantee no obstruction to local extension.

\medbreak\noindent{\it Notations.}  The structure sheaf of a complex space $X$ is denoted by ${\mathcal O}_X$.  The full ideal sheaf of a subvariety $Y$ of $X$ is denoted by ${\mathcal I}_Y$.  When $Y$ is contained in another subvariety $Z$ of $X$, whether ${\mathcal I}_Y$ is the sheaf on $Z$ or on $X$ will be specified if it is not clear from the context.  The space $H^0\left(X, V\right)$ of all holomorphic sections of a bundle $V$ over a complex space $X$ is also denoted by $\Gamma\left(X, V\right)$.

\medbreak ${\mathbb C}$ denotes the set of all complex numbers.  ${\mathbb R}$ denotes the set of all real numbers. ${\mathbb N}$ denotes the set of all positive integers.  ${\mathbb Z}$ denotes the set of all integers.  ${\mathbb Q}$ denotes the set of all rational numbers.
$\overline{\mathbb Q}$ denotes the algebraic closure of ${\mathbb Q}$.  $\Delta$ denotes the open unit $1$-disk and $\Delta_R$ denotes the open disk in ${\mathbb C}$ with center $0$ and radius $R$.
The punctured open unit $1$-disk $\Delta-\left\{0\right\}$ is also denoted by $\Delta^*$.  The punctured complex line ${\mathbb C}-\left\{0\right\}$ is also denoted by ${\mathbb C}^*$.

\medbreak Two conventions for the curvature of a metric $e^{-\varphi}$ of a line bundle will be used.  One is $\sqrt{-1}\partial\bar\partial\varphi$ and the other is $\frac{\sqrt{-1}}{2\pi}\partial\bar\partial\varphi$.  Which one is being used will be either clear from the context or explicitly specified.  The reason is that when comparing two curvatures, it does not matter which convention is being used and it is less cumbersome without the factor $2\pi$ in the denominator.  However, when the curvature as a closed positive current is compared to the closed positive current of integration on a complex submanifold, the correct normalizing constant with the factor $2\pi$ in the denominator needs to be used.

\bigbreak

\bigbreak\noindent{\bf Part I. Extension Problem}.

\bigbreak For the extension problem we first present the approach from the hyperbolic geometry of the punctured disk which motivates the use of special cut-off functions and additional weight functions, making the extension theorem of Ohsawa-Takegoshi a direct application of the usual $L^2$ estimates of $\bar\partial$.  Then we explain how the two-weight-function proof given in [Siu1996] can be naturally interpreted in terms of the hyperbolic geometry of the punctured disk.

\bigbreak

\bigbreak\noindent{\bf \S1. Extension from hyperbolic geometry of punctured disk.}

\bigbreak\noindent(I.1.1) {\it Standard Preliminaries for $\bar\partial$ Estimates.}  The standard general preliminaries for the $L^2$ estimates of $\bar\partial$ is that for a pseudoconvex domain $\Omega$ spread over a complex Euclidean space and plurisubharmonic functions $\varphi$ and $\psi$, the inequality
$$
\left\|\bar\partial g\right\|_{\varphi+\psi}^2+\left\|\bar\partial^* g\right\|_{\varphi+\psi}^2\geq
\left(\Theta_\psi\cdot g,\,g\right)_{\varphi+\psi}\leqno{({\rm I}.1.1.1)}
$$
holds for $(0,1)$-forms $g$ which is $L^2$ with respect to the weight function $e^{-\varphi-\psi}$ and the Euclidean metric and is in the domains of $\bar\partial$ and $\bar\partial^*$,
where $\Theta_\psi=\sqrt{-1}\partial\bar\partial\psi$ is (up to a universal normalizing constant) the curvature form of the metric $e^{-\psi}$ of the trivial line bundle; $\left\|\cdot\right\|_{\varphi+\psi}$ means the $L^2$ norm with respect to the weight function $e^{-\varphi-\psi}$ and the Euclidean metric; $\Theta_\psi\cdot g$ means contraction of $\Theta_\psi$ and $g$ with respect to the Euclidean metric, and $\left(\cdot,\,\cdot\right)_{\varphi+\psi}$ means the $L^2$ inner product with respect to the weight function $e^{-\varphi-\psi}$ and the Euclidean metric.  Note that here the only use of the curvature $\Theta_\varphi=\sqrt{-1}\partial\bar\partial\varphi$ (up to a universal normalizing constant) of the metric $e^{-\varphi}$ is $\Theta_\varphi\geq 0$.

\medbreak Let $\left(\Theta_\psi\right)^{-1}$ be the $(1,1)$-form whose representative Hermitian matrix is the inverse of that of the $(1,1)$-form $\Theta_\psi$.  To solve the $\bar\partial$-equation $\bar\partial u=v$ for a function $u$ on $\Omega$ subject to the compatibility condition $\bar\partial v=0$ and to use the approach of weak solutions and the Riesz representation theorem to obtain a solution $u$ with the estimate $\left\|u\right\|_{\varphi+\psi}\leq C$, one takes the inner product of the equation $\bar\partial u=v$ with the test $(0,1)$-form $g$ in the kernel of $\bar\partial$ to get $\left(g,\bar\partial u\right)_{\varphi+\psi}=\left(g,v\right)_{\varphi+\psi}$ and uses $\left(g,\bar\partial u\right)_{\varphi+\psi}=\left(\bar\partial^*g,u\right)_{\varphi+\psi}$ to reduce the problem to the boundedness of the functional $\bar\partial^*g\mapsto\left(g,v\right)_{\varphi+\psi}$ and then to the estimate
$$\left|\left(g,v\right)_{\varphi+\psi}\right|\leq C\left\|\bar\partial^*g\right\|_{\varphi+\psi}=
C\left(\left\|\bar\partial^*g\right\|_{\varphi+\psi}^2+\left\|\bar\partial g\right\|_{\varphi+\psi}^2\right)^{\frac{1}{2}}.$$
From the Cauchy-Schwarz inequality
$$\left|\left(g,v\right)_{\varphi+\psi}\right|_{\varphi+\psi}\leq
\left(\Theta_\psi\cdot g,\,g\right)_{\varphi+\psi}^{\frac{1}{2}}
\left(\left(\Theta_\psi\right)^{-1}\cdot v,\,v\right)_{\varphi+\psi}^{\frac{1}{2}}
$$
and the inequality (I.1.1.1) it follows that
$$\left|\left(g,v\right)_{\varphi+\psi}\right|_{\varphi+\psi}\leq C\left(\left\|\bar\partial^*g\right\|_{\varphi+\psi}^2+\left\|\bar\partial g\right\|_{\varphi+\psi}^2\right)^{\frac{1}{2}}$$
with
$$
C=\left(\left(\Theta_\psi\right)^{-1}\cdot v,\,v\right)_{\varphi+\psi}^{\frac{1}{2}}
$$
and the $\bar\partial u=v$ equation can be solved with the estimate
$$
\left\|u\right\|_{\varphi+\psi}\leq\left(\left(\Theta_\psi\right)^{-1}\cdot v,\,v\right)_{\varphi+\psi}^{\frac{1}{2}}.\leqno{({\rm I}.1.1.2)}
$$

\bigbreak We consider the following context of the original theorem of Ohsawa-Takegoshi [Ohsawa-Takegoshi1987].  

\bigbreak\noindent(I.1.1.3) {\it Theorem of Ohsawa-Takegoshi.} Let the coordinates of ${\mathbb C}^n\times{\mathbb C}$ be $z=\left(z_1,\cdots,z_n\right)\in{\mathbb C}$ and $w\in{\mathbb C}$ and $\Omega$ be a Stein domain in ${\mathbb C}^n\times\Delta$ and $\varphi$ be a plurisubharmonic function on $\Omega$.  Then there exists some positive constant $C$ such that any holomorphic function $f$ on $\Omega\cup\left\{w=0\right\}$ with $\int_{\Omega_0}\left|f\right|^2e^{-\varphi}$ finite can be extended to a holomorphic function $F$ on $\Omega$ with $\int_\Omega\left|F\right|^2e^{-\varphi}\leq C\int_{\Omega\cup\left\{w=0\right\}}\left|f\right|^2e^{-\varphi}$.

\bigbreak\noindent(I.1.2) {\it Cut-Off and Weight Function from Hyperbolic Geometry.}  Here we first show how to use the hyperbolic geometry of the puncture open unit $1$-disk to construct cut-off functions and additional blowup weight functions so that the original theorem of Ohsawa-Takegoshi is simply a straightforward application of the usual method of $L^2$ estimate of $\bar\partial$.  Then we explain how the the two-weight-function approach presented in [Siu1996] is related to the approach of the hyperbolic geometry of the puncture open unit $1$-disk.  Let us first write down the cut-off functions and the additional blowup weight functions.  Let $0<r_1<r_2<1$. We start out with a smooth function $0\leq\Lambda(x)\leq 1$ for $x\in{\mathbb R}$ such that
 $\Lambda\left(\log\log\frac{1}{|\zeta|^2}\right)=1$ for $0<|\zeta|<r_1$ and $\Lambda\left(\log\log\frac{1}{|\zeta|^2}\right)=0$ for $r_2<|\zeta|<1$.  Let
$$\Lambda_1(w)=\Lambda\left(\log\log\frac{1}{|w|^2}\right),\quad
\Lambda_m(w)=\Lambda_1\left(\left|w\right|^{\frac{1}{m}}\right)\quad{\rm for\ }m\in{\mathbb N}
$$
so that $\bar\partial\Lambda_m(w)$ is supported on $\left(r_1\right)^m\leq|w|\leq\left(r_2\right)^m$.  For $m\in{\mathbb N}$ we use the cut-off function $\Lambda_m(w)$
and the additional blowup weight function
$$
\frac{\log\left(1+\left|w\right|^{\frac{2}{m}}\right)}{|w|^2}.
$$
We will explain later how such choices come from the hyperbolic geometry of the punctured disk.  Of course one can go through the arguments of analysis without paying any attention to how such choices come from the hyperbolic geometry of the punctured disk.

\bigbreak\noindent(I.1.3) {\it Power Map and Invariance under Cyclic Group Action.}  In verifying the required $L^2$ estimates of $\bar\partial$ involving both the above cut-off function and additional blowup weight function the easiest way is to use the change of variables $\zeta\mapsto w=\zeta^m$ to go the space of the $\zeta$ variable to get rid of the exponent $\frac{1}{m}$ do the estimates.  The solution of $\bar\partial$ equation in the space of the $\zeta$ variable differs from solution of $\bar\partial$ equation  in the space of the $w$ variable because of the action of the cyclic group of order $m$ defined by multiplying $\zeta$ by an $m$-th root of unity.  If we only consider the solutions of the $\bar\partial$ equation in the space of the $\zeta$ variable which are appropriately invariant with respect to the action cyclic group of order $m$ so that they correspond are the pull-backs of the solutions of the $\bar\partial$ equation in the space of the $w$ variable, then instead of going back and forth between the space of the $\zeta$ variable and the space of $w$ variable when the $\bar\partial$ equation is solved in the space of the $w$ variable and the estimates are done in the space of $\zeta$ variable, we can do both the solution of the $\bar\partial$ equation and also the estimates in the space of the $\zeta$ variable by imposing the condition of invariance in solving the $\bar\partial$ equation in the $\zeta$ variable.  We will take this option of working in the space of the $\zeta$ variable under the condition of invariance.  This is only a matter of expediency and does not change the arguments in any way.   Before we solve the $\bar\partial$ equation, we would like to first track carefully how the
integrals change in the change of variables $\zeta\mapsto w=\zeta^m$, especially the role played by the constant factor $m$.  As we will explain below, keeping track of the role played by the constant factor $m$ corresponds to the technique of the twisted $\bar\partial$ operator in the two-weight-function approach discussed in [Siu1996].
 
\bigbreak\noindent(I.1.4) {\it Transformation of Integral Under Cyclic Cover Map.}  Let $U(w)$ be a function on $\Delta$.  We are interested in the transformation of the integral
$$
\int_\Omega\left|U(w)\right|^2\frac{|dw|^2}{|w|^2}
$$
under the map $\pi_m:\zeta\mapsto w=\zeta^m$.  Let $w=re^{i\theta}$ and $\zeta=\rho e^{i\phi}$.  The volume form
$$
\frac{|dw|^2}{|w|^2}
$$
actually means
$$
\frac{\sqrt{-1}\,dw\wedge d\bar w}{|w|^2}
$$
which is equal to
$$
\frac{2rdr d\theta}{r^2},
$$
because
$$
\sqrt{-1}\,d\left(x+\sqrt{-1}\,y\right)\wedge d\overline{\left(x+\sqrt{-1}\,y\right)}=2dx\wedge dy.
$$
The volume form
$$
\frac{2rdr d\theta}{r^2}=2d\left(\log r\right)d\theta
$$
under $r=\rho^m$ and $\theta=m\phi$ becomes
$$
2d\left(m\log\rho\right)\left(md\phi\right)
=m^2\frac{2\rho
d\rho d\phi}{\rho^2}=
m^2\frac{|d\zeta|^2}{|\zeta|^2}.
$$
The domain whose restriction in the variable $\phi$ is given by $0\leq\phi\leq 2\pi$ corresponds to $m$ times the domain whose restriction in the variable $\theta$ given by $0\leq\theta\leq 2\pi$. So, if we compare the integral over $0\leq\phi\leq 2\pi$ and the integral over $0\leq\theta\leq 2\pi$ (with the other variables suppressed), we should divide the integral over $0\leq\phi\leq 2\pi$ by $m$.  Thus
$$
\int_\Omega\left|U(w)\right|^2\frac{|dw|^2}{|w|^2}=m
\int_\Omega\left|U(\zeta^m)\right|^2\frac{|d\zeta|^2}{|\zeta|^2}.
$$
The key point is that the factor in front of the integral on the right-hand side is $m$ instead of $m^2$ even though
$$
\frac{|dw|^2}{|w|^2}=m^2\frac{|d\zeta|^2}{|\zeta|^2}.
$$

\bigbreak\noindent{\it Computation of Derivative of Cut-Off Function.}  Since
$$\Lambda_1(w)=\Lambda\left(\log\log\frac{1}{|w|^2}\right),\quad
\Lambda_m(w)=\Lambda_1\left(\left|w\right|^{\frac{1}{m}}\right)\quad{\rm for\ }m\in{\mathbb N}
$$
(where  $0\leq\Lambda(x)\leq 1$ is smooth for $x\in{\mathbb R}$ with 
$\Lambda\left(\log\log\frac{1}{|\zeta|^2}\right)=1$ for $0<|\zeta|<r_1<1$ and $\Lambda\left(\log\log\frac{1}{|\zeta|^2}\right)=0$ for $0<r_2<|\zeta|<1$, it follows that
$$
\bar\partial\left(\Lambda_m(w)\right)=\Lambda^\prime\left(\log\log\frac{1}{\left|w\right|^{\frac{2}{m}}}\right)
\frac{d\bar w}{w\log\frac{1}{|w|^2}}
$$
and
$$
\left|\bar\partial\left(\Lambda_m(w)\right)\right|^2
=\left|\Lambda^\prime\left(\log\log\frac{1}{\left|w\right|^{\frac{2}{m}}}\right)\right|^2\,
\frac{\left|d\bar w\right|^2}{|w|^2\left(\log\frac{1}{|w|^2}\right)^2}.
$$
The function $\bar\partial\Lambda_m(w)$ is supported on $\left(r_1\right)^m\leq|w|\leq\left(r_2\right)^m$ and $\Lambda_m(w)=1$ for $0<|w|<\left(r_1\right)^m$ and $\Lambda_m(w)=0$ for $\left(r_2\right)^m<|w|<1$.

\bigbreak\noindent(I.1.5) {\it Use of Constant Factors to Shift Blow-Up from Norm of Right-Hand Side to Norm of Solution.}   Now we do the standard process of exhausting the Stein domain by relative compact Stein subdomains so that we can apply the cut-off function to the pullback by holomorphic retraction of the given holomorphic function $f$ and solve the $\bar\partial$ equation with a blowup weight function.  For some subsequence $m_\nu$ in ${\mathbb N}$ indexed by $\nu\in{\mathbb N}$ there are relatively compact Stein subdomains $\Omega_\nu$ of $\Omega$ for $\nu\in{\mathbb N}$ such that (i) $\Omega_\nu$ is relatively compact in $\Omega_{\nu+1}$; (ii) $\cup_{\nu\in N}\Omega_\nu=\Omega$; and (iii) there is a retraction $\sigma_{\nu+1}$ from some open neighborhood $U_{\nu+1}$ of $\Omega_{\nu+1}$ in $\Omega$ to $\Omega_{\nu+1}\cap\left\{w=0\right\}$ with $\Omega_\nu\cap\left\{|w|<\left(t_2\right)^{m_\nu}\right\}$ contained in $U_{\nu+1}$.  Let $\hat f_{\nu+1}$ be the holomorphic function on $U_{\nu+1}$ which is the pullback of $f|_{\Omega_{\nu+1}\cap\left\{w=0\right\}}$ through $\sigma_{\nu+1}$.  Let $f_\nu$ be the holomorphic function on $\Omega_\nu\cap\left\{|w|<\left(t_2\right)^{m_\nu}\right\}$ which is the restriction of $\hat f_\nu$ to $\Omega_\nu\cap\left\{|w|<\left(t_2\right)^{m_\nu}\right\}$.

\medbreak By replacing $\Omega$ by $\Omega_\nu$ and $m$ by $m_\nu$ we can assume without loss of generality that our given holomorphic function $f$ on $\Omega\cap\left\{w=0\right\}$ is the restriction of a holomorphic function on $\Omega\cap\left\{|w|<\left(t_2\right)^{m_\nu}\right\}$ which is independent of the coordinate $w$ and which we also denote by $f$.  All we need to do is to obtain a holomorphic function $F$ on $\Omega$ whose $L^2$ norm on $\Omega$ with respect to the weight function $e^{-\varphi}$ is dominated by the $L^2$ norm of $f$ on $\Omega\cap\left\{w=0\right\}$ with respect to the weight function $e^{-\varphi}$ times a positive constant which is independent of $\Omega$.

\medbreak Denote by $\pi_m$ the map $\zeta\mapsto w=\zeta^m$ from ${\mathbb C}^n\times{\mathbb C}$ to itself.
Consider on the Stein domain $\left(p_m\right)^{-1}\left(\Omega\right)$ in the domain space of $\pi_m$ the $\bar\partial$-closed $(0,1)$-form
$$
\bar\partial\left(\pi_m^*\Lambda_m\right)\left(\pi_m^*f\right)
$$
as the right-hand side of the $\bar\partial$-equation on $\left(p_m\right)^{-1}\left(\Omega\right)$ and use the weight function
$$
\frac{e^{-\pi_m^*\varphi}}{|\zeta|^2\left(1+|\zeta|^2\right)^2}
$$
to solve the $\bar\partial$ equation with invariance with respect to the action of the cyclic group of order $m$ defined by multiplying $\zeta$ by an $m$-th root of unity.  
The curvature of the weight function
$$
\frac{1}{\left(1+|\zeta|^2\right)^2}
$$
is given by
$$
\partial_\zeta\partial_{\bar\zeta}\left(-\log\frac{1}{\left(1+|\zeta|^2\right)^2}\right)=\frac{2}{\left(1+|\zeta|^2\right)^2}
$$
whose reciprocal is $\left(\frac{2}{\left(1+|\zeta|^2\right)^2}\right)^{-1}$.  So we should consider
the supremum of
$$
\frac{1}{|\zeta|^2\left(1+|\zeta|^2\right)^2}\cdot\left(\frac{2}{\left(1+|\zeta|^2\right)^2}\right)^{-1}=
\frac{1}{2|\zeta|^2}
$$
on $r_1\leq|\zeta|\leq r_2$ which is
$$
\frac{1}{2\left(r_1\right)^2}.
$$
Here the reciprocal $\left(\frac{2}{\left(1+|\zeta|^2\right)^2}\right)^{-1}$ of the curvature $
\partial_\zeta\partial_{\bar\zeta}\left(-\log\frac{1}{\left(1+|\zeta|^2\right)^2}\right)$ is used because of (I.1.1.2).
Let
$$
A_{r_1,r_2}=\int_{r_1\leq|\zeta|\leq r_2}\left|\Lambda^\prime\left(\log\log\frac{1}{|\zeta|^2}\right)\right|^2
\frac{|d\zeta|^2}{|\zeta|^2\left(\log\frac{1}{|\zeta|^2}\right)^2}.
$$
Note that for any $\eta>1$ we can choose $\Lambda$ such that
$$
\displaylines{\left|\Lambda^\prime\left(\log\log\frac{1}{|\zeta|^2}\right)\right|\leq
\frac{\eta}{\log\log\frac{1}{\left|r_1\right|^2}-\log\log\frac{1}{\left|r_2\right|^2}}\cr
\leq
\frac{\eta}{\log\log\frac{\left|r_2\right|^2}{\left|r_1\right|^2}}\cr}
$$
so that $\Lambda\left(\log\log\frac{1}{|\zeta|^2}\right)=1$ for $0<|\zeta|<r_1$ and $\Lambda\left(\log\log\frac{1}{|\zeta|^2}\right)=0$ for $r_2<|\zeta|<1$ and the function $\bar\partial\Lambda_m(w)$ is supported on $\left(r_1\right)^m\leq|w|\leq\left(r_2\right)^m$.
Then
$$
\displaylines{\int_{\left(p_m\right)^{-1}\left(\Omega\right)}\left|\bar\partial\left(\pi_m^*\Lambda_m\right)\left(\pi_m^*f\right)\right|^2
\frac{e^{-\pi_m^*\varphi}}{|\zeta|^2\left(1+|\zeta|^2\right)^2}\left(\partial_\zeta\partial_{\bar\zeta}\left(-\log\frac{1}{1+|\zeta|^2}\right)\right)^{-1}\cr
\leq\frac{A_{r_1,r_2}}{2m\,\left(r_1\right)^2}
\int_{\Omega\cap\left\{w=0\right\}}\left|f\right|^2e^{-\varphi}.}
$$
Now we solve the $\bar\partial$-equation $\bar\partial\hat u=\pi_m^*\bar\partial\left(\Lambda_m f\right)$ on $\pi_m^{-1}(\Omega)$ with the weight function
$$
\frac{e^{-\pi_m^*\varphi}}{|\zeta|^2\left(1+|\zeta|^2\right)^2}.
$$
for the unknown $\hat u$ with invariance under the action of the cyclic group of order $m$ defined by multiplying $\zeta$ by an $m$-root of unity.  We have the $L^2$ estimate
$$
\int_{\pi_m^{-1}(\Omega)}\left|\hat u\right|^2\frac{e^{-\pi_m^*\varphi}}{|\zeta|^2\left(1+|\zeta|^2\right)^2}|d\zeta|^2\leq
\frac{A_{r_1,r_2}}{2m\left(r_1\right)^2}
\int_{\Omega\cap\left\{w=0\right\}}\left|f\right|^2e^{-\varphi}.
$$
The function $\hat u$ is holomorphic in an open neighborhood of $\pi_m^{-1}\left(\Omega\cap\left\{w=0\right\}\right)$ in $\pi_m^{-1}\Omega$ and vanishes on $\pi_m^{-1}\left(\Omega\cap\left\{w=0\right\}\right)$.  The invariance of $\hat u$ means that there exists a function $u$ on $\Omega$ such that $\hat u$ is the pullback of $u$ by the map $\pi_m^{-1}$.  By the observation at the very beginning concerning the transformation of integrals under $w=\zeta^m$, we have
$$
\int_{\pi_m^{-1}(\Omega)}\left|\hat u\right|^2\frac{e^{-\pi_m^*\varphi}}{|\zeta|^2\left(1+|\zeta|^2\right)^2}|d\zeta|^2=
\frac{1}{m}\,\int_\Omega\left|u\right|^2\frac{e^{-\varphi}}{|w|^{\frac{2}{m}}\left(1+|w|^{\frac{2}{m}}\right)^2}|dw|^2.
$$
The key point is the factor $\frac{1}{m}$ in front of the integral on the right-hand side.  As we will see later in (I.2.2) and (I.2.4) this factor $\frac{1}{m}$ corresponds to the shifting of the blowup weight function in the two-weight-function approach, with the big difference that what is being shifted here is just a constant instead of a blowup weight function in (I.2.2) and (I.2.4).
Thus we have the solution $u$ of the $\bar\partial$-equation $\bar\partial u=\bar\partial\Lambda_m f$ on $\Omega$ and
$$
\int_\Omega\left|u\right|^2\frac{e^{-\varphi}}{|w|^{\frac{2}{m}}\left(1+|w|^{\frac{2}{m}}\right)^2}|dw|^2
\leq \frac{A_{r_1,r_2}}{2\left(r_1\right)^2}
\int_{\Omega\cap\left\{w=0\right\}}\left|f\right|^2e^{-\varphi}.
$$
Thus we can set $F=\Lambda_m f-u$ in the limiting situation and get
$$
\int_\Omega\left|F\right|^2 e^{-\varphi}|dw|^2
\leq \frac{2\,A_{r_1,r_2}}{\left(r_1\right)^2}
\int_{\Omega\cap\left\{w=0\right\}}\left|f\right|^2e^{-\varphi},
$$
because the limit of
$$
\frac{1}{|w|^{\frac{2}{m}}\left(1+|w|^{\frac{2}{m}}\right)^2}
$$
as $m\to\infty$ becomes $\frac{1}{4}$ for $w\not=0$.  The constant $C$ in
$$
\int_\Omega\left|F\right|^2 e^{-\varphi}|dw|^2
\leq C
\int_{\Omega\cap\left\{w=0\right\}}\left|f\right|^2e^{-\varphi}
$$
can be taken to be the infimum, over all $0<r_1<r_2<1$, of
$$
\frac{\pi}{\left(r_1\right)^2\left(\log\log\frac{\left|r_2\right|^2}{\left|r_1\right|^2}\right)}
\left(\frac{1}{\log\frac{1}{r_2}}-\frac{1}{\log\frac{1}{r_1}}\right),
$$
because
$$
\displaylines{\int_{r_1<|\zeta|<r_2}\frac{|d\zeta|^2}{|\zeta|^2\left(\log\frac{1}{|\zeta|^2}\right)^2}=
2\pi\int_{r=r_1}^{r_2}\frac{rdr}{r^2\left(\log\frac{1}{r^2}\right)^2}\cr
=\frac{\pi}{2}\int_{r=r_1}^{r_2}\frac{dr}{r\left(\log\frac{1}{r}\right)^2}=\frac{\pi}{2}\left[\frac{1}{\log\frac{1}{r}}\right]_{r=r_1}^{r_2}
=\frac{1}{\log\frac{1}{r_2}}-\frac{1}{\log\frac{1}{r_1}}.}
$$

\bigbreak\noindent(I.1.6) {\it Direct Argument Without Using Cyclic Branched Cover.}  As mentioned above, the introduction of the cyclic branched cover $\zeta\mapsto w=\zeta^m$ is just a matter of expediency.  One can also directly argue without using the cyclic branched cover $\pi_m:\zeta\mapsto w=\zeta^m$ by using the weight function
$$
\frac{\log\left(1+\left|w\right|^{\frac{2}{m}}\right)}{|w|^2}
$$
and the cut-off function
$$\Lambda_1(w)=\Lambda\left(\log\log\frac{1}{|w|^2}\right),\quad
\Lambda_m(w)=\Lambda_1\left(\left|w\right|^{\frac{1}{m}}\right)\quad{\rm for\ }m\in{\mathbb N},
$$
because the measure used for the norm of $\hat u$ is
$$
\frac{\log\left(1+\left|\zeta\right|^2\right)|d\zeta|^2}{|\zeta|^2}
$$
which up to a constant depending on $m$ is equal to
$$
\frac{\log\left(1+\left|w\right|^{\frac{2}{m}}\right)|dw|^2}{|w|^2}.
$$
We can now summarize the result together with the explicit constant in the estimate as follows.

\bigbreak\noindent(I.1.7) {\it Theorem.}  Let the coordinates of ${\mathbb C}^n\times{\mathbb C}$ be $z=\left(z_1,\cdots,z_n\right)\in{\mathbb C}$ and $w\in{\mathbb C}$ and $\Omega$ be a Stein domain in ${\mathbb C}^n\times\Delta$ and $\varphi$ be a plurisubharmonic function on $\Omega$.  Then there exists some positive constant $C$ such that any holomorphic function $f$ on $\Omega\cup\left\{w=0\right\}$ with $\int_{\Omega_0}\left|f\right|^2e^{-\varphi}$ finite can be extended to a holomorphic function $F$ on $\Omega$ with $\int_\Omega\left|F\right|^2e^{-\varphi}\leq C\int_{\Omega\cup\left\{w=0\right\}}\left|f\right|^2e^{-\varphi}$, where the positive constant 
$C$ can be taken to be the infimum, over all $0<r_1<r_2<1$, of
$$
\frac{\pi}{\left(r_1\right)^2\left(\log\log\frac{\left|r_2\right|^2}{\left|r_1\right|^2}\right)}
\left(\frac{1}{\log\frac{1}{r_2}}-\frac{1}{\log\frac{1}{r_1}}\right).
$$

\bigbreak\noindent(I.1.8) {\it Hyperbolic Geometry of Punctured Disk.}  The hyperbolic metric for the punctured disk $\Delta^*=\Delta-\left\{0\right\}$ with coordinate $w$ is
$$
\frac{|dw|^2}{|w|^2\left(\log|w|^2\right)^2}=\frac{|dw|^2}{|w|^2\left(\log\frac{1}{|w|^2}\right)^2}.
$$
The hyperbolic metric for the puncture disk plays a role in the choice of our cut-off function in two ways.  For the radial direction $r=|w|$ the metric is given by
$$
\frac{dr}{r\log\frac{1}{r}}=d\log\log\frac{1}{r}.
$$
This form it is invariant under $s\to r=s^m$ for any positive number $m$, because
$$
\displaylines{d\log\log\frac{1}{r}=d\log\log\frac{1}{s^m}=d\log\left(m\log\frac{1}{s}\right)\cr=d\left(\log m+\log\log\frac{1}{s}\right)=d\log\log\frac{1}{s}.}
$$
The first role played by the hyperbolic metric for the punctured disk in the choice of the cut-off function is that this invariance of the metric in the radial direction under the power map $r\to s=r^m$ suggests setting the relation between the cut-off function $\Lambda_m(w)$ and $\Lambda_1(w)$ to be
$\Lambda_m(w)=\Lambda\left(|w|^{\frac{1}{m}}\right)$ for any positive integer $m$ so that the {\it hyperbolic} radial distance between the inner circle $|w|=\left(t_1\right)^m$ and the outer circle $|w|=\left(t_2\right)^m$ of the annulus $\left(t_1\right)^m\leq|w|\leq\left(t_2\right)^m$ is independent of $m$.

\medbreak The second role played by the hyperbolic metric for the punctured disk in the choice of the cut-off function is that the differential of the function $\log\log\frac{1}{r}$ is the hyperbolic radial distance function $\frac{dr}{r\log\frac{1}{r}}$.  That is the reason why we choose a smooth function $\Lambda_1(w)=\Lambda\left(\log\log\frac{1}{|w|^2}\right)$, where $0\leq\Lambda(x)\leq 1$ is smooth for $x\in{\mathbb R}$ with
$\Lambda\left(\log\log\frac{1}{|\zeta|^2}\right)=1$ for $0<|\zeta|<r_1<1$ and $\Lambda\left(\log\log\frac{1}{|\zeta|^2}\right)=0$ for $0<r_2<|\zeta|<1$.

\medbreak In the choice of the additional blowup weight function
$$
\frac{\log\left(1+\left|w\right|^{\frac{2}{m}}\right)}{|w|^2}
$$
the power map $\zeta\mapsto w=\zeta^m$ mentioned above in the form $s\to r=s^m$ is the reason for using the function $\left|w\right|^{\frac{2}{m}}$.

\bigbreak\noindent(I.1.9) {\it Comparison with Standard Cut-Off and Blowup Weight Functions.}  Finally we would like to compare our choice of cut-off functions and additional weight functions from the hyperbolic geometry of the punctured disk with the standard ones used in the $L^2$ estimates of $\bar\partial$ which comes from the Euclidean geometry of the open $1$-disk.  The standard sequence of additional weight functions is chosen to be $\log\frac{1}{|w|^2+\varepsilon}$ for $\varepsilon>0$ (approaching $0$) and the sequence of cut-off functions is chosen to be $\chi\left(\frac{|w|^2}{\varepsilon^2}\right)$ for $\varepsilon>0$ (approaching $0$), where $0\leq\chi(t)\leq 1$ is a smooth function of a single real variable for $0\leq t<1$ which is identically $1$ in a neighborhood of $0$ and has compact support in $[0,1)$.  The parameter $\varepsilon$ corresponds to our parameter $m\in{\mathbb N}$ and the function $\chi$ plays the role of $\Lambda$.  We use different symbols in both viewpoints to avoid confusion.  

\medbreak The big difference between two approaches lies in the kind of rescaling used for the sequence of cut-functions.  For the usual approach from the Euclidean geometry the rescaling is done by multiplication of the coordinate by a constant $w\mapsto\frac{w}{\varepsilon}$ whereas for our approach from the hyperbolic geometry of the punctured disk the rescaling is done by the power map $\zeta\mapsto w=\zeta^m$. The advantage of the viewpoint of hyperbolic geometry is that the difficulty, in the Euclidean geometric approach, of the blow-up orders inside the integrals is now transformed, in the hyperbolic geometric approach, to a factor which is a power of $m$ and then can be shifted easily from one side of the estimate to the other to guarantee that the right-hand side of the estimate is finite after normalization by an appropriate factor of a power of $m$.

\bigbreak

\bigbreak\noindent{\bf \S2. A priori estimates for extension and shifting of blowup weight function.}

\bigbreak We now explain how the two-weight-function approach given in [Siu1996] is intimately related to the hyperbolic geometry of the punctured disk even though on the surface it does not seem so.  We will do it by discussing (i) the {\it a priori} estimate for the extension problem involving the closed positive $(1,1)$-current defined by $w=0$ and (ii) the problem and technique of shifting the blowup weight function from the norm for the right-hand side of the $\bar\partial$ equation to the norm for the solution of the $\bar\partial$ equation.  The latter will be discussed first, because it is the most natural way to get into the relation between the two-weight-function approach and the hyperbolic geometry of the punctured disk.  We start by looking at the standard cut-off function and the standard additional weight function from the historic example of producing holomorphic sections of line bundles with prescribed values (and jet values) in order to prove Kodaira's embedding theorem [Kodaira1954].

\medbreak\noindent(I.2.1) {\it Standard Cut-Off Function and Additional Weight Function.}  In Kodaira's proof of his embedding theorem [Kodaira1954], for the construction of holomorphic sections of a sufficiently positive line bundle with prescribed value at a point $P_0$ he used the blow-up at $P_0$, which after the later introduction of $L^2$ estimates of $\bar\partial$ by Morrey [Morrey1958], Kohn [Kohn1963-64], H\"ormander [H\"ormander1965], {\it et al}, was translated into the equivalent form of the use of a singular metric of strictly positive curvature current and singularity behavior $\approx\frac{1}{|w|^{2n}}$ on a coordinate chart $w=\left(w_1,\cdots,w_n\right)$ centered at $P_0$.  The standard procedure starts out with a local holomorphic section $f$ with the prescribed value $P_0$ defined on the ball $B_\varepsilon$ of radius $\varepsilon>0$ centered at $P_0$ and also with a smooth function $0\leq\rho_\varepsilon\leq 1$ with support in $B_\varepsilon$ which is identically $1$ in neighborhood of $P_0$.  From the solution $u_\varepsilon$ of the $\bar\partial$ equation with $f\bar\partial\rho_\varepsilon$ as the right-hand side, the global holomorphic section $F_\varepsilon$ with the prescribed value at $P_0$ is constructed as $f\rho_\varepsilon-u_\varepsilon$.  For the case $n=1$ (which we will now assume to simplify the discussion without sacrificing the key point of the argument) the significance of the singularity behavior $\approx\frac{1}{|w|^2}$ is that the dominant part $\frac{1}{2\pi\sqrt{-1}}\partial\bar\partial\log\left(\frac{1}{|w|^2}\right)$ of its curvature current is the Dirac delta $\delta_0$ at $w=0$.  The Dirac delta is what makes the solution of the problem possible.

\medbreak Here we need only use some fixed $\varepsilon$.  However, if there is a need to pass to limit as $\varepsilon\to 0$ (as in the case of extension results of Ohsawa-Takegoshi type), the norm for $F_\varepsilon$, which depends on $\varepsilon$, may blow up of the order $\frac{1}{\varepsilon}$ as $\varepsilon\to 0$.

\medbreak\noindent(I.2.2) {\it Need to Shift Blowup Function for Extension from Fiber.}   We now rehash part of the discussion in (I.1.3) and (I.1.5) by using the standard cut-off function and additional weight function.  Suppose $\Omega$ is a bounded Stein domain in ${\mathbb C}^n\times{\mathbb C}$ with $w$ as the coordinate of the second factor ${\mathbb C}$.  Suppose $\varphi$ is a plurisubharmonic function on $\Omega$ and $f$ is a holomorphic function on $\Omega\cap\left\{w=0\right\}$ with a finite $L^2$ norm $\left\|f\right\|_\varphi$ with respect to the weight function $e^{-\varphi}$ which we seek to extend to a holomorphic function $F$ on $\Omega$ whose $L^2$ norm with respect to $e^{-\varphi}$ is bounded by $\left\|f\right\|_\varphi$ times a positive constant which depends only on $\Omega$.

\medbreak For the approach of using the standard cut-off function and additional weight function, we can assume without loss of generality that for some $\varepsilon>0$ the holomorphic function $f$ is assumed to be defined on some open neighborhood of the topological closure of $\Omega\cap\left\{|w|<\varepsilon\right\}$ in ${\mathbb C}^n\times{\mathbb C}$ (as in (I.1.5)) and that $\varphi$ is smooth and defined on an open neighborhood of the topological closure of $\Omega$ in ${\mathbb C}^n\times{\mathbb C}$.  We take  $0 < \lambda  < 1$  and take a cut-off function  $\chi $
so that  $\chi (\xi )$  is identically  1  on  $\left\{\xi  \leq  \lambda\right\} $  and
$\xi $  is supported in  $\left\{\xi<1\right\}$.  Let  $\chi _\varepsilon (w) =
\chi ({|w|^2\over \varepsilon ^2})$.  Consider  $v_\varepsilon  = {1\over w}
\left(\overline \partial  \left(\chi _\varepsilon  f\right)\right) = {1\over w} \left(\overline \partial
\chi _\varepsilon \right) f$.  We solve the  $\overline \partial $-equation
$\overline \partial u_\varepsilon  = v_\varepsilon $ and set  $F = \chi _\varepsilon f -
w u_\varepsilon $ so that  $\overline \partial F = 0$  and  $F$  agrees with  $f$  on
$w = 0$.  The difficulty is to keep track of the estimates and make sure the final estimate can be made to be independent of $\varepsilon$.   Note that here the use of $v_\varepsilon=\frac{1}{w}\,\partial  \left(\chi _\varepsilon  f\right)$ is equivalent to using the additional weight function $\frac{1}{|w|^2}$ without division by $w$.

\medbreak Let  $U_\varepsilon  = \Omega  \cap  \{ \chi _\varepsilon  \neq  0 \}$.  Then
$$
\int _{U_\varepsilon} {|\overline \partial \chi _\varepsilon |^2\over |w|^2} |f|^2e^{-\varphi}
$$
is of the order  $1 \over {\varepsilon ^2}$, because of the factor $\frac{1}{w}$ in the definition
  of $v_\varepsilon$.  On the other hand, the solution $u_\varepsilon$ of
the $\bar\partial$ equation is used in $F = \chi _\varepsilon f -
w u_\varepsilon $ only after multiplication by the factor $w$.  The blowup weight function $\frac{1}{|w|^2}$ in the norm of the right-hand side $v_\varepsilon$ of the $\bar\partial$ equation is the reason for the difficulty and at the same time
the factor $w$ for the solution $u_\varepsilon$ of the $\bar\partial$ equation in its use in $F = \chi _\varepsilon f -
w u_\varepsilon $ is not being taken advantage of.  What is needed is to shift this blowup weight function $\frac{1}{|w|^2}$ from the norm of the right-hand side of the $\bar\partial$ equation to the norm of the solution of the $\bar\partial$ equation.  Note that it is precisely this blowup function $\frac{1}{|w|^2}$ (whose curvature provides the Dirac delta and which we would like to shift) that makes the extension problem solvable.  The two-weight-function method in [Siu1996] actually is to implement this shift of the blowup function.   As mentioned in (I.1.5) and (I.1.9) this step of shifting the blowup function corresponds to the trivial task of shifting a constant which is a power of $m$ in the approach of the hyperbolic geometry of the punctured disk.  We now explain how this shift of the blowup function in the two-weight-function method is related to the hyperbolic geometry of the punctured open $1$-disk.

\medbreak Technically, since products of currents occur in the discussion, any rigorous definition and argument involving such a product would need the use of smoothing of each individual factor current first and then take the limit after forming the product of the smoothings of individual factor currents.  For the simplicity of our discussion we will first carry out the multiplication of currents without the details of rigorous justification and then will discuss the smoothing process later.

\bigbreak\noindent(I.2.3) {\it Constant Negative Curvature of Metric of Punctured Unit Disk.} The metric for the punctured open unit $1$-disk $\Delta-\left\{0\right\}=\left\{0<|w|<1\right\}$ (with coordinate $w$) is $\frac{|dw|^2}{|w|^2\left(\log|w|^2\right)^2}$.  Its curvature on $\Delta-\left\{0\right\}$ is given by
$$
\partial\bar\partial\log\left(\frac{dw\wedge d\bar w}{\left|w\right|^2\left(\log|w|^2\right)^2}\right)=
\frac{2\,dw\wedge d\bar w}{|w|^2\left(\log|w|^2\right)^2}
$$
and is a negative constant.  For this note what is important for us is not its negative constant curvature on $\Delta-\left\{0\right\}$, but its behavior as $w\to 0$ as expressed in terms of the computation of its curvature on the full $\Delta$ in the sense of currents.  From
$$
\bar\partial\log\frac{1}{\left(\log|w|^2\right)^2}=-\frac{2}{\log |w|^2}\,
\bar\partial\log |w|^2
$$
on $\Delta$ in the sense of currents it follows that
$$
\displaylines{\partial\bar\partial\log\frac{1}{\left(\log|w|^2\right)^2} =
-\frac{2}{\log |w|^2}\partial\bar\partial\log|w|^2+\frac{2 \partial \log  |w|^2 \wedge
\overline {\partial \log  |w|^2}}
{\left(\log  |w|^2\right)^2}\cr
=\left(-\frac{2}{\log |w|^2}\right)\left(-2\pi\sqrt{-1}\right)\delta_0+\frac{2\,dw\wedge d\bar w}{|w|^2\left(\log|w|^2\right)^2},}
$$
where the differential version
$$\frac{1}{2\pi\sqrt{-1}}\,\bar\partial\left(\frac{1}{w}\right)=\delta_0$$ of Cauchy integral formula is used to get the Dirac delta $\delta_0$ (at the origin) in the equation.  In order to get the Dirac delta $\delta_0$ (up to a positive constant factor) as one term, we multiply both sides by $-\log|w|^2$ to get
$$
\left(-\log|w|^2\right)\,\sqrt{-1}\,\partial\bar\partial\left(-\log\left(\log\frac{1}{|w|^2}\right)\right)
=2\pi\delta_0+\frac{\sqrt{-1}\,dw\wedge d\bar w}{|w|^2\left(-\log|w|^2\right)}.\leqno{({\rm I}.2.3.1)}
$$
Note that $-\log|w|^2$ is nonnegative for $|w|<1$.  This formula (I.2.3.1) takes the place of the  formula
$\frac{\sqrt{-1}}{2\pi}\partial\bar\partial\left(-\log\left(\frac{1}{|z|^2}\right)\right)=\delta_0$ in the use of the $\bar\partial$ equation for the construction of holomorphic sections in the proof of Kodaira's embedding theorem.  The Dirac delta $\delta_0$ in the formula (I.2.3.1) is what makes the extension possible.

\medbreak The second term on the right-hand side, as well the first factor on the left-hand side, of formula (I.2.3.1) is what makes possible the shifting of the blowup function from the norm of the right-hand side of the $\bar\partial$ equation to the norm of the solution of the equation.  The expression on the left-hand side of formula (I.2.3.1) to which $\partial\bar\partial$ is applied is the $-\log$ of the additional weight function $\log\frac{1}{|w|^2}$ to be chosen for use, as we will see later.

\medbreak\noindent(I.2.3.2) The geometric meaning of the three terms in the formula (I.2.3.1) is as follows.  The left-hand side of (I.2.3.1) is the product of the metric $\log\frac{1}{|w|^2}$ with its curvature $\sqrt{-1}\,\partial\bar\partial\left(-\log\left(\log\frac{1}{|w|^2}\right)\right)$ (up to some universal positive constant factor).  The last term on the right-hand side of (I.2.3.1) is the absolute-value-square of the connection $\partial\log\left(\log\frac{1}{|w|^2}\right)$ of the metric $\log\frac{1}{|w|^2}$ times the metric $\log\frac{1}{|w|^2}$.

\medbreak The smoothing version of formula (I.2.3.1) (which we will need to use later to rigorously define products of currents) will be
$$
\displaylines{({\rm I}.2.3.3)_\varepsilon\qquad\qquad\left(\log\frac{1}{|w|^2+\varepsilon^2}\right)\,\sqrt{-1}\,\partial\bar\partial\log\left(\log\frac{1}{|w|^2+\varepsilon^2}\right)
\hfill\cr\hfill=\frac{\varepsilon^2\left(\sqrt{-1}\,dw\wedge d\bar w\right)}{\left(|w|^2+\varepsilon^2\right)^2}+\frac{\sqrt{-1}\,dw\wedge d\bar w}{|w|^2\left(\log\frac{1}{|w|^2+\varepsilon^2}\right)}\qquad\qquad}
$$
as $\varepsilon\to 0$, when the metric $\log\frac{1}{|w|^2}$ is replaced by its smoothing $\log\frac{1}{|w|^2+\varepsilon^2}$.

\medbreak What is actually needed from the formula (I.2.3.1) is the following inequality version
$$
\left(-\log|w|^2\right)\,\sqrt{-1}\,\partial\bar\partial\left(-\log\left(\log\frac{1}{|w|^2}\right)\right)
\geq c_0\delta_0+\frac{\sqrt{-1}\,dw\wedge d\bar w}{|w|^2\left(-\log|w|^2\right)}\leqno{({\rm I}.2.3.4)}
$$
for some positive number $c_0$.  By using the geometric interpretation given in (I.2.3.2), we can rewrite the inequality (I.2.3.4) as
$$
e^{-\kappa}\Theta_\kappa\geq c_0\delta_0+e^{-\kappa}\left(\sqrt{-1}\,\omega_\kappa\wedge\overline{\omega_\kappa}\right),\leqno{({\rm I}.2.3.5)}
$$
where $e^{-\kappa}$ is the metric $\log\frac{1}{|w|^2}$ and $\Theta_\kappa$ is the curvature $$\sqrt{-1}\,\partial\bar\partial\kappa=\sqrt{-1}\,\partial\bar\partial\left(-\log\left(\log\frac{1}{|w|^2}\right)\right)$$ of the metric $e^{-\kappa}$ (up to a standard normalizing factor of $\frac{1}{2\pi}$ when some other convention for the definition of curvature is used) and $\omega_\kappa$ is the connection $$\partial\kappa=\partial\log\left(\log\frac{1}{|w|^2}\right)=\frac{dw}{w\left(-\log|w|^2\right)}$$ of the metric $e^{-\kappa}$.  The following inequality from the hyperbolic geometry of the punctured disk is also important in determining the final constant in the estimate later
$$
\left|e^{-\kappa}\omega_\kappa\right|^2\leq\gamma\left|dw\right|^2\leqno{({\rm I}.2.3.6)}
$$
with $\gamma=\frac{1}{|w|^2}$.  This inequality is actually an identity, but what matters to us is only its inequality version.

\bigbreak\noindent(I.2.4) {\it Functional Analysis Setup for Shifting Blowup Function.}  To better understand the technique of shifting the blowup function from the right-hand side of the equation to its solution, we discuss it first from the formulation in functional analysis.  For the problem of extending a holomorphic function on $\Omega\cap\left\{w=0\right\}$ to $\Omega$ with $L^2$ estimates with respect to the weight function $e^{-\varphi}$, the ideal situation is that the curvature current $\Theta_\varphi$  dominates some positive multiple of the Dirac delta $\delta_0$ at $w=0$.  Here we use Dirac delta $\delta_0$ to mean both (i) the $(1,1)$-current on ${\mathbb C}$ defined by evaluation at the origin and (ii) the $(1,1)$-current on ${\mathbb C}^n\times{\mathbb C}$ defined by the codimension-one hyperplane $w=0$.   Which one is being meant will be clear from the context.

\medbreak With the given background metric $e^{-\varphi}$, in general we do not have the domination of $\delta_0$ by $\Theta_\varphi$ as $(1,1)$-currents on $\Omega$.  We introduce another metric $e^{-\kappa}$ whose curvature $\Theta_\kappa$ dominates positive multiple of the Dirac delta $\delta_0$ so that for a test form $g$ we have
$$
\displaylines{
({\rm I}.2.4.1)\qquad\qquad\left\|\bar\partial g\right\|^2_{\varphi+\kappa}+\left\|\bar\partial^*_{\varphi+\kappa} g\right\|^2_{\varphi+\kappa}\geq
\left<\left(\Theta_\varphi+\Theta_\kappa\right) g, g\right>_{\varphi+\kappa}\hfill\cr\qquad\qquad\qquad=\left(e^{-\kappa}\,\Theta_\varphi\,g, \,g\right)_\varphi+
\left(e^{-\kappa}\,\Theta_\kappa\, g, \,g\right)_\varphi,}
$$
where the subscript $\varphi+\kappa$ in $\left\|\bar\partial g\right\|_{\varphi+\kappa}$ indicates that the metric $e^{-(\varphi+\kappa)}$ is used in the definition of the norm and the subscript $\varphi$ in $\left(e^{-\kappa}\,\Theta_\varphi\,g, \,g\right)_\varphi$ indicates that the metric $e^{-\varphi}$ is used in the definition of the inner product.

\medbreak Though we have the advantage of the domination of a positive multiple of $\delta_0$ by $\Theta_\kappa$ from the use of the metric $e^{-\kappa}$, we need to get rid of the metric $e^{-\kappa}$ in the norms and the inner products and also change $\bar\partial^*_{\varphi+\kappa}$ to $\bar\partial^*_\varphi$, while keeping the domination of a positive multiple of $\delta_0$ by $\Theta_\kappa$ on the right-hand side, because our final estimate should only involve the weight function $e^{-\varphi}$.   In order to achieve this feat, we need some appropriate modifications of the estimate (I.2.4.1). For our modifications of (I.2.4.1) we assume
$$
e^{-\kappa}\Theta_\kappa\geq c_0\delta_0+e^{-\kappa}\left|\omega_\kappa\right|^2,\leqno{({\rm I}.2.4.2)}
$$
where $c_0$ is some universal positive constant and $\omega_\kappa=\partial\kappa$ is the connection of the metric $e^{-\kappa}$ and is a $(1,0)$-current.  This formula is precisely the inequality version (II.2.3.5) of the formula (II.2.3.1) according to the geometric interpretation of the terms of (II.2.3.1) given in (II.2.3.2).  For our modifications of (I.2.4.1) we also assume that we have a nonnegative function $\gamma$ such that $\left|e^{-\kappa}\omega_\kappa\right|^2\leq\gamma\Theta_\varphi$.  From $\bar\partial^*_{\varphi+\kappa}g=\bar\partial^*_\varphi g+\omega_\kappa\cdot g$ it follows that
$$
\left\|\bar\partial^*_{\varphi+\kappa} g\right\|^2_{\varphi+\kappa}=\left\|e^{-\frac{\kappa}{2}}\,\bar\partial^*_\varphi g\right\|^2_\varphi+2{\rm Re}\left(e^{-\kappa}\omega_\kappa\cdot g,\,\bar\partial_\varphi^* g\right)_\varphi+\left\|e^{-\frac{\kappa}{2}}\left(\omega_\kappa\cdot g\right)\right\|^2_{\varphi}.\leqno{({\rm I}.2.4.3)}
$$
Putting together the inequalities (I.2.4.1) and (I.2.4.2) and the identity (I.2.4.3) and the inequality $\left|e^{-\kappa}\omega_\kappa\right|^2\leq\gamma\Theta_\varphi$, we obtain
$$\left\|e^{-\frac{\kappa}{2}}\bar\partial g\right\|^2_\varphi+\left\|\left(\gamma+e^{-\kappa}\right)^{\frac{1}{2}}\bar\partial^*_\varphi g\right\|^2_\varphi\geq
c_0\left(\delta_0\,g\,g\right)_\varphi.\leqno{({\rm I}.2.4.4)}
$$
Of course, here when we have the product of currents, smoothing is needed for rigorous interpretation, which we will do later.  We introduce $\alpha=e^{-\frac{\kappa}{2}}$ and $\beta=\left(\gamma+e^{-\kappa}\right)^{\frac{1}{2}}$ to rewrite the above inequality (I.2.4.4) as
$$\left\|\alpha\,\bar\partial g\right\|^2_\varphi+\left\|\beta\,\bar\partial^*_\varphi g\right\|^2_\varphi\geq
c_0\left(\delta_0\,g,\,g\right)_\varphi.\leqno{({\rm I}.2.4.5)}
$$
To get to the solution of the extension problem, we still need one more technique, namely the technique of twisted $\bar\partial$ operators which we are going to explain.

\bigbreak\noindent(I.2.5) {\it Twisted $\bar\partial$ Operators.}  There are two techniques of $\bar\partial$ estimates which use several weight functions.  The first technique uses different weight functions in the usual complex of $\bar\partial$ operators.  For this technique, in the case of different weight functions (or metrics) $e^{-\varphi_\nu}$ for $\nu=0,1,2$, we study, in the following complex
$$
{\mathscr L}^2_{(0,0)}\left(e^{-\varphi_0}\right)\stackrel{T}{\longrightarrow}{\mathscr L}^2_{(0,1)}\left(e^{-\varphi_1}\right)\stackrel{S}{\longrightarrow}
{\mathscr L}^2_{(0,2)}\left(e^{-\varphi_2}\right),
$$
the solution of $Tu=f$ with $f$ satisfying the compatibility condition $Sf=0$, where ${\mathscr L}^2_{(0,\nu)}\left(e^{-\varphi_\nu}\right)$ means $(0,\nu)$-forms which are $L^2$ with respect to the weight function $e^{-\varphi_\nu}$ for $\nu=0,1,2$ and $T$ and $S$ are defined by $\bar\partial$.

\medbreak The second technique involves twisted $\bar\partial$ operators and uses two different sets of weight functions $e^{-\varphi_\nu}$ and $e^{-\kappa_\nu}$ for the space of $(0,\nu)$-forms to construct the following commutative diagram of two rows of complexes.
$$
\begin{matrix}
{\mathscr L}^2_{(0,0)}\left(e^{-\varphi_0}\right)&\stackrel{T}{\longrightarrow}&{\mathscr L}^2_{(0,1)}\left(e^{-\varphi_1}\right)&\stackrel{S}{\longrightarrow}&
{\mathscr L}^2_{(0,2)}\left(e^{-\varphi_2}\right)\cr
\quad\downarrow\Phi_0&&\downarrow\Phi_1&&\downarrow\Phi_2\cr
{\mathscr L}^2_{(0,0)}\left(e^{-\kappa_0}\right)&\stackrel{\tilde T}{\longrightarrow}&{\mathscr L}^2_{(0,1)}\left(e^{-\kappa_1}\right)&\stackrel{\tilde S}{\longrightarrow}&
{\mathscr L}^2_{(0,2)}\left(e^{-\kappa_2}\right),\cr\end{matrix}$$
where $\Phi_\nu$ is defined by multiplication by a positive function $\alpha_\nu$ and where $\tilde T$ and $\tilde S$ are defined to make the above diagram commutative.  In other words, $\tilde T(u)=\alpha_1\,T\left(\frac{u}{\alpha_0}\right)$ and $\tilde S(v)=\alpha_2\,S\left(\frac{v}{\alpha_1}\right)$.  When $S$ and $T$ are $\bar\partial$ operators, the operators $\tilde T$ and $\tilde S$ are {\it twisted $\bar\partial$ operators}.

\bigbreak\noindent(I.2.6) {\it Change of Complex Structures for Bundles.}  A twisted $\bar\partial$ operator can be also be described as the $\bar\partial$ operator for a new complex structure for the line bundle involved.  For a $(1,0)$-form $\omega$ which is $\partial$-closed, we can define a new operator $\bar\partial_\omega$ by $\bar\partial_\omega(u)=\bar\partial u+\bar\omega u$.  Since $\omega$ is $\partial$-closed, we can locally write $\omega=\partial\alpha$ for some function $\alpha$.  The new operator $\bar\partial_\omega$ can be rewritten in terms of the function $\alpha$ as a twisted $\bar\partial$ operator which sends $u$ to
$e^{-\alpha}\bar\partial\left(e^\alpha u\right)$.  When $u$ is considered a local section of a line bundle, a holomorphic local frame $f$ of the original complex structure (defined by the vanishing of $\bar\partial$) is multiplied by $e^{-\alpha}$ to become a holomorphic local frame of the new complex structure (defined by the vanishing of $\bar\partial_\omega$).  Later in (II.2.5) we use this view of the change of complex structures for bundles.

\bigbreak\noindent(I.2.7) {\it Solution from Estimate.}  We now return to our discussion of the extension problem.  Let $T$ denote the twisted $\bar\partial$ operator defined by $T(u)=\bar\partial\left(\beta\,u\right)$ for $(0,0)$-forms and let $S$ denote the twisted $\bar\partial$ operator defined by $S(v)=\alpha\,\bar\partial(v)$ for $(0,1)$-forms. Then the inequality (I.2.4.5) becomes
$$\left\|S g\right\|^2_\varphi+\left\|T^* g\right\|^2_\varphi\geq
c_0\left(\delta_0\,g\,g\right)_\varphi.\leqno{({\rm I}.2.7.1)}
$$
We need to interpret this inequality as the limit, as the positive number $\varepsilon\to 0$, of an inequality for corresponding smooth $\varphi_\varepsilon$, $\kappa_\varepsilon$, and $\gamma_\varepsilon$ with $\delta_0$ interpreted as the limit of
$$\frac{\sqrt{-1}}{2\pi}\partial\bar\partial\log\left(|w|^2+\varepsilon^2\right)=
\frac{\varepsilon^2}{\left(|w|^2+\varepsilon^2\right)^2}\left(\frac{
\sqrt{-1}\,dw\wedge d\bar w}{2\pi}\right)$$
so that
$$\left\|S_\varepsilon g\right\|^2_{\varphi_\varepsilon}+\left\|T_\varepsilon^* g\right\|^2_{\varphi_\varepsilon}\geq
c_0\left(\frac{\varepsilon^2}{\left(|w|^2+\varepsilon^2\right)^2}
\left(\frac{
\sqrt{-1}\,dw\wedge d\bar w}{2\pi}\right)\cdot g,\,g\right)_{\varphi_\varepsilon},\leqno{({\rm I}.2.7.2)_\varepsilon}
$$
where $S_\varepsilon(v)=e^{-\frac{\kappa_\varepsilon}{2}}\bar\partial(v)$ and $T_\varepsilon(u)=\bar\partial\left(\left(\gamma_\varepsilon+e^{-\kappa_\varepsilon}\right)^{\frac{1}{2}}\,u\right)$.
In this smoothing process we have to assume that for $\hat\varepsilon\leq\varepsilon$ the monotonicity condition $\varphi_\varepsilon\leq\varphi_{\hat\varepsilon}$ holds as usually the case in the smoothing procedure for plurisubharmonic functions.

\medbreak We follow the standard procedure, as in (I.1.5), of taking a sequence of Stein open subsets $\Omega_\nu$ of $\Omega$ exhausting $\Omega$ with the topologically closure $\bar\Omega_\nu$ of $\Omega_\nu$ relatively compact in $\Omega_{\nu+1}$ such that for each $\nu$ there exists a holomorphic retraction $\tau_\nu$ from an open neighborhood $\tilde\Omega_\nu$ of $\left\{w=0\right\}\cap\bar\Omega_\nu$ in $\Omega_{\nu+1}$ with $\Omega_\nu\cap\pi^{-1}\left(\left\{|w|\leq\varepsilon_\nu\right\}\right)$ contained in $\tilde\Omega_\nu$ for some $\varepsilon_\nu>0$ with $\lim_{\nu\to\infty}\varepsilon_\nu=0$.  We extend $f|_{\Omega_\nu}$ to some $\tilde f_\nu$ on $\Omega_\nu\cap\pi^{-1}\left(\left\{|w|<\varepsilon_\nu\right\}\right)$ by the holomorphic retraction $\tau_\nu$.  Then we solve the twisted $\bar\partial$ equation $$\bar\partial\left(\beta_{\varepsilon_\nu} u_{\nu,\delta}\right)=\frac{1}{w}\left(\bar\partial\chi_{\varepsilon_\nu}\right)\tilde f_\varepsilon$$ on $\Omega_\nu$ and get the estimate
$$
\int_{\Omega_\nu}\left|u_\nu\right|^2e^{-\varphi_{\varepsilon_\nu}}\leq \frac{1}{\,c_0\,}\int_{\Omega_\nu}e^{-\varphi_{\varepsilon_\nu}}
\frac{\left(|w|^2+\varepsilon_\nu^2\right)^2}{\varepsilon_\nu^2}\,\frac{\left|\bar\partial\chi_{\varepsilon_\nu}\right|^2}{|w|^2}\left|f\right|^2.$$
Finally the extension $F$ of $f$ is the limit of $\chi_{\varepsilon_\nu}\tilde f_\nu-w\beta_{\varepsilon_\nu} u_\nu$ when $\nu$ goes to $\infty$.

\medbreak Thus any holomorphic function $f$ on $\Omega\cap\left\{w=0\right\}$ with finite $L^2$ norm with respect to $e^{-\varphi}$ can be extended to a holomorphic function $F$ on $\Omega$ with the $L^2$ estimate
$$
\int_\Omega\left|F\right|^2e^{-\varphi}\leq\frac{1}{\,c_0\,} \sup_\Omega\left(|w|\beta\right)\int_{\Omega\cap\left\{w=0\right\}}\left|f\right|^2 e^{-\varphi}\leqno{({\rm I}.2.7.3)}
$$
if we assume that the curvature current $\Theta_\varphi$ dominates $\sqrt{-1}\,dw\wedge d\bar w$, because, by $({\rm I}.2.3.3)_\varepsilon$ from the hyperbolic geometry of the punctured open $1$-disk, the choices which we use are $e^{-\kappa}=\log\frac{1}{|w|^2}$ with $e^{-\kappa_\varepsilon}=\log\frac{1}{|w|^2+\varepsilon^2}$ as its smoothing and $\gamma=\frac{1}{|w|^2}$ with $\gamma_\varepsilon=\frac{1}{|w|^2+\varepsilon^2}$ as its smoothing.  Note that when the condition $\Theta_\varphi\geq\sqrt{-1}\,dw\wedge d\bar w$ is not satisfied, in order to guarantee that
$\left|e^{-\kappa}\omega_\kappa\right|^2\leq\gamma\Theta_\varphi$ holds, it is necessary to replace $e^{-\varphi}$ by $e^{-\varphi}e^{-|w|^2}$ first and then to get rid of the factor $e^{-|w|^2}$ at the end with a change of the constant $\frac{1}{\sqrt{c_0\,}}$ in (I.2.7.3).

\medbreak We now summarize as follows our discussion on the {\it a priori} estimate for the extension problem and the technique of shifting blowup functions.

\bigbreak\noindent(I.2.8) {\it Theorem.} Let $\Omega$ be a Stein open subset of ${\mathbb C}^n\times\Delta$ and $\varphi$ be a plurisubharmonic function on $\Omega$.  Let $c_0$ be a positive number.  Let $\left\|\cdot\right\|_\varphi$ denote $L^2$ norm on $\Omega$ with respect to the weight function $e^{-\varphi}$.

\medbreak\noindent(a) Let $\alpha$ and $\beta$ be nonnegative functions on $\Omega$.  If $\left\|Sg\right\|^2_\varphi+\left\|T^*g\right\|^2_\varphi\geq c_0\left(\delta_0\,g,\,g\right)_\varphi$ for test $(0,1)$-form $g$, where $S(v)=\alpha\bar\partial(v)$ and $T(u)=\bar\partial(\beta u)$, then a holomorphic function $f$ on $\Omega\cap\left\{w=0\right\}$ with finite $L^2$ norm with respect to $e^{-\varphi}$ can be extended to a holomorphic $F$ on $\Omega$ with
$$
\int_\Omega\left|F\right|^2 e^{-\varphi}\leq\frac{1}{c_0}\sup\left(|w|\beta\right)\int_{\Omega\cap\left\{w=0\right\}}\left|f\right|^2 e^{-\varphi}.
$$

\medbreak\noindent(b) If there is a metric $e^{-\kappa}$ on $\Omega$, with connection $\omega_\kappa=\partial\kappa$ and curvature current $\Theta_\kappa$, which satisfies
$$
e^{-\kappa}\Theta_\kappa\geq c_0\delta_0+e^{-\kappa}\left|\omega_\kappa\right|^2,\leqno{({\rm I}.2.8.1)}
$$
and if $\gamma$ is a nonnegative function on $\Omega$ such that
$$\left|e^{-\kappa}\omega_\kappa\right|^2\leq\gamma\Theta_\varphi,\leqno{({\rm I}.2.8.2)}
$$
then
the inequality $\left\|Sg\right\|^2_\varphi+\left\|T^*g\right\|^2_\varphi\geq c_0\left(\delta_0\,g,\,g\right)_\varphi$ holds for test $(0,1)$-form $g$ with
$\alpha=e^{-\frac{\kappa}{2}}$ and $\beta=\left(\gamma+e^{-\kappa}\right)^{\frac{1}{2}}$.

\medbreak\noindent Here the inequalities involving products of currents are to be interpreted by applying smoothing processes to the currents so that all the relevant inequalities hold for the smoothings of the currents, subject to the two conditions that the smoothing of $\varphi$ dominates $\varphi$ and the smoothing of the Dirac delta $\delta_0$ has to be chosen to be
$$\frac{\sqrt{-1}}{2\pi}\partial\bar\partial\log\left(|w|^2+\varepsilon^2\right)=
\frac{\varepsilon^2}{\left(|w|^2+\varepsilon^2\right)^2}\left(\frac{
\sqrt{-1}\,dw\wedge d\bar w}{2\pi}\right)$$
which approaches $\delta_0$ as $\varepsilon\to 0$.

\bigbreak

\bigbreak\noindent{\bf Part II.  Holomorphic Family of Flat Bundles}.

\bigbreak In Part II we first prove an extension result of Ohsawa-Takegoshi type in a special situation where the curvature current is not even semipositive but only semipositive along the fibers with the additional assumption that its curvature current is bounded from below and its connection is bounded (Theorem (II.2.10)).  Then we use the extension result and apply the arithmetic argument of Gelfond-Schneider [Gelfond1934, Schneider1934], Lang [Lang1962, Lang1965, Lang1966], Bombieri [Bombieri1970, Bombieri-Lang1970], Briskorn [Brieskorn1970], and Simpson [Simpson1993] to the two descriptions of the moduli of flat line bundles (one by the use of embedding into projective space and the other by the use of transition functions) to prove the following structure result for the subvariety in the moduli space of flat bundles defined by the dimension of the space of flatly twisted pluricanonical sections.

\bigbreak\noindent(II.1) {\it Theorem.} Let $X$ be a compact complex algebraic manifold and $m$ be a positive integer.  Let $Z_{m,q}$ be the set of all flat line bundles $F$ on $X$ such that $\dim_{\mathbb C}\Gamma\left(X, mK_X+F\right)\geq q$.  Then $Z_{m,q}$ is regular and is a finite union of translates of abelian subvarieties by torsion elements in the abelian variety of all flat line bundles.  In particular, there exists some $F_0\in Z_{m,q}$ which is a torsion element, {\it i.e.,} $NF_0$ is the trivial line bundle for some positive integer $N$.

\bigbreak\noindent{\bf \S1. Natural metrics of family of flat line bundles.}

\bigbreak\noindent(II.1.1) {\it Explicit Description of Transition Functions and Metric.}  Let $X$ be a compact complex algebraic manifold of complex dimension $n$.  From the long cohomology exact sequence
$$
H^1\left(X,{\mathbb Z}\right)\stackrel{\iota}{\longrightarrow}H^1\left(X, {\mathcal O}_X\right)\stackrel{\exp\left(\cdot\right)}{\longrightarrow}H^1\left(X,{\mathcal O}_X^*\right)\stackrel{c_1(\cdot)}{\longrightarrow} H^2\left(X,{\mathbb Z}\right)
$$
of the short exact sequence
$$
0\to {\mathbb Z}\hookrightarrow{\mathcal O}_X\stackrel{\exp\left(\cdot\right)}{\longrightarrow}{\mathcal O}_X^*\to 0
$$
of sheaves (where ${\mathbb Z}\hookrightarrow{\mathcal O}_X$ is defined by multiplication by $2\pi\sqrt{-1}$, and $\exp(\cdot)$ is the exponential map, and $c_1(\cdot)$ is the first Chern class map) it follows that the abelian variety
$$A=H^1\left(X, {\mathcal O}_X\right)\left/{\rm Im}\left(H^1\left(X, {\mathbb Z}\right)\stackrel{\iota}{\longrightarrow}H^1\left(X, {\mathcal O}_X\right)\right)\right.
$$
is the kernel of $H^1\left(X,{\mathcal O}_X^*\right)\stackrel{c_1(\cdot)}{\longrightarrow} H^2\left(X,{\mathbb Z}\right)$.

\medbreak Every point $a$ of $A$ is represented by an element $\tilde a$ of $H^1\left(X, {\mathcal O}_X\right)$ whose image in $H^1\left(X,{\mathcal O}_X^*\right)$ under the exponential map $\exp(\cdot)$ defines a holomorphic line bundle $F^{(a)}$ over $X$ whose first Chern class is zero and therefore is flat.  By Hodge decomposition $H^1\left(X, {\mathcal O}_X\right)$ is isomorphic to $\Gamma\left(X,\Omega_X^1\right)$ by complex conjugation.   The element $\tilde a$ corresponds to the complex-conjugate of some holomorphic $1$-form $\omega_a$ on $X$.  Fix a finite open cover ${\mathcal U}=\left\{U_j\right\}_{j\in J}$ of $X$ so that the intersection of any finite number of elements of ${\mathcal U}$ is simply connected and Stein. We are going to explicitly describe as follows the flat line bundle $F^{(a)}$ in terms of constant transition functions with respect to ${\mathcal U}$ which are obtained from the constants of integration of the indefinite integral of $\omega_a$. There exists a holomorphic function $f_{j,\omega_a}$ on $U_j$ such that $df_{j,\omega_a}=\omega_a$ on $U_j$ and $f_{j,\omega_a}-f_{k,\omega_a}=c_{kj,\omega_a}$ for some constant $c_{jk,\omega}$ on $U_j\cap U_k$.  The flat line bundle $F^{(a)}$ is given by the transition functions $\left\{\exp\left(\overline{c_{jk,\omega_a}}\right)\right\}_{j,k}$.  Note that the complex conjugate $\overline{c_{jk,\omega_a}}$ of $c_{jk,\omega_a}$ is used, because the element $\tilde a$ of $H^1\left(X, {\mathcal O}_X\right)$ corresponds to the complex conjugate $\overline{\omega_a}$ of the element $\omega_a$ of $\Gamma\left(X,\Omega_X^1\right)$ under the Hodge decomposition.

\medbreak We would like to remark that if we replace $\exp\left(\cdot\right)$ by $\exp\left(2\pi\sqrt{-1}\left(\cdot\right)\right)$ and at the same time replace the map ${\mathbb Z}\hookrightarrow{\mathcal O}_X$ of multiplication by $2\pi\sqrt{-1}$ by just the usual inclusion ${\mathbb Z}\hookrightarrow{\mathcal O}_X$ ({\it i.e.,} of multiplication by $1$), the arguments and the conclusions which we are presenting will not be affected, because our only concern is the results (II.2.7) and (II.2.8) where the factor $2\pi\sqrt{-1}$ plays no r\^ole.  However, in the discussion of \S5 where the moduli defined by transition functions are used, the presence or absence of factors like $2\pi\sqrt{-1}$ in the exponent of the transition functions would make a great difference, because there our interest is in whether the transition functions can be chosen to be roots of unity, with the transition functions coming from the consideration of
${\rm Hom}\left(\pi_1(X), {\mathbb C}^*\right)={\rm Hom}\left(H_1\left(X,{\mathbb Z}\right), {\mathbb C}^*\right)$.

\medbreak We have a family of flat line bundles $F^{(a)}$ over $X$ parametrized by $a\in A$ with flat transition functions as described above.  This family $\left\{F^{(a)}\right\}_{a\in A}$ defines a holomorphic line bundle ${\mathcal F}$ over $X\times A$ such that the restriction of ${\mathcal F}$ to $X\times\left\{a\right\}$ can be naturally identified with $F^{(a)}$.  Consider the holomorphic family $\pi:X\times A\to A$ with the same fiber $X$ under the natural projection onto the second factor.   We are interested in applying the extension result of Ohsawa-Takegoshi type to the holomorphic line bundle ${\mathcal F}$ over $X\times A$ in the sense that for $a\in A$ we study the condition under which an element of $\Gamma\left(X,F^{(a)}+K_X\right)$ can be extended to $\Gamma\left(\pi^{-1}(W), {\mathcal F}+K_{X\times A}\right)$ for some open neighborhood $W$ of $a$ in $A$.

\medbreak\noindent(II.1.1.1) Note that such an extension to an open subset of the base manifold $A$ is not always possible, because though each $F^{(a)}$ is flat and admits a metric with zero curvature, yet the holomorphic line bundle ${\mathcal F}$ in general does not admit a metric whose curvature is semi-positive, preventing a direct application of an extension of Ohsawa-Takegoshi type.  A simple counter-example is the special case where the canonical line bundle $K_X$ of $X$ is trivial, for example, when $X$ is itself an abelian variety.  If we choose $a\in A$ such that $F^{(a)}$ to be the trivial line bundle on $X$, then the element of $\Gamma\left(X,K_X+F^{(a)}\right)$ represented by a nonzero constant function on $X$ cannot be extended to an element of $\Gamma\left(\pi^{-1}(W), {\mathcal F}+K_{X\times A}\right)$ for any open neighborhood $W$ of $a$ in $A$, otherwise there exists some $b\in W$ with the flat line bundle $F_b$ nontrivial and yet admitting a nowhere zero holomorphic section over $X$.

\medbreak Let $\tilde A$ be the universal cover of $A$.  There is a natural way of constructing a metric on $F^{(a)}$ with zero curvature so that the construction yields a metric on the pullback $\tilde{\mathcal F}$ of ${\mathcal F}$ to $X\times\tilde A$.  We are going to explicitly write down this metric of $\tilde{\mathcal F}$ and determine the extent to which its curvature on $X\times\tilde A$ fails to be semi-positive. From $f_{j,\omega_a}-f_{k,\omega_a}=c_{kj,\omega_a}$ it follows that $$\exp\left(\overline{f_{j,\omega_a}}\right)=\exp\left(\overline{f_{k,\omega_a}}\right)\exp\left(\overline{c_{kj,\omega_a}}\right).$$
A metric for $F^{(a)}$ can be defined by a smooth positive function $h_j$ on $U_j$ such that $h_k\left|\exp\left(\overline{c_{kj,\omega_a}}\right)\right|^2=h_j$ on $U_j\cap U_k$.  Thus we can set
$h_j=\exp\left(2{\rm Re}\,\overline{f_{j,\omega_a}}\right)$ to get a metric of $F^{(a)}$.

\medbreak Let $\omega_1,\cdots,\omega_\ell$ be a ${\mathbb C}$-basis of holomorphic $1$-forms on $X$ and let $t_1,\cdots,t_\ell$ be the holomorphic Euclidean coordinates for $\tilde A$, the flat line bundle corresponding to $\left(t_1,\cdots,t_\ell\right)$ is $F^{(a)}$ with $\overline{\omega_a}=\sum_{\nu=1}^\ell t_\nu\overline{\omega_\nu}$.  Again, note that we use the complex-conjugate $\overline{\omega_\nu}$ of $\omega_\nu$ and yet use $t_\nu$ instead of its complex-conjugate $\overline{t_\nu}$, because the complex structure for $A$ is induced from $H^1\left(X,{\mathcal O}_X\right)$ and yet $H^1\left(X,{\mathcal O}_X\right)$ is isomorphic to the complex-conjugate of $\Gamma\left(X,\Omega_X^1\right)$.
We can use $\left\{\exp\left(\sum_{\nu=1}^\ell t_\nu\overline{c_{jk,\omega_\nu}}\right)\right\}_{j,k\in J}$ as the transition functions for $\tilde{\mathcal F}$ with respect to the covering $\left\{U_j\times\tilde A\right\}_{j\in J}$.  We can use as metric the positive function\ \  $\exp\left(2{\rm Re}\sum_{\nu=1}^\ell t_\nu\overline{f_{j,\omega_\nu}}\right)$ on $U_j\times\tilde A$.  Since $f_{jk,\omega_\nu}$ is a holomorphic function on $U_j$, the occurrence of the product of its complex conjugate $\overline{f_{jk,\omega_\nu}}$ and $t_\nu$ in the definition of the metric for $\tilde{\mathcal F}$ prevents the curvature of the metric to be semipositive on $X\times\tilde A$.

\medbreak If we can change $t_\nu$ to a holomorphic function of $\overline{t_\nu}$ in the metric, then the curvature of this metric would be zero on $X\times\tilde A$.  Unfortunately this cannot be done.  However, if we restrict ourselves to a circle $t_\nu\overline{t_\nu}=r^2$ for some positive constant $r$, then $t_\nu=r^2\left/\overline{t_\nu}\right.$ depends holomorphically on $\overline{t_\nu}$.  This observation motivates us to consider the following situation.  First take a flat line bundle $L^{(0)}$ on $X$ with constant transition functions $g^{(0)}_{jk}$ with respect to the cover ${\mathcal U}$ of $X$.  Let $h^{(0)}_j$ be a positive constant for $U_j$ which defines a metric of zero curvature for the flat line bundle $L^{(0)}$.   Then we take a circle in $\tilde A$ of radius $1$ on the complex line in $\tilde A$ defined by a holomorphic $1$-form $\omega$ such that the center of the circle corresponds to the flat line bundle $L^{(0)}$ and a variable point on the circle is parametrized by $\tau\in{\mathbb C}$ with $|\tau|=1$ which corresponds to the flat line bundle $L^{(\tau)}$ on $X$ whose transition functions are given by $g^{(0)}_{jk}\exp\left(t\overline{c_{jk}}\right)$ on $U_j\cap U_k$, where $c_{jk}$ is the constant $f_k-f_j$ on $U_j\cap U_k$ with $df_j=\omega$ on $U_j$.  The metric for $L^{(\tau)}$ is given by the positive function $h^{(0)}_j\exp\left(2{\rm Re}\left(\tau\overline{f_j}\right)\right)$ on $U_j$.  We can extend the definition of $L^{(\tau)}$ and its metric to all $t\in{\mathbb C}$ without the restriction $|t|=1$.  We put the family of flat line bundles $L^{(\tau)}$ on $X$ for $\tau\in{\mathbb C}$ together to form a holomorphic line bundle ${\mathcal L}$ over $X\times{\mathbb C}$.

\bigbreak\noindent(II.1.2) {\it Modification of Metrics to Yield Nonnegative Curvature.}  Take an ample nonsingular divisor $D$ in $X$ whose associated line bundle over $X$, also denoted by $D$, admits a smooth metric $h_D$ with positive curvature which is $\geq\frac{\sqrt{-1}}{2\pi}\omega\wedge\bar\omega$ on $X$.  Let $s_D=\left\{s_{j,D}\right\}_{j\in J}$ be the canonical section of the line bundle $D$ over $X$ so that the divisor of $s_D$ is the divisor $D$, where $s_{j,D}$ is the holomorphic function on $U_j$ representing $s_D$ on $U_j$.  We use $\tilde D$ to denote the divisor $D\times{\mathbb C}$ in $X\times{\mathbb C}$ as well as its associated line bundle on $X\times{\mathbb C}$ and give the line bundle $\tilde D$ the metric $h_{\tilde D}$ obtained by pulling back $h_D$ by the natural projection $X\times{\mathbb C}\to X$ onto the first factor.  Let $s_{\tilde D}=\left\{s_{j,\tilde D}\right\}_{j\in J}$ be the canonical section of the line bundle $\tilde D$ over $X\times{\mathbb C}$ so that the divisor of $s_{\tilde D}$ is the divisor $\tilde D=D\times{\mathbb C}$, where $s_{j,\tilde D}$ is the holomorphic function on $U_j\times{\mathbb C}$ representing $s_{\tilde D}$ on $U_j\times{\mathbb C}$.  Write $h_{\tilde D}=e^{-\chi_j}$ on $U_j\times{\mathbb C}$.  Since
$\omega=\sum_{k=1}^n\frac{\partial f_j}{\partial z_k}\,dz_k$ on $U_j$, it follows that for any number $\eta>0$ the curvature of the metric
$$
\displaylines{\left(h_{\tilde D}\right)^\eta h^{(0)}_j\exp\left(2{\rm Re}\left(\bar\tau f_j\right)\right)e^{-\frac{1}{2\eta}\left(|\tau|^2-1\right)}\cr=
h^{(0)}_j\exp\left(-\eta\chi_j+2{\rm Re}\left(\bar\tau f_j\right)-\frac{1}{2\eta}\left(|\tau|^2-1\right)\right)}
$$
on $U_j\times{\mathbb C}$ is nonnegative and the curvature of the metric
$$
\displaylines{\left(h_{\tilde D}\right)^{2\eta} h^{(0)}_j\exp\left(2{\rm Re}\left(\bar\tau f_j\right)\right)e^{-\frac{1}{2\eta}\left(|\tau|^2-1\right)}\cr=
h^{(0)}_j\exp\left(-2\eta\chi_j+2{\rm Re}\left(\bar\tau f_j\right)-\frac{1}{2\eta}\left(|\tau|^2-1\right)\right)}
$$
on $U_j\times{\mathbb C}$ dominates $\frac{\sqrt{-1}}{4\pi\eta}\,d\tau\wedge d\bar\tau$ (when the normalizing constant in the definition of curvature is chosen so that the curvature of the metric $e^{-\psi}$ is $\frac{\sqrt{-1}}{2\pi}\partial\bar\partial\psi$).  We are going to apply the technique of section extension of Ohsawa-Takegoshi type to the line bundle ${\mathcal L}+K_{X\times{\mathbb C}}$ over the Stein manifold $X\times{\mathbb C}-\tilde D$ and to use the metric
$$
\left(h_{\tilde D}\right)^{2\eta} h^{(0)}_j\exp\left(2{\rm Re}\left(\bar\tau f_j\right)\right)e^{-\frac{1}{2\eta}\left(|\tau|^2-1\right)}\left|s_{j,\tilde D}\right|^{4\eta}.\leqno{({\rm II}.1.2.1)}
$$
for ${\mathcal L}$ on $\left(X-D\right)\times{\mathbb C}=X\times{\mathbb C}-\tilde D$.  Note that we have to add the factor $\left|s_{j,\tilde D}\right|^{4\eta}$ in order to allow the use of $\left(h_{\tilde D}\right)^{2\eta}$ and still keep the entire expression (II.1.2.1) as a metric for ${\mathcal L}$ on $X\times{\mathbb C}-\tilde D$.  Because $s_{\tilde D,j}$ has zeroes at the divisor $\tilde D$, we have to exclude the hypersurface $\tilde D$ from $X\times{\mathbb C}$ to consider $X\times{\mathbb C}-\tilde D$ instead of all of $X\times{\mathbb C}$ so that the curvature current of the metric in (II.1.2.1) is nonnegative on $X\times{\mathbb C}-\tilde D$ even though the factor $\left|s_{j,\tilde D}\right|^{4\eta}$ occurs in (II.1.2.1).

\bigbreak

\bigbreak\noindent{\bf \S2. Extension with controllably negative curvature.}

\bigbreak\noindent(II.2.1) {\it Problem of Extension of Section on Initial Fiber.}  Recall that our notation $\Delta_R$ for $R>0$ means the open disk of radius $R$ in ${\mathbb C}$ centered at the origin with $\Delta_1$ simply abbreviated to $\Delta$.  Fix some $0<\eta_0<1$.  Assume that for some complex number $\tau_0$ in the boundary $\partial\Delta$ of $\Delta$ there exists some nonzero element $s_{\tau_0}\in\Gamma\left(X,K_X+L^{\left(\tau_0\right)}\right)$.  We would like to study the extension problem, namely under what condition it is possible to find a sequence of pairs $\tau_\nu\in\partial\Delta-\left\{\tau_0\right\}$ and
$s_{\tau_\nu}\in\Gamma\left(X, K_X+L^{(\tau_\nu)}\right)$ for $\nu\in{\mathbb N}$ such that $\tau_\nu\to\tau_0$ as $\nu\to\infty$ and $s_{\tau_\nu}$ approaches $s_{\tau_0}$ as $\nu\to\infty$.  As observed earlier in (I.1.1.1) in general this is not possible without any additional condition even in the simple case where $X$ is an abelian variety and $L^{\left(\tau_0\right)}$ is the trivial line bundle.  The key tool we have at our disposal is that because $\tau_0\in\partial\Delta$, when $\tau=\tau_0$ the weight function $e^{-\frac{1}{\eta}\left(\left|\tau\right|^2-1\right)}$ is identically equal to $1$ for all values of $\eta>0$ so that the $L^2$ norm of the element $s_{\tau_0}\in\Gamma\left(X,K_X+L_{\tau_0}\right)$ with respect to the metric (II.1.2.1) is uniformly bounded by some positive number $C_0$ for all $0<\eta\leq\eta_0$.

\medbreak Fix $R>1$ and we will consider the line bundle ${\mathcal L}$ on the product of $X$ and the disk $\Delta_R$.  Since the curvature current of the metric (II.1.2.1) dominates $\frac{\sqrt{-1}}{4\pi\eta}\,d\tau\wedge d\bar\tau$, the Ohsawa-Takegoshi type extension result allows us to extend the restriction to $X-D$ of the element $s_{\tau_0}$ of $\Gamma\left(X,K_X+L^{\left(\tau_0\right)}\right)$ to an element
$\tilde s_\eta\in\Gamma\left((X-D)\times\Delta_R, K_{X\times{\mathbb C}}+{\mathcal L}\right)$ for $0<\eta\leq\eta_0$ such that
$$
\int_{(X-D)\times\Delta_R}\left|\tilde s_\eta\right|^2\left(h_{\tilde D}\right)^{2\eta} h^{(0)}_j\exp\left(2{\rm Re}\left(\bar\tau f_j\right)\right)e^{-\frac{1}{2\eta}\left(|\tau|^2-1\right)}\left|s_{j,\tilde D}\right|^{4\eta}\leq C\eta,\leqno{({\rm II}.2.1.1)}
$$
where $C$ is a constant which is independent of $\eta$.  Note that the factor $\eta$ on the right-hand side of (II.2.1.1) comes from the $\eta$ in the denominator of $\frac{\sqrt{-1}}{4\pi\eta}\,d\tau\wedge d\bar\tau$ which is dominated by the curvature current of the metric (II.1.2.1).

\bigbreak\noindent(II.2.2) {\it Annulus Shrinking to Circle.}  For $0<\eta\leq\eta_0$ let $$G_\eta=\left\{\,\tau\in{\mathbb C}\,\Big|\,1-\eta<\left|\tau\right|<1+\eta\,\right\}$$ denote the open annulus in ${\mathbb C}$ of inner radii $1-\eta$ and outer radius $1+\eta$ centered at the origin.  We replace $\Delta_R$ in (II.2.1.1) by $G_\eta$ (for $0<\eta\leq\eta_0$) to get the weaker inequality
$$
\int_{(X-D)\times G_\eta}\left|\tilde s_\eta\right|^2\left(h_{\tilde D}\right)^{2\eta} h^{(0)}_j\exp\left(2{\rm Re}\left(\bar\tau f_j\right)\right)e^{-\frac{1}{2\eta}\left(|\tau|^2-1\right)}\left|s_{j,\tilde D}\right|^{4\eta}\leq C\eta,\leqno{({\rm II}.2.2.1)}
$$
so that on $G_\eta$ the function $\frac{1}{\eta}\left(\left|\tau\right|-1\right)$ is bounded by $1$. Thus on $X\times G_\eta$ the metric $\left(h_{\tilde D}\right)^{2\eta} h^{(0)}_j\exp\left(2{\rm Re}\left(\bar\tau f_j\right)\right)e^{-\frac{1}{2\eta}\left(|\tau|^2-1\right)}\left|s_{j,\tilde D}\right|^{4\eta}$ from (II.1.2.1), which will be used as a weight function in an integral, has a bound on $X\times G_\eta$ which is independent of $0<\eta\leq\eta_0$.
We divide both sides of the equation (II.2.2.1) by $\eta$ and let $\eta\to 0$ (through an appropriate sequence) to obtain
$$
\int_{(X-D)\times \partial\Delta}\left|\tilde s_0\right|^2 h^{(0)}_j\exp\left(2{\rm Re}\left(\bar\tau f_j\right)\right)\leq\hat C
$$
for some constant $\hat C$ which is equal to $C$ times some universal positive number,
where $\tilde s_0$ is the limit of $\tilde s_\eta$ as $\eta\to 0$ through an appropriate sequence from the weak compactness.  The difficulty is that the limit $\tilde s_0$ is only $L^2$ on $X\times\partial\Delta$ (though holomorphic along the direction of $X$) and, as a result, we cannot conclude that we can find a sequence $\tau_\nu\in\partial\Delta$ (for $\nu\in{\mathbb N}$) approaching $\tau_0$ such that the restriction $\tilde s_0|_{X\times\left\{\tau_\nu\right\}}\in\Gamma\left(X,K_X+L^{\left(\tau_\nu\right)}\right)$ of $\tilde s_0$ to $X\times\left\{\tau_\nu\right\}$ (after the
natural identification of $X\times\left\{\tau_\nu\right\}$ with $X$) converges to $s_{\tau_0}$ as $\nu\to\infty$.

\medbreak In order to guarantee the convergence of $\tilde s_0|_{X\times\left\{\tau_\nu\right\}}\in\Gamma\left(X,K_X+L^{\left(\tau_\nu\right)}\right)$ to $s_{\tau_0}$ as $\nu\to\infty$, we replace the $L^2$ norm by a stronger norm.  The stronger norm is the $L^2$ of derivatives up to order $1$ for differentiation along the vector field of unit-angular-speed rotation of ${\mathbb C}$ (whose coordinate is $\tau$).  In [Popvici2005] Popovic obtained $L^2$ Sobolev estimates for extension results of Ohsawa-Takegoshi type even for the general case of Hermitian vector bundles (whose $L^2$ case without derivatives was treated earlier by Manivel [Manivel1993] as a generalization of the original extension result of Ohsawa-Takegoshi [Ohsawa-Takegoshi1987]).  However, such general formulations cannot be directly applied to our special situation of a line bundle whose curvature is not even semi-positive and whose domain involves an annulus with outer and inner radii approaching $1$ in passing to limit.  We will first comment on the standard general procedure for $L^2$ estimates of $\bar\partial$ with norms involving derivatives, but we will not simply apply the standard general procedure.  Instead we will explicitly work out the derivative estimate for our special situation at hand by using a vector bundle of rank $2$ on the domain involving an annulus.

\bigbreak\noindent(II.2.3) {\it Estimates for Norm of Differentiation by Vector Field.}  We now comment on the standard general procedure for $L^2$ estimates of $\bar\partial$ with norms involving derivatives and its adaption to our situation, but we will not follow through with the implementation of its details in our case, because the explicit derivation of the derivative estimate for our our special situation at hand which we will give later in (II.2.6) is much simpler and clearer.

\medbreak For the complex coordinate $\tau=re^{\sqrt{-1}\theta}$ of ${\mathbb C}$, from
$$\tau\frac{d}{d\tau}=\frac{d}{d\log\tau}=\frac{1}{2}\left(\frac{\partial}{\partial\log r}-\sqrt{-1}\,\frac{\partial}{\partial\theta}\right)=\frac{1}{2}\left(r\frac{\partial}{\partial r}-\sqrt{-1}\,\frac{\partial}{\partial\theta}\right)$$ we introduce the vector field $\xi={\rm Re}\left(2\sqrt{-1}\,\tau\frac{d}{d\tau}\right)=\frac{\partial}{\partial\theta}$ of rotation of ${\mathbb C}$ with unit angular speed which is the real part of the holomorphic vector field $2\sqrt{-1}\,\tau\frac{d}{d\tau}$ of ${\mathbb C}^*={\mathbb C}-\left\{0\right\}$. For a $\left(K_{X\times{\mathbb C}}+{\mathcal L}\right)$-valued $(0,\nu)$-form $g$ on $X\times G_\eta$ we use the norm $\left\|g\right\|^2_{X\times G_\eta}+\gamma\left\|\nabla_\xi g\right\|^2_{X\times G_\eta}$ with $\gamma$ being a positive constant, where $\left\|\cdot\right\|_{X\times G_\eta}$ means the $L^2$ norm on $X\times G_\eta$ and $\nabla_\xi$ denotes the covariant differentiation with respect to the vector field by using the connection from the metric of the line bundle.  The vector field $\xi$ is chosen for the following two reasons.
\begin{itemize}
\item[(i)]  The domain $G_\eta$ is invariant under the vector field $\xi$ so that the condition of $g$ being in the domain of $\bar\partial^*$ is preserved under differentiation by $\xi$.
\item[(ii)] The factor $e^{-\frac{1}{2\eta}\left(|\tau|^2-1\right)}$ of the metric used in the $L^2$ estimate (II.1.2.1) is invariant under $\xi$ so that after differentiation there is no new contribution from $\frac{1}{\eta}$ which would blow up when we take the limit $\eta\to 0$.
\end{itemize}
In the solution of the equation $Tu=f$ subject to $Sf=0$ with the $L^2$ estimate
$$
\left\|Sg\right\|^2+\left\|T^*g\right\|^2\geq c_0\left\|g\right\|^2
$$
for the test function $g$ (where $c_0$ is a positive number), the standard technique of substituting the $L^2$ norm $\left\|g\right\|$ by a stronger norm involving the $L^2$ of derivatives (for example, by $\left\|g\right\|_\xi$ with $\left\|g\right\|_\xi^2=\left\|g\right\|^2+\gamma\left\|\nabla_\xi g\right\|^2$ in our case) is to apply the $L^2$ estimate to $\nabla_\xi g$ (instead of just $g$) to get
$$
\left\|S\left(\nabla_\xi g\right)\right\|^2+\left\|T^*\left(\nabla_\xi g\right)\right\|^2\geq c_0\left\|\nabla_\xi g\right\|^2
$$
(when $\nabla_\xi g$ belongs to the domain of $T^*$ for $g$ in the domain of $T^*$) and to use the commutators $\left[\nabla_\xi,S\right]$ and
$\left[\nabla_\xi,T^*\right]$ to transform the estimate to
$$
\left\|\nabla_\xi\left(Sg\right)\right\|^2+\left\|\nabla_\xi \left(T^*g\right)\right\|^2\geq c_0\left\|\nabla_\xi g\right\|^2-C\left\|g\right\|^2
$$
(where $C$ is a positive number) to finally obtain
$$
\left\|Sg\right\|_\xi^2+\left\|T^*g\right\|_\xi^2\geq\frac{c_0}{2}\left\|g\right\|_\xi^2
$$
for some positive constant $c_1$ when the positive constant $\gamma$ is chosen to be sufficiently small relative to $C$, for example, no more than $\frac{c_0}{2C}$.  When the argument involving the curvature boundary behavior of the metric of the punctured open unit $1$-disk of negative curvature is applied with the vector field, there is an additional order of $\varepsilon$ in the denominator from the cut-off function with the size of its support of the order $\varepsilon$.  In order to get from this argument a section $\tilde s_0$ of $K_{X\times{\mathcal C}}+{\mathcal L}$ over $(X-D)\times\partial\Delta$ with
$$
\int_{(X-D)\times \partial\Delta}\left(\left|\tilde s_0\right|^2+\left|\nabla_\xi\tilde s_0\right|^2\right)h^{(0)}_j\exp\left(2{\rm Re}\left(\bar\tau f_j\right)\right)\leq\hat C,
$$
we need to assume that the given element $s_{\tau_0}\in\Gamma\left(X,K_X+L^{\left(\tau_0\right)}\right)$ is already known to be extendible to an element of
$$\Gamma\left(X\times{\mathbb C},\left({\mathcal O}_{X\times{\mathbb C}}\left/\left({\mathcal I}_{X\times\left\{\tau_0\right\}}\right)^2\right.\right)\left(K_{X\times{\mathbb C}}+{\mathcal L}\right)\right),$$ where ${\mathcal I}_{X\times\left\{\tau_0\right\}}$ is the ideal sheaf of $X\times\left\{\tau_0\right\}$ on
$X\times{\mathbb C}$.  In other words, we need to know first that the section $s_{\tau_0}$ is known to be extendible from over the point $\tau_0$ in ${\mathbb C}$ to over the $1$-jet of ${\mathbb C}$ at $\tau_0$ ({\it i.e.,} over the double point of ${\mathbb C}$ at $\tau_0$).

\medbreak Since a function on a circle which is $L^2$ and whose first derivative is also $L^2$ is continuous on the circle, we now can conclude the following.  For $\tau_0\in\partial\Delta$ every element $s_{\tau_0}\in\Gamma\left(X,K_X+L^{\left(\tau_0\right)}\right)$ can be extended to a continuous section of $K_{X\times{\mathbb C}}+{\mathcal L}$ over $X\times\partial\Delta$ if $s_{\tau_0}$ is assumed to be extendible to an element of
$\Gamma\left(X\times{\mathbb C},\left({\mathcal O}_{X\times{\mathbb C}}\left/\left({\mathcal I}_{X\times\left\{\tau_0\right\}}\right)^2\right.\right)\left(K_{X\times{\mathbb C}}+{\mathcal L}\right)\right)$.

\medbreak Now we put this conclusion together with the following simple statement.  If a subvariety $E$ of ${\mathbb C}^p$ has the property that $E$ contains every complex line in ${\mathbb C}^p$ which is tangential to $E$ at some point of $E$, then $E$ itself must be a (locally finite) union of translates of  ${\mathbb C}$-linear subspaces of the ${\mathbb C}$-vector space ${\mathbb C}^p$.  We obtain the following proposition.

\medbreak\noindent(II.2.4) {\it Proposition (Local Linearity of Moduli of Flat Bundles with Dimension for Space of Canonically Twisted Sections No Less Than Prescribed Number).}  Let $X$ be a compact complex algebraic manifold and $A$ be the abelian variety of all flat line bundles on $X$.  Let $q$ be a positive integer and let $Z_q$ be the set of points $a\in A$ such that $\dim_{\mathbb C}\Gamma\left(X, K_X+F^{(a)}\right)\geq q$ (where $F^{(a)}$ is the flat line bundle on $X$ which corresponds to the point $a$ of $A$).  Then $Z_q$ is regular and is a finite union of translates of abelian subvarieties of the abelian variety $A$.

\bigbreak\noindent(II.2.5) {\it Alternative Description of Relation of Moduli of Flat Bundles and Holomorphic $1$-Forms by Connection.}  In the above derivation of the linearity of moduli of flat bundles characterized by the dimension for the space of canonically twisted sections no less than a prescribed number, we use a somewhat elaborate argument of applying an extension result of Ohsawa-Takegoshi type with an $L^2$ norm involving first-order derivatives.  The main reason for the need of the somewhat elaborate argument is that the natural metric of $\tilde{\mathcal F}$ is defined by the positive function\ \  $\exp\left(2{\rm Re}\sum_{\nu=1}^\ell t_\nu\overline{f_{j,\omega_\nu}}\right)$ on $U_j\times\tilde A$, where the holomorphic function $f_{j,\omega_\nu}$ occurs in the form of its complex-conjugation whereas the coordinates $t_\nu$ of the universal $\tilde A$ of the abelian variety $A$ occur as themselves without complex-conjugation.  The need to use the complex-conjugation of $f_{j,\omega_\nu}$ is explained above abstractly using the Hodge decomposition which gives the isomorphism between $H^1\left(X,{\mathcal O}_X\right)$ and $\Gamma\left(X,\Omega_X^1\right)$.  In order to corroborate the need for the use of complex-conjugation of is indeed correct, we would like to argue from another perspective by directly and explicitly using holomorphic $1$-forms to give different complex structures to a smooth flat line bundle.  This paragraph is included solely for the sake of affirmation and fuller geometric understanding and is not an essential part of the arguments of this note.

\medbreak The transition functions of the flat line bundle $L^{(0)}$ on $X$ which we start out with are constant functions $g_{jk}^{(0)}$ on $U_j\cap U_k$ and the non identically zero holomorphic $1$-form on $X$ which we start out with is $\omega$ so that $\omega=df_j$ on $U_j$ for some holomorphic function $f_j$ on $U_j$.  For a complex variable $\tau\in{\mathbb C}$ we would like to write down the transition functions of the flat line bundle $L^{(\tau)}$ which is obtained by perturbing the connection for the definition of $\bar\partial$ of $L^{(0)}$ by using $\tau\omega$.
There are local holomorphic frames $F_j$ of $L^{(0)}$ on $U_j$ ({\it i.e.,} each each $F_j$ is a holomorphic section of $L^{(0)}$ on $U_j$) with $F_j=g_{jk}^{(0)}F_k$ on $U_j\cap U_k$.
We give $L^{(0)}$ another complex structure to make it into another holomorphic line bundle $L^{(\tau)}$ by defining a local smooth section $s$ of $L^{(\tau)}$ to be holomorphic if and only if $\bar\partial s=\tau s\,\bar\omega$.  In other words, the new differentiation in the $(0,1)$-direction for local smooth sections of $L^{(\tau)}$ is $\bar\partial_\tau=\bar\partial -\tau\,\bar\omega$.  Note that the underlying smooth structure of $L^{(\tau)}$ is still $L^{(0)}$.  We have to use $\bar\omega$ in the definition of $\bar\partial_\tau$ instead of $\omega$ because of type considerations.  The new differentiation $\bar\partial_\tau$ in the $(0,1)$-direction is indeed integrable, because
$$
\displaylines{\bar\partial_\tau\bar\partial_\tau s=\left(\bar\partial -\tau\,\bar\omega\right)
\left(\bar\partial -\tau\,\bar\omega\right)s\cr=\bar\partial\bar\partial s-\tau\,\bar\partial\left(\bar\omega\,s\right)
-\tau\,\bar\omega\wedge\bar\partial s+\tau^2
\bar\omega\wedge\bar\omega\,s\cr=\bar\partial\bar\partial s-\tau\left(\bar\partial\bar\omega\right)s+\tau\,\bar\omega\wedge\bar\partial s
-\tau\,\bar\omega\wedge\bar\partial s+\tau^2
\bar\omega\wedge\bar\omega\,s=0,}
$$
because the holomorphic $1$-form $\omega$ on the compact K\"ahler manifold $X$ must be $d$-closed.  To write down the transition functions for $L^{(\tau)}$ we construct a nowhere zero local holomorphic sections $s_{j,\tau}$ of $L^{(\tau)}$ on $U_j$ by solving the equation $\bar\partial_\tau s_{j,\tau}\equiv 0$ on $U_j$ which is the same as the equation
$\bar\partial s_{j,\tau}=\tau s_{j,\tau}\bar\omega$ on $U_j$.  One natural solution is $s_{j,\tau}=e^{\tau\overline{f_j}}F_j$ on $U_j$ so that the transition functions $\left\{g_{jk,\tau}\right\}_{j,k\in J}$ for $L^{(\tau)}$ with respect to the cover $\left\{U_j\right\}_{j\in J}$ are given by
$$g_{jk,\tau}=\frac{\ e^{\tau\overline{f_j}}F_j\ }{\ e^{\tau\overline{f_k}}F_k}=
e^{\tau\left(\overline{f_j}-\overline{f_k}\right)}g_{jk}^{(0)}=e^{-\tau\overline{c_{jk}}}g_{jk}^{(0)}$$
on $U_j\cap U_k$ which is a constant.  Note that here we cannot replace $\tau$ by $\bar\tau$, because the transition function $g_{jk,\tau}$ must depend holomorphically on $\tau$.

\medbreak If $h^{(0)}_j$ is a positive constant function on $U_j$ with $h^{(0)}_j\left|g_{jk}^{(0)}\right|^2=h^{(0)}_k$ so that $\left\{h_j^{(0)}\right\}$ defines a metric for $L^{(0)}$ with zero curvature, then $h^{(0)}_je^{-2{\rm Re}\left(\tau\overline{f_j}\right)}$ on $U_j$ defines a metric $h_j^{(\tau)}$ on $U_j$ for $L^{(\tau)}$ with zero curvature, because
$$
h_j^{(0)} e^{-2{\rm Re}\left(\tau\overline{f_j}\right)}\left|e^{-\tau\overline{c_{jk}}}g_{jk}\right|^2=h_j^{(0)} e^{-2{\rm Re}\left(\tau\overline{f_j}\right)}\left|e^{\tau\left(\overline{f_j}-\overline{f_k}\right)}g_{jk}\right|^2=
h^{(0)}_k e^{-2{\rm Re}\left(\tau\overline{f_k}\right)}
$$
on $U_j\cap U_k$.  This affirms once more that the product of $\tau$ and $\overline{f_j}$ occurs in the metric with this explicit description without the use of the abstract Hodge decomposition.

\bigbreak\noindent(II.2.6) {\it Direct Explicit Derivation of $L^2$ Estimates Involving First-Order Derivatives.}  We use the s in (II.2.5).  We use the same notation to denote the vector field  $\frac{\partial}{\partial\tau}$ on ${\mathbb C}$ with complex coordinate $\tau$ and for its lifting to the product space $X\times\Delta_R$.  The covariant derivative with respect to vector field $\frac{\partial}{\partial\tau}$ is denoted by $\nabla_\tau$.  The holomorphic line bundle ${\mathcal L}$ on $X\times\Delta_R$ is obtained by putting together the line bundles $L^{(\tau)}$ over $X\times\left\{\tau\right\}$ for $\tau\in\Delta_R$ so that with respect to the covering $\left\{U_j\times\Delta_R\right\}$ the transition function from $U_k\times\Delta_R$ to $U_j\times\Delta_R$ for ${\mathcal L}$ is $g_{jk,\tau}=e^{-\tau\overline{c_{jk}}}g_{jk}$ and the metric for ${\mathcal L}$ on $U_j\times\Delta_R$ is $h_j=h_j^{(0)}e^{-2{\rm Re}\left(\tau\overline{f_j}\right)}$.

\medbreak By using, for a holomorphic local frame $s$ of ${\mathcal L}$, the pair $\left(s,\partial_\tau s\right)$ as a new local frame of a vector bundle of rank $2$, we obtain from ${\mathcal L}$ a holomorphic vector bundle ${\mathcal V}$ of rank $2$ over $X\times\Delta_R$, which is the $1$-jet bundle of ${\mathbb L}$ consisting only of $1$-jets in the direction of the coordinate $\tau$.  Note that this is possible, because $X\times\Delta_R$ is a product space and $1$-jets in the direction of the second factor $\Delta_R$ of the product $X\times\Delta_R$ is well-defined.  Explicitly, the transition functions of ${\mathcal V}$ with respect to the covering $\left\{U_j\times\Delta_R\right\}$ can be described as follows by using the covariant derivative $\nabla_\tau s$ of a local section $s$ of ${\mathcal L}$ and taking advantage of the fact that $\nabla_\tau s$ is well-defined as a local section of ${\mathcal L}$ and its computation by differentiation can be done in any coordinate chart with the result independent of the coordinate chart used.  When a local section $s$ of ${\mathcal L}$ is represented by $s_j$ with respect to the coordinate chart $U_j\times\Delta_R$, we have $$\nabla_\tau s=\partial_\xi s+\left(\partial_\tau\log h_j\right)s_j=\partial_\tau s_j-\left(2\overline{f_j}\right)s_j.$$  The pair $\left(s_j,\partial_\tau s_j\right)$ is related to the pair $\left(s_j,\nabla_\tau s_j\right)$ by
$$
\left(\begin{matrix}s_j\cr\nabla_\tau s_j\end{matrix}\right)
=\left(\begin{matrix}1&0\cr-2\overline{f_j}&1\end{matrix}\right)\left(\begin{matrix}s_j\cr\partial_\tau s_j\end{matrix}\right).\leqno{({\rm II}.2.6.1)}
$$
It follows from $s_j=e^{-\tau c_{jk}}s_k$ and $\nabla_\tau s_j=e^{-\tau c_{jk}}\nabla_\tau s_k$ on $\left(U_j\cap U_k\right)\times\Delta_R$ that the transition function $G_{jk}$ for ${\mathcal V}$ from $U_k\times\Delta_R$ to $U_j\times\Delta_R$ is given by the $2\times 2$ matrix of functions
$$
\displaylines{G_{jk}=e^{-\tau c_{jk}}\left(\begin{matrix}1&0\cr-2\overline{f_j}&1\end{matrix}\right)^{-1}
\left(\begin{matrix}1&0\cr-2\overline{f_k}&1\end{matrix}\right)\cr=e^{-\tau c_{jk}}\left(\begin{matrix}1&0\cr2\overline{f_j}-2\overline{f_k}&1\end{matrix}\right).}
$$
We can use, as the metric $H_j$ for ${\mathcal V}$ on $U_j\times\Delta_R$, the positive-definite $2\times 2$ matrix
$$
\displaylines{H_j=h_j\left(\begin{matrix}1&0\cr-2\overline{f_j}&1\end{matrix}\right)^t
\overline{\left(\begin{matrix}1&0\cr-2\overline{f_j}&1\end{matrix}\right)}\cr
=h_j\left(\begin{matrix}1+4\left|f_j\right|^2&-2\overline{f_j}\cr-2f_j&1\end{matrix}\right).}
$$
Now we can repeat the arguments in (II.2.1) and (II.2.2) with the line bundle ${\mathcal L}$ replaced by the vector bundle ${\mathcal V}$ of rank $2$, because we can add to the metric of ${\mathcal V}$ the metric (or weight function)
$$\left(h_{\tilde D}\right)^{\hat m\eta} e^{-\frac{1}{2\eta}\left(|\tau|^2-1\right)}\left|s_{j,\tilde D}\right|^{2\hat m\eta}
$$
of the trivial line bundle as an additional factor
with some appropriate positive integer $\hat m$ (which is $2$ in the case of ${\mathcal L}$ earlier but needs to be increased for ${\mathcal V}$).   The following additional argument is needed before passing to the limit $\eta\to 0$.

\medbreak The vector ${\mathcal V}$ is simply a $1$-jet bundle obtained from the line bundle ${\mathcal L}$ by using only the jet direction of the vector field $\frac{\partial}{\partial\tau}$.  Introduce the simply connected open subset
$$G_{\eta,\tau_0}=\left\{\,\tau\in{\mathbb C}\,\Big|\,1-\eta<\left|\tau\right|<1+\eta,\,{\rm Re}\left(\tau\left(\tau_0\right)^{-1}\right)>0\right\}$$
of $G_\eta$ which is one half of $G_\eta$ containing $\tau_0$.  Starting from $\tau=\tau_0$ we integrate the holomorphic section of ${\mathcal V}$ on $\left(X-A\right)\times G_\eta$ with $L^2$ estimate (analogous to (II.2.2.1)), with respect to the variable $\tau$ over the simply connected set $G_{\eta,\tau_0}$, to obtain a holomorphic section $\hat s_\eta$ of ${\mathcal L}$ on $\left(X-A\right)\times G_{\eta,\tau_0}$ with $L^2$ estimate (analogous to (II.2.2.1)).  Let $\widehat{\partial\Delta}_{\tau_0}$ be $\partial\Delta\cap\left\{{\rm Re}\left(\tau\left(\tau_0\right)^{-1}\right)>0\right\}$.  By passing to limit as $\eta\to 0$, we obtain from $\hat s_\eta$ a section $\hat s_0$ of ${\mathcal L}$ on $\left(X-A\right)\times\widehat{\partial\Delta}_{\tau_0}$ which is holomorphic along the factor $X-A$ and is only continuous along the factor $\widehat{\partial\Delta}_{\tau_0}$.  From the $L^2$ estimate of $\hat s_0$ which is obtained from that of $\hat s_\eta$, we conclude that $\hat s_0$ can be extended to a section of ${\mathcal L}$ on $X\times\widehat{\partial\Delta}_{\tau_0}$ which is holomorphic along the factor $X$ and is only continuous along the factor $\widehat{\partial\Delta}_{\tau_0}$.  This finishes the additional argument.

\bigbreak\noindent(II.2.6.2) {\it Remark on Integrating Section of Jet Bundle.}  In general the section $\sigma_1$ of a $1$-jet bundle $J_1(E)$ of a line bundle $E$ is the $1$-jet of a section $\sigma$ of $E$. That is why in the last part of the preceding argument we need to start from $\tau=\tau_0$ to integrate the holomorphic section of ${\mathcal V}$ on $\left(X-A\right)\times G_\eta$, with respect to the variable $\tau$ over the simply connected set $G_{\eta,\tau_0}$, to obtain a holomorphic section $\hat s_\eta$ of ${\mathcal L}$ on $\left(X-A\right)\times G_{\eta,\tau_0}$.  It is necessary to perform the integration on a simply connected open subset $G_{\eta,\tau_0}$ of $G_\eta$ in order to get a well-defined single-valued section $\hat s_\eta$ of ${\mathcal L}$ on $\left(X-A\right)\times G_{\eta,\tau_0}$.

\medbreak We would like to remark that, in the context of a holomorphic family of compact complex algebraic manifolds over the open unit $1$-disk and a holomorphic line bundle over the family space, if a holomorphic section on a fiber over a point is extendible to holomorphic sections of the fibers over a circular arc containing the point, even with only continuous dependence on the varying point on the arc, then the section can be extended to a holomorphic section over the entire open unit $1$-disk.  Let us summarize the result.

\bigbreak\noindent(II.2.7) {\it Proposition (Extension of Flatly Twisted Canonical Section Known to Extend to Over Double Point).}  Let $X$ be a compact complex algebraic manifold of complex dimension $n$ and $A$ be the abelian variety of all flat line bundles of $X$ and let ${\mathcal F}$ be the universal flat line bundle over $X\times A$ such that for $a\in A$ the restriction of ${\mathcal F}$ to $X\times\left\{a\right\}$ is the flat line bundle $F^{(a)}$ on $X$ corresponding to the point $a$ of $A$.  Let $\tilde\pi:{\mathbb C}^\ell\to A$ be the universal covering map of $A$.  Let $\tilde a_0$ and $\tilde b$ be distinct points of $\tilde\pi$.  Let $a_0=\tilde\pi\left(\tilde a_0\right)$ and $b=\tilde\pi\left(\tilde b\right)$ and $R>1$ such that the map $\sigma_{a_0,b,R}$ from $\Delta_R$ to $A$ defined by $\tilde\pi\left(\tilde a_0+\tau\left(\tilde b-\tilde a_0\right)\right)$ for $t\in\Delta_R$ is a (holomorphic) embedding of $\Delta_R$ into $A$ with image $C_{a_0,b,R}$.  Let $s_{a_0}$ be an element of
$\Gamma\left(X,K_X+F^{\left(a_0\right)}\right)$ which is naturally identified with
$$\Gamma\left(X\times C_{a_0,b,R}, \left({\mathcal O}_{X\times C_{a_0,b,R}}\left/{\mathcal I}_{X\times\left\{a_0\right\}}\right.\right)\left(K_{X\times C_{a_0,b,R}}+{\mathcal F}\right)\right),
$$
where ${\mathcal I}_{X\times\left\{a_0\right\}}$ is the ideal sheaf for $X\times\left\{a_0\right\}$ on $X\times C_{a_0,b,R}$.  Assume that $s_{a_0}$ can be extended to an element of
$$\Gamma\left(X\times C_{a_0,b,R}, \left({\mathcal O}_{X\times C_{a_0,b,R}}\left/\left({\mathcal I}_{X\times\left\{a_0\right\}}\right.\right)^2\right)\left(K_{X\times C_{a_0,b,R}}+{\mathcal F}\right)\right).
$$
In other words, extension of $s_{a_0}$ is assumed possible from the fiber of the projection $X\times C_{a_0,b,R}\to C_{a_0,b,R}$ at the point $a_0$ of the base manifold $C_{a_0,b,R}$ to over the double point of the base manifold $C_{a_0,b,R}$ at $a_0$.  Then $s_{a_0}$ can be extended to an element of
$$\Gamma\left(X\times C_{a_0,b,R}, {\mathcal O}_{X\times C_{a_0,b,R}}\left(K_{X\times C_{a_0,b,R}}+{\mathcal F}\right)\right).
$$
That is, $s_{a_0}$ can be extended from the fiber of the projection $X\times C_{a_0,b,R}\to C_{a_0,b,R}$ at the point $a_0$ of the base manifold $C_{a_0,b,R}$ to over the entire base manifold $C_{a_0,b,R}$.

\medbreak For the arguments in the proof of Proposition (II.2.7) extension results of Ohsawa-Takegoshi type are applied to only to the special case of a holomorphic family of compact complex algebraic manifolds which is a trivial product.  Such extension results of Ohsawa-Takegoshi type can be applied to a holomorphic family of compact complex manifolds which may not be a trivial product.  When we apply the same arguments to a holomorphic family of compact complex manifolds which may not be a trivial product, we obtain the following more general form of Proposition (II.2.7).

\bigbreak\noindent(II.2.8) {\it Proposition (Extension of Flatly Twisted Canonical Section Known to Extend to Over Double Point for General Holomorphic Family).}  Let $A$ be an abelian variety with universal cover $\tilde\pi:{\mathbb C}^\ell\to A$ and for $a\in A$ let ${\mathcal F_A}^{(a)}$ be the flat line bundle on $A$ which corresponds to the point $a$ of $A$ when $A$ is naturally identified with the moduli space of all flat line bundles over $A$.  Let ${\mathfrak p}$ be a point of $A$ and $\hat A$ be the translate of an $\hat\ell$-dimensional ${\mathbb C}$-linear subspace of $A$ by ${\mathfrak p}$. Let $\hat G$ be a bounded open convex subset of $\hat A$ such that $\tilde\pi$ maps $\hat G$ biholomorphically onto its image $G$ in $A$.
Let $\pi:{\mathcal X}\to G$ be a holomorphic family of compact complex algebraic manifolds of complex dimension $n$ over $G$ with fiber $X_a=\pi^{-1}(a)$ and let ${\mathcal F}$ be a holomorphic line bundle on $X$ such that for $a\in G$ the restriction $F^{(a)}$ of ${\mathcal F}$ to $X_a$ is equal to $\pi^*\left(F_A^{(a)}\right)$.
Let $\tilde a_0$ and $\tilde b$ be distinct points of $\hat G$.  Let $a_0=\tilde\pi\left(\tilde a_0\right)$ and $b=\tilde\pi\left(\tilde b\right)$ and $R>1$ such that the map $\sigma_{a_0,b,R}$ from $\Delta_R$ to $A$ defined by $\tilde\pi\left(\tilde a_0+\tau\left(\tilde b-\tilde a_0\right)\right)$ for $\tau\in\Delta_R$ is a (holomorphic) embedding of $\Delta_R$ into $A$ with image $C_{a_0,b,R}$ contained in $G$.  Let $s_{a_0}$ be an element of
$\Gamma\left(X_{a_0},K_{X_{a_0}}+F^{\left(a_0\right)}\right)$ which is naturally identified with
$$\Gamma\left(\pi^{-1}\left(C_{a_0,b,R}\right), \left({\mathcal O}_{\pi^{-1}\left(C_{a_0,b,R}\right)}\left/{\mathcal I}_{X_{a_0}}\right.\right)\left(K_{\pi^{-1}\left(C_{a_0,b,R}\right)}+{\mathcal F}\right)\right),
$$
where ${\mathcal I}_{X_{a_0}}$ is the ideal sheaf for $X_0$ on $\pi^{-1}\left(C_{a_0,b,R}\right)$.  Assume that $s_{a_0}$ can be extended to an element of
$$\Gamma\left(\pi^{-1}\left(C_{a_0,b,R}\right), \left({\mathcal O}_{\pi^{-1}\left(C_{a_0,b,R}\right)}\left/\left({\mathcal I}_{X_{a_0}}\right.\right)^2\right)\left(K_{\pi^{-1}\left(C_{a_0,b,R}\right)}+{\mathcal F}\right)\right).
$$
In other words, extension of $s_{a_0}$ is assumed possible from the fiber of $$\pi^{-1}\left(C_{a_0,b,R}\right)\to C_{a_0,b,R}$$ at the point $a_0$ of the base manifold $C_{a_0,b,R}$ to over the double point of the base manifold $C_{a_0,b,R}$ at $a_0$.  Then $s_{a_0}$ can be extended to an element of
$$\Gamma\left(\pi^{-1}\left(C_{a_0,b,R}\right), {\mathcal O}_{\pi^{-1}\left(C_{a_0,b,R}\right)}\left(K_{\pi^{-1}\left(C_{a_0,b,R}\right)}+{\mathcal F}\right)\right).
$$
That is, $s_{a_0}$ can be extended from fiber of the projection $$\pi^{-1}\left(C_{a_0,b,R}\right)\to C_{a_0,b,R}$$ at the point $a_0$ of the base manifold $C_{a_0,b,R}$ to over the entire base manifold $C_{a_0,b,R}$.

\bigbreak\noindent(II.2.9) {\it Section Extension Obstruction from Negativity of Curvature.} The extension result given in (II.2.8) for the case of curvature current with controllable negativity immediately raises the question of how to quantitatively describe the obstruction of section extension in terms of the negativity of the curvature current.  Since extendibility can be interpreted as deformability, to imitate the theory of obstruction to deformation, we can formulate quantitatively the obstruction to section extension in terms of an ideal sheaf on the base manifold whose zero-set consists only of the point of the base from whose fiber the given section is to be extended.

\medbreak\noindent(II.2.9.1) {\it Definition (Extendibility Ideal).}  Let $\pi:{\mathcal Y}\to T$ be a holomorphic family of compact complex algebraic manifolds with a Stein manifold $T$ as its base manifold.  Let ${\mathcal E}$ be a holomorphic line bundle on ${\mathcal Y}$.  Let $t_0\in T$ and let ${\mathfrak m}_{T,t_0}$ be the maximum ideal sheaf on $T$ of the point $t_0$.  An ideal sheaf ${\mathcal J}$ on $T$ whose zero-set is the singleton set $\left\{t_0\right\}$ is called an {\it extendibility ideal} at $t_0$ for $\pi:{\mathcal Y}\to T$ and the line bundle ${\mathcal E}$ if the image
$$
\Gamma\left({\mathcal Y}, {\mathcal E}+K_{\mathcal Y}\right)\to \Gamma\left({\mathcal Y}, \left({\mathcal O}_{\mathcal Y}\left/\pi^*\left({\mathfrak m}_{T,t_0}\right)\right)\left({\mathcal E}+K_{\mathcal Y}\right)\right)\right)
$$
contains the image of
$$
\Gamma\left({\mathcal Y}, \left({\mathcal O}_{\mathcal Y}\left/\pi^*\left({\mathcal J}\right)\right)\left({\mathcal E}+K_{\mathcal Y}\right)\right)\right)\to \Gamma\left({\mathcal Y}, \left({\mathcal O}_{\mathcal Y}\left/\pi^*\left({\mathfrak m}_{T,t_0}\right)\right)\left({\mathcal E}+K_{\mathcal Y}\right)\right)\right).
$$

\medbreak\noindent(II.2.9.2) {\it Problem of Quantitative Description of Extension Obstruction by Integral of Negative Part of Curvature Current.}  Suppose in the definition of extendibility ideal in (II.2.9.1) the holomorphic line bundle ${\mathcal E}$ is given a metric whose curvature current is $\Theta_{\mathcal E}$. The problem is to find a formula for a positive number $q$ in terms of some integral expression involving the negative part of $\Theta$ such that there exists an extendibility ideal ${\mathcal J}$ with $\dim_{\mathbb C}{\mathcal O}_T\left/{\mathcal J}\right.\leq q+1$.  For example, the integral expression could be analogous to the integral expressions introduced by Demailly in his theory of holomorphic Morse inequality [Demailly1985].

\medbreak The techniques described above, in (II.2.1) through (II.2.6), of using an additional twisting by an ample line bundle with metrics uniformly bounded on annuli with difference of outer and inner radii approaching zero and of using the circle group action can straightforwardly be adapted to give the following extension theorem.

\bigbreak\noindent(II.2.10) {\it Theorem.}  Let $\pi:X\to\Delta$ be a holomorphic family of compact complex algebraic manifolds over the open unit $1$-disk $\Delta$ whose fiber $X_t=\pi^{-1}(t)$ is of complex dimension $n$.  Let ${\mathcal L}$ be a holomorphic line bundle over $X$ with metric $e^{-\varphi}$ such that 

\medbreak\noindent(i) the curvature current of the metric $e^{-\varphi|_{X_t}}$ of the restriction $L_t$ of ${\mathcal L}$ to $X_t$ is nonnegative for each $t\in\Delta$,

\medbreak\noindent(ii) the connection $\partial\varphi$ of the metric $e^{-\varphi}$ is locally uniformly bounded on $X$, and 

\medbreak\noindent(iii) the curvature current $\sqrt{-1}\,\partial\bar\partial\varphi$ of $e^{-\varphi}$ is $\geq -\omega$ for some smooth positive $(1,1)$-current on $X$.

\medbreak\noindent Let $s_0$ be a holomorphic section of $L_0$ on $X_0$ with $\left|s_0\right\|^2 e^{-\varphi}$ locally integrable on $X_0$ such that $s_0$ can be extended to a section $s_1$ of ${\mathcal L}$ over $\left(X,{\mathcal O}_X\left/\pi^{-1}\left({\mathfrak m}_{\Delta,0}^2\right)^2\right.\right)$, where ${\mathfrak m}_{\Delta,0}$ is the maximal ideal of $\Delta$ at $0$.  If the covariant derivative at $0$ of $s_1$ with respect to the coordinate of $\Delta$ and with respect to the metric $e^{-\varphi}$ of ${\mathcal L}$ is locally $L^2$ on $X_0$, then $s_0$ can be extended to a holomorphic section of ${\mathcal L}$ over any relatively compact open neighborhood of $X_0$ in $X$.

\bigbreak

\bigbreak\noindent{\bf \S3. Local linearity of subvariety defined by minimum number of independent flatly twisted pluricanonical sections.}

\medbreak\noindent(II.3.1) {\it Proposition (Linearity of Moduli of Flat Bundles with Dimension for Space of Pluricanonically Twisted Sections No Less Than Prescribed Number).}  Let $X$ be a compact complex algebraic manifold and $A$ be the abelian variety of all flat line bundles on $X$.  Let $m$ and $q$ be positive integers and let $Z_{m,q}$ be the set of points $a\in A$ such that $\dim_{\mathbb C}\Gamma\left(X, mK_X+F^{(a)}\right)\geq q$ (where $F^{(a)}$ is the flat line bundle on $X$ which corresponds to the point $a$ of $A$).  Then $Z_{m,q}$ is regular and is a finite union of translates of abelian subvarieties of the abelian variety $A$.

\medbreak\noindent{\it Proof.}  We let $\hat Z_{m,q}$ be the set of points $a$ of $A$ such that the complex dimension of $\Gamma\left(X, mK_X+mF^{(a)}\right)$ is $\geq q$.  It suffices to show that $\hat Z_{m,q}$ is regular and is a finite union of translates of abelian subvarieties of the abelian variety $A$.  Let $A\approx{\mathbb C}^\ell$ be the universal of $A$ with the universal covering map $\tilde\pi:\tilde A\to A$.  Let $C$ be a local regular curve inside $\hat Z_{m,q}$ and let $a_0$ be a point of $C$.  Let $\hat A$ be a translate of a linear complex line in $\tilde A$ such that $\tilde\pi$ maps some point $\tilde a_0$ of $\hat A$ to $a_0$ and the image $\tilde\pi\left(\hat A\right)$ of $\hat A$ in $A$ is the tangent line of the local complex curve $C$ at $a_0$.  Choose $R>1$ and a point $\tilde b$ of $\hat A$ distinct from $\tilde a_0$ such that the map $\sigma_{\tilde a_0,\tilde b,R}$ from $\Delta_R$ to $A$ defined by $\tilde\pi\left(\tilde a_0+\tau\left(\tilde b-\tilde a_0\right)\right)$ for $\tau\in\Delta_R$ is a (holomorphic) embedding of $\Delta_R$ into $A$ with image $\tilde\Delta_R$ and there is a biholomorphic map $\Phi:C\to\tilde\Delta_R$ (possibly after we replace $C$ by some open neighborhood of $a_0$ in $C$) which fixes $a_0$ and also fixes the tangent direction of $C$ at $a_0$.  For our proof it suffices to show that $\tilde\Delta_R$ is contained in $\hat Z_{m,q}$ for any local regular curve $C$ inside $\hat Z_{m,q}$ and for any point $a_0$ of $C$.

\medbreak There exist $q$ linearly independent elements $$s_1(t),\cdots,s_q(t)\in\Gamma\left(X,mK_X+mF^{(t)}\right)$$ for each $t\in C$ (possibly after we replace $C$ by some open neighborhood of $a_0$ in $C$) which depends on $t$ holomorphically.  We take the graph of $\left(s_1(t)\right)^{\frac{1}{m}}$ and then resolve its singularities.  For clarity of presentation first let us assume that we can do this simultaneous for all $t\in C$ so that
\begin{itemize}
\item[(i)] for every $t\in C$ we have a compact complex algebraic manifold $X_t^\natural$ and a holomorphic map $\pi_t^\natural: X_t^\natural\to X$ which outside a proper subvariety of $X_t^\natural$ is an $m$-sheeted analytic cover;
\item[(ii)] $\left(s_1(t)\right)^{\frac{1}{m}}$ defines an element $s_1^\natural(t)$ of $\Gamma\left(X^\natural_t, K_{X^\natural_t}+\left(\pi^\natural_t\right)^*\left(F^{(t)}\right)\right)$;
\item[(iii)] the (multiplicatively defined) direct image $\left(\pi^\natural_t\right)_*\left(s_1^\natural(t)\right)$ (in the sense that the value of $\left(\pi^\natural_t\right)_*\left(s_1^\natural(t)\right)$ at a generic point of $X$ is the product of the values of $s_1^\natural(t)$ at the $m$ distinct points of $X_t^\natural$ in the fiber of that generic point) is equal to $s_1(t)$;
\item[(iv)] the holomorphic maps $\pi_t^\natural: X_t^\natural\to X$ for $t\in C$ can be put together to form a proper holomorphic map $\tilde\pi^\natural$ from a complex manifold ${\mathcal X}^\natural$ to $X\times C$.
\end{itemize}
We denote by $F^{\natural(t)}$ the pullback $\left(\pi^\natural\right)^*\left(F^{(\Phi(t))}\right)$ to $X_t^\natural$ of the line bundle $F^{(\Phi(t))}$ on $X$ under $\pi_t^\natural: X_t^\natural\to X$ and we put together $F^{\natural(t)}$ to form a line bundle ${\mathcal F}^\natural$ on ${\mathcal X}^\natural$.  Let $\pi^\natural:{\mathcal X}^\natural\to\tilde\Delta_R$ be the composite of the map $\tilde\pi^\natural: {\mathcal X}^\natural\to X\times C$ and the natural projection $X\times C\to C$ onto the second factor and the map $\Phi:C\to\tilde\Delta_R$.

 \medbreak We now apply Proposition (II.2.8) to the holomorphic family $\pi^\natural: {\mathcal X}^\natural\to X\times\tilde\Delta_R$ and the line bundle ${\mathcal F}^\natural$ on ${\mathcal X}^\natural$ and the element $s_1^\natural\left(a_0\right)$ of $\Gamma\left(X^\natural_{a_0}, K_{X^\natural_t}+F^{\natural\left(a_0\right)}\right)$ which we seek to extend to an element $\tilde s_1^\natural$ of
 $\Gamma\left({\mathcal X}^\natural, K_{{\mathcal X}^\natural}+{\mathcal F}^\natural\right)$.

 \medbreak Since we have a holomorphic family of $s_1^\natural(t)\in\Gamma\left(X^\natural_t, K_{X^\natural_t}+\left(\pi^\natural_t\right)^*\left(F^{(t)}\right)\right)$ for $t\in C$ and since $\tilde\Delta_R$ is tangential to $C$ at $a_0$, it follows that the element $s_1^\natural\left(a_0\right)$ of $\Gamma\left(X^\natural_{a_0}, K_{X^\natural_t}+F^{\natural\left(a_0\right)}\right)$ can be extended to an element of $\Gamma\left({\mathcal X}^\natural, \left({\mathcal O}_{{\mathcal X}^\natural}\left/\left({\mathcal I}_{X_{a_0}^{\natural}}\right)^2\right.\right)
 \left(K_{{\mathcal X}^\natural}+{\mathcal F}^\natural\right)\right)$ over the double point of the base $\tilde\Delta_R$ at $a_0$, where ${\mathcal I}_{X_{a_0}^{\natural}}$ is the ideal sheaf of $X_{a_0}^\natural$ on ${\mathcal X}^\natural$.  By Proposition (II.2.8) we conclude that $s_1^\natural\left(a_0\right)$ can be extended to an element $\tilde s_1^\natural$ of
 $\Gamma\left({\mathcal X}^\natural, K_{{\mathcal X}^\natural}+{\mathcal F}^\natural\right)$.

 \medbreak For $t\in C$ let $\tilde s_1^\natural(t)$ be the element of $\Gamma\left(X^\natural_t, K_{X_t^\natural}+F^{\natural(t)}\right)$ defined by the element $\tilde s_1^\natural$ of
 $\Gamma\left({\mathcal X}^\natural, K_{{\mathcal X}^\natural}+{\mathcal F}^\natural\right)$.
 We have to construct from the element $\tilde s_1^\natural(t)$ of $\Gamma\left(X^\natural_t, K_{X_t^\natural}+F^{\natural\left(t\right)}\right)$ an element $\hat s_1^\natural(t)$ of $\Gamma\left(X,mK_X+mF^{(t)}\right)$.  We divide the element $\tilde s_1^\natural(t)$ of $\Gamma\left(X^\natural_t, K_{X_t^\natural}+F^{\natural(t)}\right)$
 by the element $s_1^\natural(t)$ of $\Gamma\left(X^\natural_t, K_{X^\natural_t}+\left(\pi^\natural_t\right)^*\left(F^{(t)}\right)\right)$ to form a meromorphic section $\frac{\tilde s_1^\natural(t)}{s_1^\natural(t)}$ of the flat line bundle $$F^{\natural(t)}-\left(\pi^\natural_t\right)^*\left(F^{(t)}\right)=
 \left(\pi^\natural_t\right)^*\left(F^{\left(\Phi(t)\right)}-F^{(t)}\right)$$ on $X^\natural_t$.  We now construct the (multiplicatively defined) direct image of $\frac{\tilde s_1^\natural(t)}{s_1^\natural(t)}$ under $\pi^\natural_t$ to get a meromorphic section $s_1^\flat$ of the flat line bundle $m\left(F^{\left(\Phi(t)\right)}-F^{(t)}\right)$ on $X$ in the sense that the value of $s_1^\flat$ at a generic point $P$ of $X$ is the product of the values of $\frac{\tilde s_1^\natural(t)}{s_1^\natural(t)}$ at the $m$ distinct points of $\left(\pi_t^\natural\right)^{-1}(P)$ in $X_t^\natural$.  Now we can define the element $\hat s_1^\natural(t)$ of $\Gamma\left(X,mK_X+mF^{\left(\Phi(t)\right)}\right)$ to be $s_1(t)\,s_1^{(t)}$.

 \medbreak We can repeat the same argument with $s_1(t)$ replaced by $s_j(t)$ for $2\leq j\leq q$ to conclude that we have a family of ${\mathbb C}$-linearly independent elements
 $$
 \hat s_1^\natural(t),\cdots,\hat s_1^\natural(t)\in\Gamma\left(X,mK_X+mF^{\left(\Phi(t)\right)}\right)
 $$
 holomorphic in $t$ for $t\in C$ sufficiently close to $a_0$.  This means that points of $\tilde\Delta_R$ which are sufficiently close to $a_0$ are contained in $\hat Z_{m,q}$.  This implies that $\tilde\Delta_R$ is contained in $\hat Z_{m,q}$.

\medbreak In an early part of the proof we have assumed that we can simultaneously resolve the singularities of the graph of $\left(s_1(t)\right)^{\frac{1}{m}}$ for all $t\in C$.  This assumption is not essential, because for the use of extension results of Ohsawa-Takegoshi type, instead of a compact complex algebraic manifold, we can replace it by the complement of an ample complex hypersurface which contains the set of singular points.  Q.E.D.

\bigbreak

\bigbreak\noindent{\bf\S4. Two descriptions of moduli, one by projective embedding and one by transition functions.}

\bigbreak\noindent(II.4.1) Let $X$ be a compact complex algebraic manifold of complex dimension $n$ and $A$ be the abelian variety of all flat line bundles of $X$, representing the moduli of flat line bundles on $X$.  There are two ways to describe the moduli of flat line bundles on $A$, one by projectively embedding $A$ into some complex projective space and the other by using the transition functions of the flat line bundles.

\medbreak Let $\Phi:A\to{\mathbb P}_N$ be an embedding of the abelian variety $A$ into a complex projective space ${\mathbb P}_N$.  The first kind of moduli for flat line bundles over $X$ (described by projective embedding) is given by the pullbacks to $A$ under $\Phi$ of meromorphic functions on ${\mathbb P}_N$.  When $X$ is defined over the algebraic closure $\bar{\mathbb Q}$ of ${\mathbb Q}$, $A$ is also defined over $\bar{\mathbb Q}$ and the map $\Phi:A\to{\mathbb P}_N$ can be chosen to be defined over $\bar{\mathbb Q}$.  The image $\Phi\left(Z_q\right)$ (with $Z_q$ from (II.2.4)) in ${\mathbb P}_N$ is defined over $\bar{\mathbb Q}$.

\medbreak Another way to describe the moduli of flat line bundles on $X$ is by using the transition functions of the line bundles with respect to some covering $\left\{U_j\right\}$ of $X$.  The transition functions can be obtained by integrating holomorphic $1$ forms on $X$.  In this sense the moduli of transition functions are related to the moduli from a projective embedding by integration.  Because of the relation by integration, the two descriptions of the moduli of flat line bundles on $X$ are transcendentally related rather than algebraically related.  One essential component of the proof of Theorem (II.1) is to use the very strong restriction on a subvariety of $A$ which are algebraic in both descriptions of the moduli of flat line bundles on $X$.  Such a very strong restriction comes from the arithmetic argument of Gelfond-Schneider [Gelfond1934, Schneider1934], Lang [Lang1962, Lang1965, Lang1966], Bombieri [Bombieri1970, Bombieri-Lang1970], Briskorn [Brieskorn1970], and Simpson [Simpson1993].  We now discuss the description of the moduli of flat line bundles on $X$ by transition functions.

\bigbreak\noindent(II.4.2) {\it Transition Functions as Local Moduli of Flat Bundles and the Computation of Cohomology Groups of Local Systems.}  On our compact complex algebraic manifold $X$ of complex dimension $n$, every flat line bundle is represented by an element of $${\rm Hom}\left(\pi_1(X), {\mathbb C}^*\right)={\rm Hom}\left(H_1\left(X,{\mathbb Z}\right), {\mathbb C}^*\right).$$  The finitely generated abelian group $H_1\left(X,{\mathbb Z}\right)$ is a direct sum of cyclic groups whose generators are $e_1,\cdots,e_{m_1}$.  Let $F$ be a flat line bundle over $X$ which is represented by the element of ${\rm Hom}\left(H_1\left(X,{\mathbb Z}\right), {\mathbb C}^*\right)$ which maps $e_j$ to some $\gamma_j\in{\mathbb C}^*$.  We use the notation $H_{\rm loc\ sys}^p\left(X,F\right)$ to denote the cohomology group when $F$ is regarded as a local system, in contrast to $H^p\left(X,{\mathcal O}_X\left(F\right)\right)$ when the sheaf ${\mathcal O}_X(F)$ of germs of holomorphic sections of $F$ is used for its computation.  When the meaning is clear from the context so that there is no room for misunderstanding, we also use $H^p\left(X,F\right)$ to denote $H^p\left(X,{\mathcal O}_X\left(F\right)\right)$.

\medbreak Let ${\mathcal U}=\left\{U_j\right\}_{1\leq j\leq\ell_0}$ be a finite open cover of $X$ such that
\begin{itemize}
\item[(i)]
for any nonnegative integer $k$ each nonempty intersection $U_{j_0}\cap\cdots\cup U_{j_\nu}$ with $1\leq j_0<\cdots<j_\nu\leq\ell_0$ is simply connected and Stein,
\item[(ii)] the transition function $g_{jk}$ for the flat bundle $F$ from $U_k$ to $U_j$ is a monomial $\gamma_1^{a_{jk,1}}\cdots\gamma_{m_1}^{a_{jk,m_1}}$ with $a_{jk,1},\cdots,a_{jk,m_1}$ being integers.
\end{itemize}
We use the finite open cover ${\mathcal U}=\left\{U_j\right\}_{1\leq j\leq\ell_0}$ to compute $H_{\rm loc\ sys}^p\left(X,F\right)$ for $1\leq p\leq 2n-1$.  Let $J_\nu$ be the set of all $(\nu+1)$-tuple $\left(j_0,\cdots,j_\nu\right)$ with strictly increasing components $1\leq j_0<\cdots<j_\nu\leq\ell_0$ such that $U_{j_0}\cap\cdots U_{j_\nu}$ is nonempty.  The $\nu$-th group of cochains $C^\nu\left({\mathcal U},L\right)$ is given by
$$
C^\nu\left({\mathcal U},L\right)=\bigoplus_{\left(j_0,\cdots,j_\nu\right)\in J_\nu}{\mathbb C}_{j_0,\cdots,j_\nu},
$$
where ${\mathbb C}_{j_0,\cdots,j_\nu}$ is just ${\mathbb C}$ but indexed by $\left(j_0,\cdots,j_\nu\right)$ in order to distinguish the different copies of ${\mathbb C}$ used as direct summands of $C^\nu\left({\mathcal U},L\right)$.  Let
$$\left\{c_{j_0,\cdots,j_\nu}\right\}_{\left(j_0,\cdots,j_\nu\right)\in J_\nu}$$
be an element of $C^\nu\left({\mathcal U},L\right)$ and let
$$\left\{d_{j_0,\cdots,j_{\nu+1}}\right\}_{\left(j_0,\cdots,j_{\nu+1}\right)\in J_{\nu+1}}\in C^{\nu+1}\left({\mathcal U},L\right)$$
be its image under the coboundary map $\delta: C^\nu\left({\mathcal U},L\right)\to C^{\nu+1}\left({\mathcal U},L\right)$.  Then
$$
d_{j_0,\cdots,j_{\nu+1}}=g_{j_0 j_1}c_{j_1,\cdots,j_{\nu+1}}+\sum_{\lambda=1}^{k+1}(-1)^\lambda c_{j_0,\cdots,j_{\lambda-1},j_{\lambda+1},\cdots,j_{\nu+1}}.
$$
The cohomology group $H_{\rm loc\ sys}^p\left(X,F\right)$ is given by the quotient group
$$
\frac{\ {\rm Ker}\left(C^p\left({\mathcal U},L\right)\stackrel{\delta}{\longrightarrow} C^{p+1}\left({\mathcal U},L\right)\right)\ }
{{\rm Im}\left(C^{p-1}\left({\mathcal U},L\right)\stackrel{\delta}{\longrightarrow} C^p\left({\mathcal U},L\right)\right)}.
$$
We denote by $I_\nu$ the number of elements in the finite set $J_\nu$.  Since the space $C^\nu\left({\mathcal U},L\right)$ can be regarded as the ${\mathbb C}$-vector space ${\mathbb C}^{I_\nu}$ whose elements are column $I_\nu$-vectors with components $c_{j_0,\cdots,j_\nu}$ for $\left(j_0,\cdots,j_\nu\right)\in J_\nu$, the ${\mathbb C}$-linear map $\delta: C^\nu\left({\mathcal U},L\right)\to C^{\nu+1}\left({\mathcal U},L\right)$ can be represented by an $I_{\nu+1}\times I_\nu$ matrix ${\mathcal A}_{F,\nu}$ whose nonzero entries are either $\pm 1$ or $\pm \gamma_1^{a_1}\cdots\gamma_{m_1}^{a_{m_1}}$ with $a_1,\cdots,a_{m_1}$ being integers.

\bigbreak

\bigbreak\noindent{\bf\S5. Algebraicity of subvariety in moduli described by transition functions.}

\medbreak We now vary the flat line bundle $F$ by varying $\gamma_1,\cdots,\gamma_{m_1}$ in the following way.  We assume that the cyclic subgroup of $H_1\left(X,{\mathbb Z}\right)$ generated by $e_j$ is infinite for $1\leq j\leq m_0$ and the cyclic subgroup of $H_1\left(X,{\mathbb Z}\right)$ generated by $e_j$ is finite for $m_0+1\leq j\leq m_1$.   Note that each $\gamma_j$ automatically is a root of unity for $m_0+1\leq j\leq m_1$, because the cyclic subgroup of $H_1\left(X,{\mathbb Z}\right)$ generated by $e_j$ is finite for $m_0+1\leq j\leq m_1$.  We vary $\gamma_1,\cdots,\gamma_{m_0}$ with $\gamma_j$ fixed for $m_0+1\leq j\leq m_1$ so that the point $\gamma=\left(\gamma_1,\cdots,\gamma_{m_0}\right)$ varies as a point of $\left({\mathbb C}^*\right)^{m_0}$.  We denote the line bundle $F$ defined by the point $\gamma\in\left({\mathbb C}^*\right)^{m_0}$ by $F_\gamma$.  We fix a point $\gamma^{(0)}\in\left({\mathbb C}^*\right)^{m_0}$.  Assume that $H^p\left(X,F_{\gamma^{(0)}}\right)$ is nonzero.  We have the following simple lemma.

\medbreak\noindent(II.5.1) {\it Lemma.}  Let $Z$ be the set of $\gamma\in\left({\mathbb C}^*\right)^{m_0}$ such that  $\dim_{\mathbb C}H_{\rm loc\ sys}^p\left({\mathcal U}, F_\gamma\right)$ is no less than $\dim_{\mathbb C}H_{\rm loc\ sys}^p\left({\mathcal U}, F_{\gamma^{(0)}}\right)$.  Then there exists some open neighborhood $U$ of $\gamma^{(0)}$ in $\left({\mathbb C}^*\right)^{m_0}$ and there exist a finite number of polynomials ${\mathcal P}_j\left(\gamma_1,\cdots,\gamma_{m_0}\right)$  (for $1\leq j\leq N$) in the variables $\gamma_1,\cdots,\gamma_{m_0}$ with coefficients in a number field generated by ${\mathbb Q}$ and a finite number of roots of unity such that their common zero-set in $U$ is equal to $Z\cap U$.

\medbreak\noindent{\it Proof.}  For $\nu=p-1$ or $p$ let $q_\nu$ be the rank of the $I_{\nu+1}\times I_\nu$ matrix ${\mathcal A}_{F_{\gamma^{(0)}},\nu}$ and let $M_\nu(\gamma)$ be the determinant of a $q_\nu\times q_\nu$ sub-matrix of ${\mathcal A}_{F_\gamma,\nu}$ for $\gamma\in\left({\mathbb C}^*\right)^{m_0}$ such that $M_\nu\left(\gamma^{(0)}\right)$ is nonzero.  Choose an open neighborhood $U$ of $\gamma^{(0)}$ in $\left({\mathbb C}^*\right)^{m_0}$ such that both $M_p(\gamma)$ and $M_{p-1}(\gamma)$ are nonzero for $\gamma\in U$.  The set of polynomials ${\mathcal P}_j\left(\gamma_1,\cdots,\gamma_{m_1}\right)$  (for $1\leq j\leq N$) are obtained by setting equal to $0$ the determinants of all $q_\nu\times q_\nu$ sub-matrix of ${\mathcal A}_{F_\gamma,\nu}$ for $q=p-1,p$.  The conditions on the coefficients of
${\mathcal P}_j\left(\gamma_1,\cdots,\gamma_{m_1}\right)$  (for $1\leq j\leq N$) are fulfilled, because the entries of $I_{\nu+1}\times I_\nu$ matrix ${\mathcal A}_{F_\gamma,\nu}$ whose nonzero entries are either $\pm 1$ or $\pm \gamma_1^{a_1}\cdots\gamma_{m_1}^{a_{m_1}}$ with $a_1,\cdots,a_{m_1}$ being integers and $\gamma_j$ being a root of unity for $m_0+1\leq j\leq m_1$.

\medbreak\noindent(II.5.2) {\it Corollary.}  Let $\gamma^{(0)}\in\left({\mathbb C}^*\right)^{m_0}$.  Let $q$ be $\dim_{\mathbb C}\Gamma\left(X,K_X+F_{\gamma^{(0)}}\right)$ and let $Z_q$ is the set of all $\gamma\in\left({\mathbb C}^*\right)^{m_0}$ such that $\dim_{\mathbb C}\Gamma\left(X,K_X+F_\gamma\right)\geq q$.  Then there exists some open neighborhood $U$ of $\gamma^{(0)}$ in $\left({\mathbb C}^*\right)^{m_0}$ and there exist a finite number of polynomials ${\mathcal P}_j\left(\gamma_1,\cdots,\gamma_{m_0}\right)$  (for $1\leq j\leq N$) in the variables $\gamma_1,\cdots,\gamma_{m_0}$ with coefficients in the algebraic closure $\overline{\mathbb Q}$ of ${\mathbb Q}$ such that their common zero-set in $U$ is equal to $Z_q\cap U$.

\medbreak\noindent{\it Proof.} For $0\leq j\leq n$ let $Z^{(j)}$ be the subvariety of all $\gamma\in\left({\mathbb C}^*\right)^{m_0}$ such that $\dim_{\mathbb C}H^j\left(X,\left(\wedge^{n-j}T^*_X\right)\otimes F_\gamma\right)$ is no less than
$\dim_{\mathbb C}H^j\left(X,\left(\wedge^{n-j}T^*_X\right)\otimes F_{\gamma^{(0)}}\right)$ so that $Z_q=Z^{(0)}$ (where $\wedge^{n-j}T^*_X$ is the exterior product of $n-j$ copies of the dual vector bundle $T_X^*$ of the tangent vector bundle $T_X$ of $X$).  Recall that $Z$ is the set of $\gamma\in\left({\mathbb C}^*\right)^{m_0}$ such that  $\dim_{\mathbb C}H_{\rm loc\ sys}^p\left({\mathcal U}, F_\gamma\right)$ is no less than $\dim_{\mathbb C}H_{\rm loc\ sys}^p\left({\mathcal U}, F_{\gamma^{(0)}}\right)$.  From Hodge decomposition with twisting by a flat line bundle it follows that
$$
H_{\rm loc\ sys}^n\left(X,F_\gamma\right)=\bigoplus_{j=0}^n H^j\left(X,\left(\wedge^{n-j}T^*_X\right)\otimes F_\gamma\right).\leqno{({\rm II}.5.2.1)}
$$
Hence $\cap_{j=1}^n Z^{(j)}\subset Z$.  We claim that $Z^{(0)}$ is a branch of $Z$.  Suppose the contrary and we are going to derive a contradiction.  Then $Z^{(0)}$ is contained in some irreducible branch germ $Z^\prime$ of $Z$ as a proper subvariety-germ of $Z^\prime$, because the subvariety-germ of $Z^{(0)}$ at $\gamma^{(0)}$ is irreducible due to the fact that the set in $A$ corresponding to $Z^{(0)}$ is regular and is a finite union of translates of abelian subvarieties of $A$.  For $1\leq j\leq n$ let $\tilde Z^{(j)}$ be the set of all $\gamma\in\left({\mathbb C}^*\right)^{m_0}$ such that $\gamma\in Z^\prime$ and $\dim_{\mathbb C}H^j\left(X,\left(\wedge^{n-j}T^*_X\right)\otimes F_\gamma\right)$ is strictly greater than $\dim_{\mathbb C}H^j\left(X,\left(\wedge^{n-j}T^*_X\right)\otimes F_{\gamma^{(0)}}\right)$.
From (II.2.5.2.1) it follows that $Z^\prime-Z^{(0)}$ is contained $\cup_{j=1}^n\tilde Z^{(j)}$, which contradicts the fact that each $Z^{(j)}$ must be a proper subvariety of $Z^\prime$ for $1\leq j\leq n$.  This finishes the verification of the claim that $Z^{(0)}$ is a branch of $Z$.  According to Lemma (II.5.1) the subvariety $Z$ is defined by a finite number of polynomials in the variables $\gamma_1,\cdots,\gamma_{m_0}$ with coefficients in $\overline{\mathbb Q}$.  It follows that the branch-germ $Z^{(0)}$ of $Z$ is
defined by a finite number of polynomials in the variables $\gamma_1,\cdots,\gamma_{m_0}$ with coefficients in $\overline{\mathbb Q}$.  Q.E.D.

\bigbreak\noindent(II.5.3) {\it Relation Between Two Kinds of Moduli Descriptions.}  Let ${\rm Alb}_X$ be the Albanese of $X$ and ${\rm alb}_X:X\to{\rm Alb}_X$ be the Albanese map so that every holomorphic $1$-form on $X$ is the pullback of a holomorphic $1$-form on ${\rm Alb}_X$.  Let $\ell$ be the complex dimension of ${\rm Alb}_X$.  Denote by $u_1,\cdots,u_\ell$ the standard basis of the ${\mathbb C}$-vector space ${\mathbb C}^\ell$ ({\it i.e.,} all the components of the column $\ell$-vector $u_j$ are zero except the $j$-th component which is $1$).
There exist a symmetric $\ell\times\ell$ matrix $Z$ with positive-definite imaginary part ${\rm Im}\,Z$ and positive integers $m_1,\cdots,m_\ell$ such that the abelian variety ${\rm Alb}_X$ is the quotient of ${\mathbb C}^\ell$ by the lattice $\Lambda$ of rank $2\ell$ generated by $u_1,\cdots,u_\ell, m_1Zu_1,\cdots,m_\ell Zu_\ell$ (see {\it e.g.,} [Weil1958, p.116, Proposition 6]).  For the variable $\zeta\in{\mathbb C}^\ell$ (regarded as a column $\ell$-vector)  we consider the standard theta function
$$
\Theta(\zeta) = \sum  _{\lambda\in {\mathbb Z}^g}
\exp\left(\pi\sqrt{-1}\left(\lambda^\prime Z\lambda + 2\lambda^\prime
\zeta\right)\right)
$$
for the symmetric matrix $Z$ with positive definite imaginary part, where the multi-index $\lambda$ is a varying column $\ell$-vector with integer components and $\lambda^\prime$ denotes its transpose as an $1\times\ell$ matrix.  The effects of
translations $\zeta\mapsto\zeta+\lambda$ and $\zeta\mapsto\zeta+Z\lambda$ for $\lambda\in{\mathbb Z}^\ell$ are given by $\Theta(\zeta+\lambda) = \Theta(\zeta)$ and  $$\Theta (\zeta + Z\lambda) =
\exp\left(\pi\sqrt{-1}\left(-\lambda^\prime Z\lambda - 2\lambda^\prime
\zeta\right)\right)\Theta (\zeta).$$
For two points $v^{(1)}, v^{(2)}$ of ${\mathbb C}^\ell$ let $\Theta_{v^{(1)},v^{(2)}}(\zeta)$ be the entire function $$\Theta\left(\zeta-v^{(1)}\right)\Theta\left(\zeta-v^{(2)}\right)\Theta\left(\zeta+v^{(1)}+v^{(2)}\right)$$ on ${\mathbb C}^\ell$. For any two pairs of points $\left(v^{(1)}, v^{(2)}\right)$ and $\left(\tilde v^{(1)}, \tilde v^{(2)}\right)$ of ${\mathbb C}^\ell$ the quotient $\Theta_{v^{(1)},v^{(2)}}(\zeta)\left/\Theta_{\tilde v^{(1)},\tilde v^{(2)}}(\zeta)\right.$ is a meromorphic function on ${\rm Alb}_X$.  Let $W_{v^{(1)},v^{(2)}}$ be the divisor in ${\rm Alb}_X$ defined by the function $\Theta_{v^{(1)},v^{(2)}}(\zeta)$.

\medbreak By the theorem of Lefschetz [Weil1958, p.130, Theorem 5] we can find a finite number of pairs of points $\left(v^{(1,j)}, v^{(2,j)}\right)$ (for $0\leq j\leq J_*$) of ${\mathbb C}^\ell$ such that the map ${\rm Alb}_X\to{\mathbb P}_{J_*}$ with
$\Theta_{v^{(1,j)}, v^{(2,j)}}$ ($0\leq j\leq J_*$) as its homogeneous components is an embedding.  In other words, the divisor $W_{v^{(1)},v^{(2)}}$ is a very ample divisor in ${\rm Alb}_X$.  We can choose two pairs of points $\left(b^{(1)}, b^{(2)}\right)$ and $\left(c^{(1)}, c^{(2)}\right)$ of ${\mathbb C}^\ell$ such that for $a\in{\mathbb C}^\ell$ the line bundle $L^{(a)}$ on ${\rm Alb}_X$ defined by $a$ is the same as the line bundle on ${\rm Alb}_X$ associated to the divisor $a+W_{b^{(1)},b^{(2)}}-W_{c^{(1)},c^{(2)}}$.  Note that when we take a point $a\in{\mathbb C}^\ell$ to form the line bundle $L^{(a)}$ on ${\rm Alb}_X$, we are regarding ${\mathbb C}^\ell$ as the universal cover $\tilde A$ of the abelian variety $A$ which is the set of all flat line bundles on $X$ and we are not regarding ${\mathbb C}^\ell$ as the universal cover of the Albanese ${\rm Alb}_X$.  The line bundles $L^{(a)}$ and $L^{\left(a^\prime\right)}$ on $X$ are the same if the two points $a$ and $a^\prime$ of ${\mathbb C}^\ell=\tilde A$ have the same image in $A$ under the quotient map ${\mathbb C}^\ell=\tilde A\to A$ which is the universal cover map.

\medbreak We are going to construct a finite open cover of ${\rm Alb}_X$ and to find the transition functions of $L^{(a)}$ with respect to the open cover.
For $\sigma\in{\mathbb C}^\ell$ let $\tilde U^{(\sigma)}$ be the open parallotope
$$\left\{\,\zeta\in{\mathbb C}^\ell\,\bigg|\,\zeta=\sigma+\sum_{j=1}^\ell\left(\lambda_j u_j+\lambda_{n+j}m_jZu_j\right)\ {\rm for\ some\ }0<\lambda_1,\cdots,\lambda_{2\ell}<\frac{1}{4}\,\right\}.
$$ in ${\mathbb C}^\ell$.
Let $U^{(\sigma)}$ be the open subset of ${\rm Alb}_X$ which is the image of the open subset $\tilde U^{(\sigma)}$ of ${\mathbb C}^\ell$ under the quotient map ${\mathbb C}^\ell\to{\mathbb C}^\ell\left/\Lambda\right.$.  Let $f^{(\sigma)}$ be the meromorphic function on $U^{(\sigma)}$ which is equal to the meromorphic function
$\Theta_{a+b^{(1)},b^{(2)}}(\zeta)\left/\Theta_{\tilde c^{(1)},\tilde c^{(2)}}(\zeta)\right.$ on $\tilde U^{(\sigma)}$ under the biholomorphic map between $\tilde U^{(\sigma)}$ and $U^{(\sigma)}$ defined by the quotient map ${\mathbb C}^\ell\to{\mathbb C}^\ell\left/\Lambda\right.$.  The divisor of $f^{(\sigma)}$ on $U^{(\sigma)}$ is equal to the restriction of the divisor $a+W_{b^{(1)},b^{(2)}}-W_{c^{(1)},c^{(2)}}$ to $U^{(\sigma)}$.  We select a finite number of points
$\sigma^{(j)}$ of ${\mathbb C}^\ell$ (for $1\leq j\leq I$) such that ${\rm Alb}_X$ is covered by $\left\{U^{\left(\sigma^{(j)}\right)}\right\}_{j=1}^I$.

\medbreak The transition functions for the line bundle $L^{(a)}$ on ${\rm Alb}_X$ with respect to $U^{\left(\sigma^{(j)}\right)}$ are given as follows.  For $1\leq j\not=k\leq I$ there exists some $\lambda_{j,k,1},\cdots,\lambda_{j,k,2\ell}\in{\mathbb Z}$ such that the intersection of
$$\tilde U^{\left(\sigma^{(j)}\right)}+\sum_{i=1}^\ell\left(\lambda_{j,k,i}u_i+\lambda_{j,k,\ell+i}m_iZu_i\right)$$
and $\tilde U^{\left(\sigma^{(k)}\right)}$ is mapped bijectively onto $U^{\left(\sigma^{(j)}\right)}\cap U^{\left(\sigma^{(k)}\right)}$ under the
quotient map ${\mathbb C}^\ell\to{\mathbb C}^\ell\left/\Lambda\right.$.  Then
$$
\frac{f^{\left(\sigma^{(k)}\right)}(\zeta)}{f^{\left(\sigma^{(j)}\right)}(\zeta)}=\exp\left(2\pi\sqrt{-1}\sum_{i=1}^\ell\lambda_{j,k,\ell+i}m_i\zeta_i\right).\leqno{({\rm II}.5.3.1)}
$$
This means that the transition function in (II.5.3.1) is of the form $\prod_{i=1}^\ell\left(w_i\right)^{q_i}$, where $w_i=e^{2\pi\sqrt{-1}\,\zeta_i}$ and $q_i\in{\mathbb Z}$ for $1\leq i\leq\ell$.

\medbreak For $1\leq j\leq I$ let $U_j$ be the open subset of $X$ which is the inverse image $\left({\rm alb}_X\right)^{-1}\left(U^{\left(\sigma^{(j)}\right)}\right)$ of $U^{\left(\sigma^{(j)}\right)}$ under the Albanese map ${\rm alb}_X:X\to {\rm Alb}_X$.  The conclusion is that we can use as as the moduli variables $\gamma_j$ (the pullback under the Albanese map ${\rm alb}_X$ of) $\prod_{i=1}^\ell\left(w_i\right)^{q_{i,j}}$ with $w_i=e^{2\pi\sqrt{-1}\,\zeta_i}$ and $q_{i,j}\in{\mathbb Z}$ for $1\leq i\leq\ell$.

\bigbreak\noindent(II.5.4) {\it Algebraic Specialization of Transcendental Numbers to Construct Manifold Defined over Algebraic Number Field.}  Our goal is to show that in (II.2.5.2) every component of the subvariety $Z_q$ in $A$ contains a torsion point of $A$.  Our given manifold $X$ is algebraic and is defined as a complex submanifold of some complex projective space ${\mathbb P}_N$ by a number of homogeneous polynomials $P_j$ ($1\leq j\leq N^\prime$) of the homogeneous coordinates of ${\mathbb P}_N$.  The coefficients of these polynomials belong to a subfield ${\mathbb Q}\left(\xi_1,\cdots,\xi_{n_t}\right)\left[\eta_1,\cdots,\eta_{n_a}\right]$ of ${\mathbb C}$ where $\xi_1,\cdots,\xi_{n_t}$ are algebraically independent over ${\mathbb Q}$ and $\eta_1,\cdots,\eta_{n_a}$ are algebraic over the purely transcendental extension field ${\mathbb Q}\left(\xi_1,\cdots,\xi_{n_t}\right)$ of ${\mathbb Q}$.  If the statement that every component of the subvariety $Z_q$ in $A$ contains a torsion point of $A$ is not true, then we can find elements $\hat\xi_1,\cdots,\hat\xi_{n_t}$ in the algebraic closure $\overline{\mathbb Q}$ of ${\mathbb Q}$ and can construct a manifold $\hat X$ from $X$,
by applying the specialization $\xi_j\mapsto\hat\xi_j$ to replace the coefficients of $P_j$ by elements of $\overline{\mathbb Q}$, such that the statement that every component of the subvariety $Z_q$ in $A$ contains a torsion point of $A$ remains untrue when $X$ is replaced by $\hat X$.  With this technique of algebraic specialization of transcendental numbers, in order to prove that every component of the subvariety $Z_q$ in $A$ contains a torsion point of $A$ it suffices to assume without loss of generality that $X$ is defined over $\overline{\mathbb Q}$ in the sense that $X$ is a complex submanifold of ${\mathbb P}_N$ defined by a number of homogeneous polynomials of the homogeneous coordinates of ${\mathbb P}_N$ whose coefficients are elements of $\overline{\mathbb Q}$.

\bigbreak\noindent(II.5.5) {\it Summary from Two Descriptions of Moduli of Flat Bundles.}  Let us summarize what we conclude from our two ways of describing the moduli for flat line bundles over $X$ (which is now assumed to be defined over $\overline{\mathbb Q}$).   Let $q$ be a positive integer.  Denote by $Z^*_q$ the set of all $a\in A$ such that $\dim_{\mathbb C}\Gamma\left(X,K_X+F^{(a)}\right)$ is no less than $q$.  From the first kind of moduli description we conclude that the image $\Phi\left(Z^*_q\right)$ of $Z^*_q$ under the embedding $\Phi:A\to{\mathbb P}_N$ is defined by a finite number of polynomials of the homogeneous coordinates of ${\mathbb P}_N$ with coefficients in $\overline{\mathbb Q}$.  Assume that $A$ is the quotient of ${\mathbb C}^\ell$ (with complex coordinates $\zeta_1,\cdots,\zeta_\ell$) by the lattice generated by the standard ${\mathbb C}$-basis $u_1,\cdots,u_\ell$ of ${\mathbb C}^\ell$ (as column $\ell$-vectors) and the $\ell$ elements $m_1Zu_1,\cdots,m_\ell Zu_\ell$ with $Z$ being an $\ell\times\ell$ symmetric matrix with positive definite imaginary part and $m_1,\cdots,m_\ell$ being positive integers.  From the second kind of moduli description we conclude that $Z^*_q$ is defined by a finite number of polynomials in the variables $e^{2\pi\sqrt{-1}\,\zeta_j}$ for $1\leq j\leq\ell$ with coefficients in $\overline{\mathbb Q}$.

\bigbreak\noindent(II.5.6) {\it Proposition (Algebraicity with Transition Functions as Variables for Subvariety of Moduli of Flat Bundles Defined by Mimimum Number of Linearly Independent Pluricanonically Twisted Sections).}  Let $X$ be a compact complex algebraic manifold defined over $\overline{\mathbb Q}$.   Let $m$ and $q$ be positive integers.  Denote by $Z_{m,q}$ the set of all $a\in A$ such that $\dim_{\mathbb C}\Gamma\left(X, mK_X+F^{(a)}\right)$ is no less than $q$.  Assume that $A$ is the quotient of ${\mathbb C}^\ell$ (with complex coordinates $\zeta_1,\cdots,\zeta_\ell$) by the lattice generated by the standard ${\mathbb C}$-basis $u_1,\cdots,u_\ell$ of ${\mathbb C}^\ell$ (as column $\ell$-vectors) and the $\ell$ elements $m_1Zu_1,\cdots,m_\ell Zu_\ell$ with $Z$ being an $\ell\times\ell$ symmetric matrix with positive definite imaginary part and $m_1,\cdots,m_\ell$ being positive integers.  Then $Z_{m,q}$ is defined by a finite number of polynomials in the variables $e^{2\pi\sqrt{-1}\,\zeta_j}$ for $1\leq j\leq\ell$ with coefficients in $\overline{\mathbb Q}$.

\medbreak\noindent{\it Proof.}  We let $\hat Z_{m,q}$ be the set of points $a$ of $A$ such that the complex dimension of $\Gamma\left(X, mK_X+mF^{(a)}\right)$ is $\geq q$.  It suffices to prove that $\hat Z_{m,q}$ is defined by a finite number of polynomials in the variables $e^{2\pi\sqrt{-1}\,\zeta_j}$ for $1\leq j\leq\ell$ with coefficients in $\overline{\mathbb Q}$.

\medbreak Consider first the special case when $\hat Z_{m,q}$ is isolated at some point $a_0$ of $A$ and we are going to prove that in a neighborhood of $a_0$ the subvariety $\hat Z_{m,q}$ is defined by a finite number of polynomials in the variables $e^{2\pi\sqrt{-1}\,\zeta_j}$ for $1\leq j\leq\ell$ with coefficients in $\overline{\mathbb Q}$.  Let $$s_1\left(a_0\right),\cdots,s_q\left(a_0\right)\in\Gamma\left(X,mK_X+mF^{\left(a_0\right)}\right)$$ be ${\mathbb C}$-linear independent defined over $\overline{\mathbb Q}$.  For $1\leq\nu\leq q$ take a resolution $\tilde\pi_\nu:\tilde X_\nu\to X$ of the singularity of the graph of $\left(s_\nu\left(a_0\right)\right)^{\frac{1}{m}}$ so that at a generic point of $X$ the fiber of $\tilde\pi_\nu:\tilde X_\nu\to X$ consists of $m$ distinct points and the pullback of $\left(s_\nu\left(a_0\right)\right)^{\frac{1}{m}}$ to $\tilde X_j$ by $\tilde\pi_\nu:\tilde X_\nu\to X$ defines an element of $\tilde s_\nu\left(a_0\right)$ of $\Gamma\left(\tilde X_\nu, K_{\tilde X_\nu}+\left(\tilde\pi_\nu\right)^*\left(F^{\left(a_0\right)}\right)\right)$.   Let $\tilde Z_\nu$ be the set of points $a\in A$ such that the complex dimension of $\Gamma\left(\tilde X_\nu, K_{\tilde X_\nu}+\left(\tilde\pi_\nu\right)^*\left(F^{\left(a\right)}\right)\right)$ is no less than the complex dimension of $\Gamma\left(\tilde X_\nu, K_{\tilde X_\nu}+\left(\tilde\pi_\nu\right)^*\left(F^{\left(a_0\right)}\right)\right)$.  By (II.5.5) each of $\tilde Z_\nu$ is defined by a finite number of polynomials in the variables $e^{2\pi\sqrt{-1}\,\zeta_j}$ for $1\leq j\leq\ell$ with coefficients in $\overline{\mathbb Q}$.

\medbreak We claim that the dimension of the subvariety germ $\cap_{\nu=1}^q\tilde Z_\nu$ is $0$ at $a_0$, which would imply that $\hat Z_{m,q}$ which consists only of the single point $a_0$ in a neighborhood of $a_0$ is defined in a neighborhood of $a_0$ by a finite number of polynomials in the variables $e^{2\pi\sqrt{-1}\,\zeta_j}$ for $1\leq j\leq\ell$ with coefficients in $\overline{\mathbb Q}$.  To prove the claim, we assume the contrary so that there is some positive branch-germ $Z^*$ of $\cap_{\nu=1}^q\tilde Z_\nu$ at $a_0$.  We can find a family of $\tilde s_\nu\left(a\right)\in\Gamma\left(\tilde X_\nu, K_{\tilde X_\nu}+\left(\tilde\pi_\nu\right)^*\left(F^{\left(a\right)}\right)\right)$ for $a\in Z^*$ which is holomorphic in $a$ and which agrees with $\tilde s_\nu\left(a_0\right)$ at $a=a_0$.  We divide the element $\tilde s_\nu\left(a\right)$ of $\Gamma\left(\tilde X_\nu, K_{\tilde X_\nu}+\left(\tilde\pi_\nu\right)^*\left(F^{\left(a\right)}\right)\right)$ by the element
$\tilde s_\nu\left(a_0\right)$ of $\Gamma\left(\tilde X_\nu, K_{\tilde X_\nu}+\left(\tilde\pi_\nu\right)^*\left(F^{\left(a_0\right)}\right)\right)$ to form a meromorphic section $\frac{\tilde s_\nu\left(a\right)}{\tilde s_\nu\left(a_0\right)}$ of the flat line bundle $\left(\tilde\pi_\nu\right)^*\left(F^{\left(a\right)}-F^{\left(a_0\right)}\right)$ on $\tilde X_\nu$.  We now construct the (multiplicatively defined) direct image of $\frac{\tilde s_\nu\left(a\right)}{\tilde s_\nu\left(a_0\right)}$ under $\tilde\pi_\nu$ to get a meromorphic section $s_\nu^\flat\left(a\right)$ of the flat line bundle $mF^{\left(a\right)}-mF^{\left(a_0\right)}$ on $X$ in the sense that the value of $s_\nu^\flat\left(a\right)$ at a generic point $P$ of $X$ is the product of the values of $\frac{\tilde s_\nu\left(a\right)}{\tilde s_\nu\left(a_0\right)}$ at the $m$ distinct points of $\left(\tilde\pi_\nu\right)^{-1}(P)$.  Now we  define the element $\hat s_\nu(a)$ of $\Gamma\left(X,mK_X+mF^{\left(a\right)}\right)$ to be the product $s_\nu\left(a_0\right)s_\nu^\flat(a)$, making points $a$ of $Z^*$ belong to $\hat Z_{m,q}$ for $a$ close to $a_0$, which contradicts the fact that $\hat Z_{m,q}$ is isolated at $a_0$.

\medbreak Now we consider the general case where $\hat Z_{m,q}$ is positive dimensional.  Since by Proposition (II.2.3.1) we know that $\hat Z_{m,q}$ is regular and is a finite union of translates of abelian subvarieties of $A$, we take any abelian subvariety $A^\prime$ of $A$ defined over $\overline{\mathbb Q}$ such that $A^\prime\cap \hat Z_{m,q}$ is $0$-dimensional.  We apply the argument to $A^\prime$ instead of $A$ to conclude that $\hat Z_{m,q}$ is defined by a finite number of polynomials in the variables $e^{2\pi\sqrt{-1}\,\zeta_j}$ for $1\leq j\leq\ell$ with coefficients in $\overline{\mathbb Q}$.  Q.E.D.

\bigbreak

\bigbreak\noindent{\bf\S6. Technique of Gelfond-Schneider, Lang, Bombieri, Brieskorn, and Simpson.}

\bigbreak We now come to the final step of the proof of Theorem (II.1).  This step applies the arithmetic argument of Gelfond-Schneider [Gelfond1934, Schneider1934], Lang [Lang1962, Lang1965, Lang1966], Bombieri [Bombieri1970, Bombieri-Lang1970], Briskorn [Brieskorn1970], and Simpson [Simpson1993] to our two descriptions of the moduli of flat line bundles, one by the use of embedding into projective space and the other by the use of transition functions.

\bigbreak\noindent(II.6.1) {\it Three Ways of Implementing Technique of Gelfond-Schneider.}  In their independently obtained solution of Hilbert's seventh problem [Hilbert1900] of the transcendence of $a^b$ with $a$ algebraic $\not=0,1$ and $b$ algebraic irrational, Gelfond and Schneider [Gelfond1934, Schneider1934] introduced a method of applying Nevanlinna's First Main Theorem [Nevanlinna1925] to a sequence of polynomials, with algebraic coefficients, of algebraically independent entire functions of finite-order growth and passing to limit.  This method originally introduced by Gelfond and Schneider [Gelfond1934, Schneider1934] for entire functions on ${\mathbb C}$ was later generalized by Lang [Lang1962, Lang1965, Lang1966] and Bombieri [Bombieri1970, Bombieri-Lang1970] to meromorphic functions on ${\mathbb C}^d$ to yield the following result of Bombieri [Bombieri1970, p.267, Theorem A].

\bigbreak\noindent(II.6.1.1) If $Z$ is the set of points of ${\mathbb C}^d$ where more than $d$ algebraically independent meromorphic functions $f_j$ on ${\mathbb C}^d$ with growth of finite order $\rho$ assume values in a number field $K$, then $Z$ is contained in an algebraic hypersurface whose degree bounded by a number in terms of $d$ and $\rho$ and the degree of $K$ over $\mathbb Q$ if every first-order partial derivative of each $f_j$ can be expressed as a rational function of all $f_j$ with coefficients in $K$.

\medbreak This is the original, first way of implementing the technique of Gelfond-Schneider.  For the purpose of our note this first way of implementation cannot be applied to our situation for the following reason.  No complex hypersurface in ${\mathbb C}^d$ can contain an additive subgroup of ${\mathbb C}^d$ generated by $d$ ${\mathbb C}$-linearly independent elements of ${\mathbb C}^d$ over ${\mathbb Z}$.  In our case we consider points of an abelian variety which is the quotient of ${\mathbb C}^d$ by a lattice $\Lambda$ of rank $2d$ and we want to conclude that some nonempty set $E$ of points in $A$ must be a finite subset of $A$ because some meromorphic functions assume algebraic values on the inverse image $\tilde E$ of $E$ in the universal cover of $A$.  With $\tilde E$ invariant under translations by elements
of $\Lambda$, this implementation of the technique of Gelfond-Schneider could only give the conclusion that the subset $E$ in $A$ is empty and cannot give the conclusion that $E$ is a finite subset of $A$.

\medbreak What we need for our situation is actually a later modification of the technique of Gelfond-Schneider by Bombieri-Lang [Bombieri-Lang1970].  This second way of implementing the technique of Gelfond-Schneider does not use the condition of expressing every first-order partial derivative of each finite-order transcendental function as a rational function of all the given transcendental functions with coefficients in the number field.  The conclusion is of a different nature.  It gives a statement limiting the growth order of density of the set of points where all the given transcendental functions assume algebraic values with heights growing with some prescribed order (see [Bombieri-Lang1970, p.8, Th.1; p.11, Th.2; p.11, Th.3]).  The differential equation in the first way of implementing the technique of Gelfond-Schneider is to make it possible to construct a polynomial of the transcendental functions with coefficients in the number field with height estimates so that the polynomial vanishes to high order at points of $Z$.  The high-order vanishing is needed to make $Z$ contained in an algebraic hypersurface in ${\mathbb C}^d$.  In the second way of implementing the technique of Gelfond-Schneider no high-vanishing order is used so that the differential equation is not needed but the conclusion is weaker and concerns only the density of the set $Z$ instead of $Z$ being contained in an algebraic hypersurface.

\medbreak In conjunction with the second way of implementing the technique of Gelfond-Schneider, we would like to remark that there are some more recent developments of it in a very general setting by Bombieri-Pila [Bombieri-Pila1989] and Pila-Wilkie [Pila-Wilkie2006].

\medbreak The third implementation of the technique of Gelfond-Schneider was introduced by Brieskorn in [Brieskorn1970] to conclude that the characteristic polynomial of Milnor's local complex Picard-Lefschetz monodromy for the cohomology group with integer coefficients in the case of a complex hypersurface of isolated singularity [Milnor1968] is a product of cyclotomic polynomials (see [Brieskorn1970, p.11, Satz 4]).  He introduced the additional technique of (non-continuous) isomorphisms of ${\mathbb C}$.  In his argument the root of his characteristic polynomial is of the form $e^{2\pi\sqrt{-1}\,\mu_j}$, where $\mu_j$ is algebraically determined by the coefficients of the complex hypersurface.  Such roots $e^{2\pi\sqrt{-1}\,\mu_j}$ are always algebraic.  If $\mu_j$ is not algebraic, it is possible to use a (non-continuous) automorphism $\phi$ of ${\mathbb C}$ such that $e^{2\pi\sqrt{-1}\,\phi\left(\mu_j\right)}$ is not algebraic, giving a contradiction.  When $\mu_j$ is known to be algebraic, $\mu_j$ must be rational otherwise the original result of Gelfond-Schneider implies that $e^{2\pi\sqrt{-1}\,\mu_j}$ is not algebraic.  This third implementation of the technique of Gelfond-Schneider is also not applicable to our situation.

\bigbreak\noindent(II.6.2) {\it Application of Technique of Gelfond-Schneider to Locate Torsion Flat Line Bundles.}  We now go back to our compact complex algebraic manifold $X$ which we can assume to be defined over the algebraic closure $\overline{\mathbb Q}$ of the field ${\mathbb Q}$ of all rational numbers, as explained in (II.5.4).  Let $A$ be the abelian variety of all flat line bundles on $X$ whose universal cover is $\tilde A$.

\medbreak We compose the universal covering map ${\mathbb C}^\ell=\tilde A\to A$ with the embedding $\Phi:A\to{\mathbb P}_N$ defined over $\overline{\mathbb Q}$ to get entire functions $\hat f_0,\cdots,\hat f_N$ on ${\mathbb C}^\ell$ which can be used as homogeneous components for the composite map ${\mathbb C}^\ell\to A\stackrel{\Phi}{\longrightarrow}{\mathbb P}_N$.  We use the quotients $f_j=\frac{\hat f_j}{\hat f_0}$ for $1\leq j\leq N$ as some of the transcendental meromorphic functions on ${\mathbb C}^\ell$ for the application of the technique of Gelfond-Schneider (in its second way of implementation as explained in (II.6.1)).

\medbreak The differentials $d\zeta_j$ of the coordinates $\zeta_j$ of ${\mathbb C}^\ell$ define holomorphic $1$-forms on $A$ which via $\Phi$ correspond to holomorphic $1$-forms $\omega_1,\cdots,\omega_\ell$ on $\Phi\left(A\right)$ expressed in terms of $f_1,\cdots,f_N,df_1,\cdots,df_N$ defined over $\overline{\mathbb Q}$.  The differential equations expressing every first-order derivative of each $f_j$ as a rational functions of $f_1,\cdots,f_N$ just come from $\frac{\Phi^*\left(\omega_j\right)}{d\zeta_k}$ being equal to the Kronecker delta $\delta_{jk}$ on ${\mathbb C}^\ell$ and from differentiating the defining functions of $\Phi(A)$ in ${\mathbb P}_N$ over $\overline{\mathbb Q}$.

\medbreak Because the complex dimension of $\Phi(A)$ is $\ell$, the maximum number of algebraically independent functions in the set $f_1,\cdots,f_N$ cannot be more than $\ell$.  In order to have more algebraically independent functions than the complex dimension $\ell$ of ${\mathbb C}^\ell$, we add the $\ell$ functions $e^{2\pi\sqrt{-1}\,\zeta_1},\cdots,e^{2\pi\sqrt{-1}\,\zeta_\ell}$ which come from the second description of moduli of flat line bundles.  Since $\frac{\partial}{\partial\zeta_j}\,e^{2\pi\sqrt{-1}\,\zeta_j}=2\pi\sqrt{-1}$ and $2\pi\sqrt{-1}$ is not an element of $\overline{\mathbb Q}$, the differential equation condition is not satisfied by the functions
$e^{2\pi\sqrt{-1}\,\zeta_1},\cdots,e^{2\pi\sqrt{-1}\,\zeta_\ell}$.  Moreover, in the discussion of the preceding paragraph, we know that we cannot hope to use the differential equation in applying the technique of Gelfond-Schneider to our situation.   Though in [Simpson1993,p.369, Proof of Proposition 3.4] Simpson used the differential equation which comes from the group composition law and the fact that the differential at the group identity is defined over $\overline{\mathbb Q}$, at least in our simple explicit analytic setting the second implementation of the technique of Gelfond-Schneider using differential equations cannot be used.

\medbreak For the application of the technique of Gelfond-Schneider (in its second implementation as explained in (II.6.1)) we now come to the question of the set $Z$ of points of $A$ where the functions $$f_1,\cdots,f_N,e^{2\pi\sqrt{-1}\,\zeta_1},\cdots,e^{2\pi\sqrt{-1}\,\zeta_\ell}$$ assume algebraic values.  We do not yet have such a set $Z$, but for a positive integer $q$ we have a subvariety $Z_q$ defined by the set of all $a\in A$ such that $\dim_{\mathbb C}\Gamma\left(X, K_X+F^{(a)}\right)\geq q$.  This subvariety is regular and is a finite union of translates of abelian subvarieties of $A$ and can be defined by equations in two different ways.  One way is that $Z_q$ the common zero-set of a finite number of polynomials of $f_1,\cdots,f_N$ with coefficients in $\overline{\mathbb Q}$.  Another way is that $Z_q$ is the common zero-set of a finite number of polynomials of $e^{2\pi\sqrt{-1}\,\zeta_1},\cdots,e^{2\pi\sqrt{-1}\,\zeta_\ell}$ with coefficients in $\overline{\mathbb Q}$.  If $Z_q$ is zero-dimensional, these two ways of defining $Z_q$ precisely imply that the functions $$f_1,\cdots,f_N,e^{2\pi\sqrt{-1}\,\zeta_1},\cdots,e^{2\pi\sqrt{-1}\,\zeta_\ell}$$ assume algebraic values at points of $Z_q$.

\medbreak In general, the complex dimension of $Z_q$ is positive.  The method is to reduce the general case to the case of zero dimension by taking the quotient of $A$ by the abelian subvariety of which a component $Z^\prime_q$ of $Z_q$ is a translate.  Write $Z^\prime_q=a^\prime+A^{\prime\prime}$, where $a^\prime\in A$ and $A^{\prime\prime}$ is an abelian subvariety of $A$.  Since $\Phi\left(Z^\prime_q\right)$ is a subvariety of ${\mathbb P}_N$ defined over $\overline{\mathbb Q}$, we can find an abelian subvariety $A^\prime$ of $A$ such that $\Phi\left(A^\prime\right)$ is a subvariety of ${\mathbb P}_N$ defined over $\overline{\mathbb Q}$ and $A^\prime$ intersects $Z^\prime_q$ at a single point $\hat a$ of $A$.  We now replace the abelian variety $A$ by the abelian variety $A^\prime$ (or equivalently consider the quotient abelian variety $A\left/A^{\prime\prime}\right.$).  Restrict the functions $$f_1,\cdots,f_N,e^{2\pi\sqrt{-1}\,\zeta_1},\cdots,e^{2\pi\sqrt{-1}\,\zeta_\ell}$$ to $A^\prime$ and then pull them back to the universal cover $\tilde A^\prime={\mathbb C}^{\ell^\prime}$ of $A^\prime$ under the universal covering map $\pi^\prime:{\mathbb C}^{\ell^\prime}\to A^\prime$ to get functions $F_1,\cdots,F_N,F_{N+1},\cdots,F_{N+\ell}$ on ${\mathbb C}^{\ell^\prime}$.  Let $E$ be the subgroup of $A^\prime$ generated by $a^\prime$.  Since
the values of the functions $$f_1,\cdots,f_N,e^{2\pi\sqrt{-1}\,\zeta_1},\cdots,e^{2\pi\sqrt{-1}\,\zeta_\ell}$$ at $a^\prime$ belong to $\overline{\mathbb Q}$, it follows from the addition formula for the abelian variety $A^\prime$ and the exponential law for the exponential function $e^{2\pi\sqrt{-1}\,\zeta_1},\cdots,e^{2\pi\sqrt{-1}\,\zeta_\ell}$ that the values of the functions $$f_1,\cdots,f_N,e^{2\pi\sqrt{-1}\,\zeta_1},\cdots,e^{2\pi\sqrt{-1}\,\zeta_\ell}$$ at $a^\prime$ belong to $\overline{\mathbb Q}$ at every point of $E$ belong also to $\overline{\mathbb Q}$.  As a result the values of the functions $F_1,\cdots,F_N,F_{N+1},\cdots,F_{N+\ell}$ at every point of $\left(\pi^\prime\right)^{-1}\left(E\right)$ belong to $\overline{\mathbb Q}$.  If $E$ is not finite, by the arguments of [Bombieri-Lang1970] the density of the subset $\left(\pi^\prime\right)^{-1}\left(E\right)$ of ${\mathbb C}^{\ell^\prime}$ has big enough growth order to give us a contradiction.  Hence we conclude that $E$ is a finite and the point $a^\prime$ of $Z_q$ is a torsion point for $A$.  This finally concludes the proof of Theorem (II.1).

\bigbreak

\bigbreak\noindent{\bf References}

\bigbreak\noindent[Andreotti-Vesentini1965]
Aldo Andreotti and Edoardo Vesentini,
Carleman estimates for the Laplace-Beltrami equation on complex manifolds.
{\it Inst. Hautes \'Etudes Sci. Publ. Math.} \textbf{25} (1965) 81–-130.

\medbreak\noindent[Angehrn-Siu1995]
Urban Angehrn and Yum-Tong Siu,
Effective freeness and point separation for adjoint bundles.
{\it Invent. Math.} \textbf{122} (1995), 291–-308.

\medbreak\noindent{Artin1968]
Michael Artin, On the Solutions of Analytic Equations, {\it Invent. Math.} \textbf{5} (1968), 277--291.

\medbreak\noindent{Artin1968]
Michael Artin, Algebraic approximation of structures over complete local rings. {\it Pub. Math. I.H.E.S.}
\textbf{36} (1969), 23--58.

\medbreak\noindent[Berndtsson1996] Bo Berndtsson, The extension theorem of Ohsawa-Takegoshi and the theorem of Donnelly-Fefferman, {\it Ann. Inst. Fourier} (Grenoble) \textbf{46} (1996), 1083 -- 1094.

\medbreak\noindent[Bombieri1970] Enrico Bombieri, Algebraic values of meromorphic maps, {\it Invent.
Math.} \textbf{10} (1970), 267-287.  Addendum. {\it Invent.
Math.} \textbf{11} (1970), 163--166.

\medbreak\noindent[Bombieri-Lang1970] Enrico Bombieri and Serge Lang,
Analytic subgroups of group varieties.
{\it Invent. Math.} \textbf{11} (1970), 1--14.

\medbreak\noindent[Bombieri-Pila1989] Enrico Bombieri and Jonathan Pila, The number of integral points on arcs and ovals, {\it Duke Math. J.} \textbf{59} (1989), 337--357.

\medbreak\noindent[Brieskorn1970]
Egbert Brieskorn,
Die Monodromie der isolierten Singularit\"aten von Hyperfl\"achen.
{\it Manuscripta Math.} \textbf{2} (1970), 103--161.

\medbreak\noindent
[Budur2009] Nero Budur, Unitary local systems, multiplier ideals, and polynomial
periodicity of Hodge numbers. {\it Adv. Math.} \textbf{221} (2009), 217--250.

\medbreak\noindent[Campana-Peternell-Toma2007]
Frederic Campana, Thomas Peternell, and Matei Toma
Geometric stability of the cotangent bundle and the universal cover of a projective manifold. arXiv:math/0405093.

\medbreak\noindent[Demailly1985]
Jean-Pierre Demailly,
Champs magn\'etiques et in\'egalit\'es de Morse pour la $d''$-cohomologie.
{\it Ann. Inst. Fourier} (Grenoble) \textbf{35} (1985), 189–-229.

\medbreak\noindent[Donnelly-Fefferman1983]
Harold Donnelly and Charles Fefferman,
$L^2$-cohomology and index theorem for the Bergman metric.
{\it Ann. of Math.} \textbf{118} (1983), 593–-618.

\medbreak\noindent[Donnelly-Xavier1984]
Harold Donnelly and Frederico Xavier,
On the differential form spectrum of negatively curved Riemannian manifolds.
{\it Amer. J. Math.} \textbf{106} (1984), 169–-185.

\medbreak\noindent[Gelfond1934] A. O. Gelfond, Sur le septi\`eme Probl\`eme de D. Hilbert. {\it Comptes Rendus Acad. Sci. URSS Moscou} \textbf{2} (1934), 1-6. {\it Bull. Acad. Sci. URSS Leningrade} \textbf{7} (1934), 623--634.

\medbreak\noindent[Hensel1897] Kurt Hensel, \"Uber eine neue Begr\"undung der Theorie der
algebraischen Zahlen. {\it Jahresbericht der Deutschen
Mathematiker-Vereinigung} \textbf{6} (1897), 83 -- 88.

\medbreak\noindent[Hilbert1900] David Hilbert, Mathematische Probleme.  {\it Nachr. K\"onigl. Ges. der Wiss. zu G\"ottingen, Math.-Phys. Klasse} (1900), 251--297.

\medbreak\noindent[H\"ormander1965]
Lars H\"ormander,
$L^{2}$ estimates and existence theorems for the $\bar\partial$ operator.
{\it Acta Math.} \textbf{113} (1965), 89–-152.

\medbreak\noindent[Kodaira1954]
Kunihiko Kodaira,
On K\"ahler varieties of restricted type (an intrinsic characterization of algebraic varieties).
{\it Ann. of Math.} \textbf{60} (1954), 28–-48.

\medbreak\noindent[Kim2010] Dano Kim,
$L^2$ extension of adjoint line bundle sections.
{\it Ann. Inst. Fourier (Grenoble)} \textbf{60} (2010), 1435–-1477.

\medbreak\noindent[Kohn1963-64] Joseph J. Kohn,
Harmonic integrals on strongly pseudo-convex manifolds. I.
{\it Ann. of Math.} \textbf{78} (1963), 112–-148; II.
{\it Ann. of Math.} \textbf{79} (1964), 450–-472.

\medbreak\noindent[Lang1962] Serge Lang,
Transcendental points on group varieties.
{\it Topology} \textbf{1} (1962), 313--318.

\medbreak\noindent[Lang1965] Serge Lang,
Algebraic values of meromorphic functions.
{\it Topology} \textbf{3} (1965), 183--191.

\medbreak\noindent[Lang1966] Serge Lang,
{\it Introduction to transcendental numbers}. Addison-Wesley Publishing Co., Reading, Mass.-London-Don Mills, Ont. 1966.

\medbreak\noindent[Manivel1993]
Laurent Manivel,
Un th\'eor\`eme de prolongement $L^2$ de sections holomorphes d'un fibr\'e hermitien.
{\it Math. Z.} \textbf{212} (1993), 107–-122.

\medbreak\noindent[Milnor1968]
John Milnor,
Singular points of complex hypersurfaces.
{\it Annals of Mathematics Studies}, No. \textbf{61} Princeton University Press, Princeton, N.J.; University of Tokyo Press, Tokyo 1968.

\medbreak\noindent[Morrey1958]
Charles B. Morrey Jr.,
The analytic embedding of abstract real-analytic manifolds.
{\it Ann. of Math.}  \textbf{68} (1958), 159–-201.

\medbreak\noindent[Nevanlinna1925]
Rolf Nevanlinna,
Zur Theorie der Meromorphen Funktionen.
{\it Acta Math.}
\textbf{46} (1925), 1--99.

\medbreak\noindent[Ohsawa-Takegoshi1987]
Takeo Ohsawa and Kensho Takegoshi,
On the extension of $L^2$ holomorphic functions.
{\it Math. Zeitschr.} \textbf{195} (1987),  197--204.

\medbreak\noindent[Paun2007]
Mihai Paun,
Siu's invariance of plurigenera: a one-tower proof.
{\it J. Differential Geom.} \textbf{76} (2007), 485--493.

\medbreak\noindent[Pila-Wilkie2006] Jonathan Pila and Alex J. Wilkie, The rational points of a definable set, {\it Duke Math. J.} \textbf{33} (2006), 591--616.

\medbreak\noindent[Popvici2005]
Dan Popovici,
$L^2$ extension for jets of holomorphic sections of a Hermitian line bundle.
{\it Nagoya Math. J.} \textbf{180} (2005), 1–-34.

\medbreak\noindent[Schneider1934] T. Schneider, Transzendenzuntersuchungen periodischer Funktionen. I, II. {\it J. reine angew. Math.} \textbf{172} (1934), 65--74.

\medbreak\noindent[Simpson1993] Carlos Simpson,
Subspaces of moduli spaces of rank one local systems.
{\it Ann. Sci. \'Ecole Norm. Sup.} \textbf{26} (1993), 361--401.

\medbreak\noindent[Siu1996] Yum-Tong Siu, The Fujita conjecture and the extension theorem of Ohsawa-Takegoshi. {\it Geometric complex analysis (Hayama, 1995)}, 577–592, World Sci. Publ., River Edge, NJ, 1996.

\medbreak\noindent[Siu1998] Yum-Tong Siu, Invariance of plurigenera.
{\it Invent. Math.} \textbf{134} (1998), 661--673.

\medbreak\noindent[Siu2002] Yum-Tong Siu,
Extension of twisted pluricanonical sections with plurisubharmonic weight and invariance of semipositively twisted plurigenera for manifolds not necessarily of general type. {\it Complex geometry} (G\"ottingen, 2000), 223--277, Springer, Berlin, 2002.

\medbreak\noindent[Siu2010] Yum-Tong Siu,
Abundance Conjecture in {\it Geometry and Analysis, Vol II.} ed. Lizhen Ji, International Press 2010, pp.271--317.
(arXiv:math/0912.0576)

\medbreak\noindent[Takayama2006]
Shigeharu Takayama,
Pluricanonical systems on algebraic varieties of general type.
{\it Invent. Math.} \textbf{165} (2006), 551-–587.

\medbreak\noindent[Varolin2008]
Dror Varolin,
A Takayama-type extension theorem.
{\it Compos. Math.} \textbf{144} (2008), 522--540.

\medbreak\noindent[Wavrik1975]
A Theorem on Solutions of Analytic Equations
with Applications to Deformations of Complex Structures, {\it Math. Ann.} \textbf{216} (1975), 127--142.

\medbreak\noindent[Weil1958]
Andr\'e Weil,
{\it Introduction \`a l'\'etude des vari\'et\'es k\"ahl\'eriennes}.
Publications de l'Institut de Math\'ematique de l'Universit\'e de Nancago, VI. Actualit\'es Sci. Ind. no.1267, Hermann, Paris 1958.

\bigbreak\noindent{\it Author's mailing address}: Department of
Mathematics, Harvard University, Cambridge, MA 02138, U.S.A.

\medbreak\noindent {\it Author's e-mail address}:
siu@math.harvard.edu

\end{document}